\documentclass[12pt,fleqn]{article}
\usepackage[cp866]{inputenc}
\usepackage[english]{babel}
\usepackage{makeidx}
\usepackage{latexsym,amsfonts,amssymb,amsmath,longtable,amsthm}
\input amssym.def
\usepackage{latexsym}
\usepackage{graphicx}
\usepackage[all]{xy}
\usepackage{amsfonts}
\usepackage{amsmath}
\usepackage{mathrsfs}
\usepackage{color}
\usepackage{ulem}
\usepackage[table]{xcolor}

\def\e{\varepsilon}

\def\0{{\mathbf{0}}}
\def\a{{\mathbf{a}}}

\def\w{{\mathbf{w}}}
\def\n{{\bf n}}
\def\ov{\boldsymbol{\omega}}
\def\u{{\mathbf{u}}}

\def\h{{\mathbf{A}}}
\def\s{{\boldsymbol{\sigma}}}
\def\D{{\mathcal{D}}}

\def\Ti{{\mathscr{T}}}

\def\Si{{\mathscr{S}}}

\def\Ha{{\mathfrak{H}}}

\def\F{{\mathcal{F}}}

\def\N{\mathbb{N}}

\def\F{{\mathcal{F}}}

\def\bfeta{\boldsymbol{\eta}}

\def\int{\intop}

\newcommand{\wotwo}{\mathop{\mathaccent23{W}^{1,2}}\nolimits}
\newcommand{\wotwotwo}{\mathop{\mathaccent23{W}^{2,2}}\nolimits}

\newtheorem{theo}{\bf Theorem}
\newtheorem{lem}{\bf Lemma}
\newtheorem{lemA}{\bf Corollary}
\newtheorem{lemr}{\bf Remark}

\newcommand{\pr}{{\bf Proof. }}

\def\Xint#1{\mathchoice
   {\XXint\displaystyle\textstyle{#1}}%
   {\XXint\textstyle\scriptstyle{#1}}%
   {\XXint\scriptstyle\scriptscriptstyle{#1}}%
   {\XXint\scriptscriptstyle\scriptscriptstyle{#1}}%
   \!\int}
\def\XXint#1#2#3{{\setbox0=\hbox{$#1{#2#3}{\int}$}
     \vcenter{\hbox{$#2#3$}}\kern-.5\wd0}}

\def\dashint{\Xint-}

\newcommand{\esssup}{\mathop{\rm ess\,sup}}

\newcommand{\ue}{{\mathbf u}}
\newcommand{\fe}{{\mathbf f}}

\newcommand{\R}{{\mathbb R}}
\newcommand{\const}{{\rm const}}
\newcommand{\dist}{\mathop{\rm dist}}
\newcommand{\arctg}{\mathop{\rm arctg}}
\newcommand{\Int}{\mathop{\rm Int}}

\newcommand{\meas}{\mathop{\rm meas}}

\newcommand{\ve}{{\mathbf v}}
\newcommand{\diam}{\mathop{\rm diam}}
\newcommand{\curl}{\mathop{\rm curl}}
\renewcommand{\div}{\mathop{\rm div}}
\renewcommand{\Pr}{\mathop{\rm Proj}\nolimits}
\renewcommand{\H}{{\mathcal H}}
\newcommand{\loc}{{\rm loc}}

\begin{document}
\title{The existence theorem   for the steady
Navier--Stokes  problem  in exterior axially symmetric 3D
domains \footnote{{\it Mathematical Subject classification\/}
(2000). 35Q30, 76D03, 76D05; {\it Key words}: three dimensional
exterior axially symmetric domains, stationary Navier Stokes
equations, boundary--value problem.}}
\author{Mikhail Korobkov\footnote{Sobolev Institute of Mathematics,
Koptyuga pr. 4, and Novosibirsk State University, Pirogova Str. 2,
630090 Novosibirsk, Russia; korob@math.nsc.ru}, Konstantin
Pileckas\footnote{Faculty of Mathematics and Informatics, Vilnius
University, Naugarduko Str., 24, Vilnius, 03225  Lithuania;
pileckas@ktl.mii.lt } and Remigio Russo\footnote{  Seconda
Universit\`a di Napoli,  via Vivaldi 43, 81100 Caserta, Italy;
remigio.russo@unina2.it}$\;$}

\maketitle

\begin{abstract} {We study the nonhomogeneous boundary value problem
for the Navier--Stokes equations of steady motion of a viscous
incompressible fluid in  a three--dimensional   exterior domain
with multiply connected boundary. We prove that this problem has a
solution for axially symmetric  domains and data (without any smallness restrictions on
the fluxes).  }
\end{abstract}

\section{Introduction}
In this paper we shall consider the Navier--Stokes problem
\begin{equation}
\label{NS}
 \left\{\begin{array}{rcl}
-\nu \Delta{\bf u}+\big({\bf u}\cdot \nabla\big){\bf u} +\nabla p
&
= & {\bf f}\qquad \hbox{\rm in }\;\;\Omega,\\[4pt]
\div\,{\bf u} &  = & 0  \qquad \hbox{\rm in }\;\;\Omega,\\[4pt]
 {\bf u} &  = & {\bf a} \qquad \hbox{\rm on }\;\;\partial\Omega,
 \\[4pt]
\lim\limits_{|x|\to+\infty}{\bf u}(x) &  = &{\bf u}_0
\end{array}\right.
\end{equation}
in the exterior domain of
 $\R^3$
\begin{equation}\label{Omega}
\begin{array}{rcl}
\Omega=\R^3\setminus\bigl(\bigcup\limits_{j=1}^N\bar\Omega_j\bigr),
\end{array}
\end{equation}
where $\Omega_i$ are bounded domains with connected $C^2$-smooth
boundaries $\Gamma_i$ and $\bar\Omega_j\cap\bar\Omega_i=\emptyset$
for $i\ne j$.  In (\ref{NS}) $\nu
>0$ is the viscosity coefficient,
${\bf u}$, $p$ are the (unknown) velocity and pressure fields,
${\bf a}$ and ${\bf u}_0$ are the (assigned) boundary data and a
constant vector respectively, $\fe$ is the body force density.

Let
\begin{equation}
\label{FluxN} {\cal F}_i= \int_{\Gamma_i}{\bf a}\cdot{\bf n}
\,dS,\quad i=1,\ldots N,
\end{equation}
 where ${\bf n}$ is the unit outward normal to $\partial\Omega$.
 Under suitable regularity hypotheses on $\Omega$ and ${\bf a}$
 and  assuming that
\begin{equation} \label{FluxNL} {\cal F}_i=0,\quad
i=1,\ldots N,
\end{equation}
 in the celebrated paper~\cite{Leray} of 1933,
 J. Leray was able to show that (\ref{NS})  has a  solution ${\bf u}$  such
 that
\begin{equation}
\label{DI}
 \int_{ \Omega}|\nabla{\bf u}|^2\,dx<+\infty,
\end{equation}
and ${\bf u}$ satisfies (\ref{NS}${}_4$) in a suitable sense for general ${\bf
u}_0$ and uniformly for ${\bf u}_0={\bf 0}$. In the fifties the
problem was reconsidered by  R. Finn \cite{Finn} and O.A.
Ladyzhenskaia \cite{Lad1}, \cite{Lad}. They showed that the solution
satisfies the condition at infinity uniformly. Moreover, the
condition~(\ref{FluxNL}) and the regularity of ${\bf a}$ have been
relaxed by requiring $\sum\limits_{i=1}^N|{\cal F}_i|$ to be sufficiently
small \cite{Finn} and    ${\bf a}\in
W^{1/2,2}(\partial\Omega)$~\cite{Lad}.

In 1973 K.I. Babenko \cite{Babenko} proved that if $({\bf u},p)$
is a  solution to (\ref{NS}), (\ref{DI}) with ${\bf u}_0\ne{\bf
0}$, then $({\bf u}-{\bf u}_0,p)$ behaves at infinity as the
solutions to the linear Oseen system. In particular\footnote{See
also \cite{Galdibook}. Here the symbol $f(x)=O(g(r))$  means that
there is a positive constant $c$ such that $|f(x)|\le c g(r)$ for
large $r$},
 \begin{equation}
 \label{CAB}
 {\bf u}(x)-{\bf u}_0=O(r^{-1}),\quad p(x)=O(r^{-2}).
 \end{equation}
 However, nothing is known, in general, on the rate of convergence at
infinity for ${\bf u}_0={\bf 0}$ \footnote{For small $\|{\bf
a}\|_{L^\infty(\partial\Omega)}$ existence of a solution $({\bf
u}, p)$ to (\ref{NS})  such that ${\bf u}=O(r^{-1})$  is a simple
consequence of Banach contractions theorem \cite{RuSta}. Moreover,
one can show that   $p=O(r^{-2})$ and the derivatives of order $k$
of ${\bf u}$  and $p$ behave at infinity as $r^{-k-1}$,
$r^{-k-2}$, respectively \cite{ST} (see also \cite{Galdibook},
\cite{NP}).}.

 One of the most important problems in the theory of the
steady--state Navier--Stokes equations  concerns  the possibility
to prove existence  of a solution to (\ref{NS}) without any
assumptions on the fluxes ${\cal F}_i$ (see, e.g.
\cite{Galdibook}). To the
best of our knowledge, the most general assumptions assuring
existence is expressed by
 \begin{equation}
 \label{IpoF}
 \sum_{i=1}^N \max_{\Gamma_i}{|{\cal F}_i|\over|x-x_i|}<8\pi\nu
 \end{equation}
(see \cite{RussoPG}), where ${\cal F}_i$ is defined by (\ref{FluxN}) and $x_i$ is a
fixed point of $\Omega_i$ (see also \cite{BOPI} for analogous
conditions in bounded domains).

The present paper is  devoted to the above question in  the
axially symmetric case. To introduce the problem we have to specify
some notations.
 Let
$O_{x_1},O_{x_2},O_{x_3}$ be coordinate axis in $\R^3$  and
$\theta=\arctg(x_2/x_1)$, $r=(x_1^2+x_2^2)^{1/2}$, $z=x_3$ be
cylindrical coordinates. Denote by $v_\theta,v_r,v_z$ the
projections of the vector ${\bf v}$ on the axes $\theta,r,z$.

A function $f$ is said to be {\it axially symmetric} if it does
not depend on~$\theta$. A vector-valued function ${\bf
h}=(h_\theta,h_r,h_z)$ is called {\it axially symmetric} if
$h_\theta$, $h_r$ and $h_z$ do not depend on~$\theta$. A
vector-valued function ${\bf h}=(h_\theta,h_r,h_z)$ is called
{\it axially symmetric with no swirl} if $h_\theta=0$, while $h_r$
and $h_z$ do not depend on~$\theta$.

Note that for axially-symmetric solutions~$\ue$ of~(\ref{NS}) the
vector~$\ue_0$ has to be parallel to the symmetry axis. The main
result of the paper is the following

\begin{theo}
 \label{kmpTh4.2} {\sl
Assume that  $\Omega\subset\R^3$ is an exterior  axially symmetric
domain~\eqref{Omega} with $C^2$-smooth boundary
$\partial\Omega$,  $\ue_0\in\R^3$ is a constant vector parallel to
the symmetry axis, and $\fe\in W^{1,2}(\Omega)\cap L^{6/5}(\Omega)$,
${\bf a}\in
W^{3/2,2}(\partial\Omega)$ are axially symmetric. Then
$(\ref{NS})$~admits at least one weak axially symmetric
solution~${\bf u}$ satisfying~$(\ref{DI})$. Moreover, if $\bf a$
and  $\fe$ are axially symmetric with no swirl, then~$(\ref{NS})$
admits at least one weak axially
symmetric solution satisfying~$(\ref{DI})$ with no swirl. }
\end{theo}

\begin{lemr}
\label{Russo2}{\rm It is well known (see, e.g.,~\cite{Lad}) that
under hypothesis of Theorem~\ref{kmpTh4.2}, every weak solution
$\ue$ of problem~(\ref{NS}) is more regular, i.e, $\ue\in
W_\loc^{2,2}(\overline\Omega)\cap W^{3,2}_\loc(\Omega)$. }
\end{lemr}

Let us emphasize that Theorem \ref{kmpTh4.2} furnishes the first existence
result without any assumption on the fluxes for the stationary Navier--Stokes problem in an
exterior domains.

Note that in the papers~\cite{kpr_a_arx}--\cite{kpr_a_ann}
existence of a solution to
problem~(\ref{NS}${}_1$)--(\ref{NS}${}_3$) in  arbitrary
$C^2$-smooth  {\bf bounded} plane or axially symmetric spatial
domain~$\Omega$ has been  proved under the sole  condition of
zero total flux through the boundary (for a historical review
in the case of bounded domains see, e.g., \cite{kpr_a_arx}, \cite{Pukhnachev}--\cite{Pukhnachev1}).

The proof of the existence theorem is based on  an a priori
estimate which we derive using the classical \textit{reductio ad
absurdum} argument of J.~Leray~\cite{Leray}, see Section~\ref{poet}.
As   well-known, after applying  Leray's argument one comes along to a solution of Euler system satisfying
zero boundary conditions.
Such solutions are studied in Section~\ref{eueq}. The essentially new part
here is the use of Bernoulli's law obtained in~\cite{korob1} for
Sobolev solutions to the Euler equations (the detailed proofs are
presented in~\cite{kpr}  and~\cite{kpr_a_arx}  for  the plane  and
the axially symmetric bounded domains, respectively). In
Section~\ref{eueq} we present the proof of the Bernoulli Law for unbounded domains
(Theorem~\ref{kmpTh2.2}, see also Lemma~\ref{lemTh2.2unb} and
Remark~\ref{rem__int}). Furthermore, we prove here that the value
of the pressure on the boundary components intersecting the
symmetry axis coincides with the value of the pressure at infinity
(see Corollary~\ref{un_cor1}). This phenomena is connected with
the fact that  the symmetry axis can be approximated by stream
lines, where the total head pressure is constant (see
Theorem~\ref{kmpTh2.3_un}).

 The results
concerning Bernoulli's law are based on the recent version of
the~Morse-Sard theorem proved by J. Bourgain, M.~Korobkov and J.
Kristensen \cite{korob}, see also subsection~\ref{SMS}. This
theorem implies, in particular, that almost all level sets of a
function $\psi\in W^{2,1}(\R^2)$ are finite unions of $C^1$-curves
homeomorphic to a~circle.

We obtain the required contradiction in Section~\ref{contrad_s}
for the case~$\ue_0=0$. The above mentioned results allow to
construct suitable subdomains (bounded by smooth stream lines) and
to estimate the $L^2$-norm of the gradient of the total head
pressure. We use here some ideas which are close (on a~heuristic
level) to the Hopf maximum principle for the solutions of elliptic
PDEs (for a more detailed explanation see
Section~\ref{contrad_s}). Finally, a contradiction is obtained
using the Coarea formula, isoperimetric inequality (see
Lemma~\ref{ax-lkr11}--\ref{ax-lkr12}), and some elementary facts
from real analysis (see Appendix). In Section~\ref{contrad-u0} we
show how to modify our arguments for the case~$\ue_0\ne{\bf0}$.

The analogous result (by quite different methods) for
symmetric exterior {\bf plane} domains was established
in~\cite{kpr_ext_plane} under the additional assumption that all
connected components of the boundary intersects the symmetry axis.

\section{Notations and preliminary results}
\setcounter{theo}{0} \setcounter{lem}{0}
\setcounter{lemr}{0}\setcounter{equation}{0}

By {\it a domain} we mean an open connected set. We use standard
notations for function spaces: $W^{k,q}(\Omega)$,
$W^{\alpha,q}(\partial\Omega)$, where $\alpha\in(0,1),
k\in{\mathbb N}_0, q\in[1,+\infty]$. In our notation we do not
distinguish function spaces for scalar and vector valued
functions; it is clear from the context whether we use scalar or
vector (or tensor) valued function spaces.

For $q\ge1$ denote by $D^{k,q}(\Omega)$ the set of functions $f\in
W^{k,q}_{\loc}$ such that $\|f\|_{D^{k,q}}(\Omega)=\|\nabla^k
f\|_{L^q(\Omega)}<\infty$. Further, $D^{1,2}_0(\Omega)$ is the
closure of the set of all smooth functions having compact supports
in $\Omega$ with respect to the norm
$\|\,\cdot\,\|_{D^{1,2}(\Omega)}$, \, and $H(\Omega)=\{{\bf v}\in
D^{1,2}_0(\Omega):\, \div {\bf v}=0\}$.

Let $\Omega\subset \R^3$ be an exterior domain with $C^2$-smooth
boundary $\partial\Omega$, defined by~\eqref{Omega}.  It is well
known that functions $u\in D^{1,2}_0(\Omega)$ belong to
$L^6(\Omega)$ and, hence tend (in some sense) to zero at infinity
(see, e.g., \cite{Lad}). Moreover, in exterior domains $\Omega$
with $C^2$-smooth boundaries any vector-field ${\bf v}\in
H(\Omega)$ can be approximated in the norm
$\|\cdot\|_{D^{1,2}(\Omega)}$ by solenoidal smooth vector-fields
with compact supports (see \cite{LadSol1}).

Working with Sobolev functions we always assume that the "best
representatives" are chosen. If $w\in L^1_{\loc}(\Omega)$, then
the best representative $w^*$ is defined by
\begin{displaymath} w^*(x)=\left\{\begin{array}{rcl}\lim\limits_{R\to
0} \dashint_{B_R(x)}{w}(z)dz, & {\rm \;if\; the\; finite\; limit\; exists;} \\[4pt]
 0 \qquad\qquad\quad & \; {\rm otherwise },
\end{array}\right.
\end{displaymath}
where $\dashint_{B_R(x)}{ w}(z)dz=\dfrac{1}{\meas(
B_R(x))}\int\limits_{B_R(x)}{ w}(z)dz$, $B_R(x)=\{y\in\R^3:
|y-x|<R\}$ is a ball of radius $R$ centered at $x$. Also we use
the notation~$B_R=B_R(0)$, $S_R=\partial B_R$.

\subsection{Extension of the boundary values }

The next lemma concerns the existence of a solenoidal extensions
of boundary values.

\begin{lem}
\label{kmpLem14.1}  {\sl Let $\Omega\subset\R^3 $ be an exterior
axially symmetric domain~$(\ref{Omega})$. If ${\bf a}\in
W^{3/2,2}(\partial\Omega)$, then there exists a solenoidal
extension ${\bf A}\in W^{2,2}(\Omega)$ of ${\bf a}$ such that
$\h(x)=\s(x)$ for sufficiently large $|x|$, where
\begin{equation}
\label{SE-s}
\boldsymbol{\sigma}(x)=-\frac{x}{4\pi|x|^3}\sum_{i=1}^N{\cal F}_i
\end{equation}
  and ${\cal F}_i$ are  defined by
$(\ref{FluxN})$.   Moreover, the following estimate
\begin{equation}
\label{ExtensionEst}
\|{\bf A}\|_{W^{2,2}(\Omega)}\leq c \|{\bf a}\|_{W^{3/2,2}(\partial\Omega)}
\end{equation}
holds. Furthermore, if $\a$ is axially symmetric (axially
symmetric with no swirl), then $\h$ is axially symmetric (axially
symmetric with no swirl) too. }\end{lem}

\pr The proof is based on a standard technique.  Let
$(\mathbf{u}_s,p_s)\in  \big( H(\Omega)\cap
W^{2,2}_{loc}(\bar\Omega)\big)\times \big(L^2(\Omega)\cap
W^{1,2}_{loc}(\bar\Omega)\big)$ be the solution (see, e.g.,
\cite{Lad}) to the Stokes system
\begin{equation}
\label{SE}
\begin{array}{r@{}l}
 \Delta\mathbf{u_s} - \nabla p_s&{}= {\bf 0} \quad\hbox{\rm in }
\Omega,\\[2pt]
\div\mathbf{u_s}  &{}= 0\,\quad\hbox{\rm in } \Omega,\\[2pt]
{\u}_s&{}={\a} \quad\hbox{\rm on } \partial\Omega,\\[2pt]
\displaystyle\lim_{r\to+\infty}\mathbf{u}_s(x) & {}={\bf 0}.
\end{array}
\end{equation}
The solution~$\mathbf{u_s}$ satisfies the estimate
\begin{equation}
\label{e1}
\|\mathbf{u_s}\|_{H(\Omega)}+\|\mathbf{u_s}\|_{W^{2,2}(\Omega')}\leq  c \|{\bf a}\|_{W^{3/2,2}(\partial\Omega)}
\end{equation}
for any bounded  $\Omega'$ with $\bar \Omega'\subset\Omega$ (see,
e.g.,  \cite{Lad}). Define
$$
{\bf v}(x)={\bf u}_s(x)-\boldsymbol{\sigma}(x).
$$
Take  $R>0$ such that $\partial\Omega \Subset B_R$. Then
\begin{equation}
\label{NFF} \displaystyle\displaystyle \int_{\partial B_R}{\bf
v}\cdot{\bf n} \,dS=0
\end{equation}

Let $\zeta=\zeta(|x|)$ be a smooth cut-off  function, equal to 1
in $B_{R}$ and vanishing outside $B_{2R}$. Then $\div(\zeta{\bf v})=\nabla\zeta\cdot{\bf v}\in
\wotwo(B_{2R}\setminus B_{R})$. Because of (\ref{NFF}) the
equation
\begin{equation} \label{div}
\div \boldsymbol{\xi}=-\div(\zeta{\bf v}), \quad
\boldsymbol{\xi}|_{\partial(B_{2R}\setminus B_{R})}=0
\end{equation}
has a solution $\boldsymbol{\xi}\in \wotwotwo(B_{2R}\setminus
B_{R})$ (see \cite{Bogovski}) and
\begin{equation}
\label{e2}
\|\mathbf{\boldsymbol{\xi}}\|_{\wotwotwo(B_{2R}\setminus
B_{R})}\leq c \|{\bf v}\|_{\wotwo(B_{2R}\setminus B_{R})}\leq c
\|{\bf a}\|_{W^{3/2,2}(\partial\Omega)}.
\end{equation}

 The field
$$
 \h=\begin{cases}
{\bf u}_s,&x\in B_R\cap\Omega,\\
 \boldsymbol{\sigma}+\boldsymbol{\xi} +\zeta{\bf v}, &x\in B_{2R}\setminus B_R,\\
\boldsymbol{\sigma},&x\in \Omega\setminus B_{2R}\\
\end{cases}
$$
is the desired solenoidal extension of~${\bf a}$ in~$\Omega$. If
$\a$ is axially symmetric (axially symmetric with no swirl), then,
according to results of Section 2.1 in \cite{kpr_a_arx}, the
extension $\mathbf{u}_s$  and the solution $\boldsymbol{\xi}$ to
the divergence equation \eqref{div} are axially symmetric (axially
symmetric with no swirl). Thus, in this case the extension $\h$ is
axially symmetric (axially symmetric with no swirl) too. $\qed$

\subsection{A uniform estimate of the pressure in the
Stokes system}

This subsection is rather technical: the estimates below are not
considered as new or sharp, but they sufficient for our purposes.

\begin{lem}
\label{LemPr1}  {\sl Let $\Omega\subset\R^3 $ be an exterior
axially symmetric domain $(\ref{Omega})$, ${\bf f}\in
L^{3/2}(\Omega_R)$ and  ${\bf a}\in W^{3/2,2}(\partial\Omega)$.
Take $R_0>0$ such that $\partial\Omega\Subset\frac12B_{R_0}$.
Suppose that $R\ge R_0$ and $\ue\in W^{1,2}(\Omega_R)$ is a
solution to the~Stokes system
\begin{equation}
\label{SE-b}
\begin{array}{r@{}l}
 \nu\Delta\mathbf{u} - \nabla p&{}= {\bf f} \quad\hbox{\rm in }
\Omega_R,\\[2pt]
\div\mathbf{u}  &{}= 0\,\quad\hbox{\rm in } \Omega_R,\\[2pt]
{\u}&{}={\a} \quad\hbox{\rm on } \partial\Omega,\\[2pt]
\ue & {}=\boldsymbol{\sigma}\quad\hbox{\rm on } \partial B_R,
\end{array}
\end{equation}
where $\Omega_R=\Omega\cap B_R$ and $\boldsymbol{\sigma}(x)$ is
defined by~(\ref{SE-s}). Then the estimates
\begin{equation}
\label{L-u1} \|\ue\|_{L^6(\Omega_R)}\le c \bigl(\|\nabla
\ue\|_{L^2(\Omega_R)}+\|\a\|_{W^{3/2,2}(\partial\Omega)}\bigr),
\end{equation}
\begin{equation}
\label{L-p1}
\begin{array}{lcr}
\|\nabla^2 {\bf u}\|_{L^{3/2}(\Omega_R)}+\|\nabla
p\|_{L^{3/2}(\Omega_R)}\le
\\
\\
\qquad\qquad\le c \bigl(\|\nabla
\ue\|_{L^2(\Omega_R)}+\|\fe\|_{L^{3/2}(\Omega_R)}+\|\a\|_{W^{3/2,2}(\partial\Omega)}\bigr)
\end{array}
\end{equation}
hold with the constant~$c$ independent of $\a,\ue,\fe,$ and $R$.
}\end{lem}

\pr Rewriting the Stokes system~(\ref{SE-b}) for the new function
$\ue'=\ue-\h$, where $\h$ is the solenoidal extension of~$\a$ from
Lemma~\ref{kmpLem14.1}, we can assume, without loss of generality,
that~$\a=0$, i.e.,
\begin{equation}
\label{SE-b0}
\begin{array}{r@{}l}
 \nu\Delta\mathbf{u} - \nabla p&{}= {\bf f} \quad\hbox{\rm in }
\Omega_R,\\[2pt]
\div\mathbf{u}  &{}= 0\,\quad\hbox{\rm in } \Omega_R,\\[2pt]
{\ue}&{}={\bf 0} \quad\hbox{\rm on } \partial\Omega_R.
\end{array}
\end{equation}
Then the first estimate~(\ref{L-u1}) follows easily from the well
known inequality
$$
\|g\|_{L^6(\R^3)}\le c\|\nabla g\|_{L^2(\R^3)}\qquad\forall g\in
C^\infty_0(\R^3)
$$
and \eqref{ExtensionEst}.

Let us prove~(\ref{L-p1}). Let $R=2^lR_0$ with $l\in\N$. Denote
$$ \omega_k=\{x\in \R^3: R_02^{k-1}\le|x|\le R_02^{k}\}, \;
k=1,2,\ldots,l; \omega_0=\{x\in\Omega: |x|\le R_0\},
$$
$$
\widehat\omega_k=\{x\in \R^3: R_02^{k-2}<|x|<R_02^{k+1}\}, \;
k=1,2,\ldots; \widehat\omega_0=\{x\in\Omega: |x|<2R_0\}.
$$
Obviously, $\omega_k\Subset\widehat\omega_k$ and
$\Omega_R=\bigcup\limits_{k=0}^l \omega_k$. Consider the
system~(\ref{SE-b0}) in $\widehat\omega_k$, $k\le l-1$. After the
scaling $y=\dfrac{x}{R_02^k}$ we get the Stokes equations in the
domain $\widehat\sigma_0=\{y:\; \frac14<|y|<2\}$:
$$
-\nu \Delta_y {\bf u}+\nabla_y \tilde p=\tilde{\bf f},\ \ \
\div\nolimits_y {\bf u}=0,
$$
where $\tilde p=2^kR_0p,\; \tilde{\bf f}=2^{2k}R_0^2{\bf f} $. Let
$\sigma_0=\{y:\; \frac12<|y|<1\}$. By the local estimate for
ADN-elliptic problems  ( see \cite{Agmon}, \cite{Solonnikov}) we
have
\begin{equation}
\label{p15}
\begin{array}{lcr}
\|{\bf u}\|_{W^{2,\frac32}(\sigma_0)} +\|\nabla_y\tilde
p\|_{L^{\frac32}(\sigma_0)}\leq\\
\\
\qquad\qquad\leq c\Big( \|\tilde
{\bf f}\|_{L^\frac32(\widehat\sigma_0)} +\|{\bf
u}\|_{L^\frac32(\widehat\sigma_0)} +\|\tilde p-\tilde
p_{0}\|_{L^\frac32(\widehat\sigma_0)}\Big),
\end{array}
\end{equation}
where $\tilde
p_0=\dfrac{1}{|\widehat\sigma_0|}\int\limits_{\widehat\sigma_0}\tilde
p(y)dy$.

Consider the functional
$$
H({\bf  w})=\int_{\widehat\sigma_0}(\tilde p(y)-\tilde p_{0})\;{\rm
div}\,{\bf w}\,dy,\quad \forall {\bf w}\in
{\mathop{\mathaccent23{W}^{1,3}}\nolimits}(\widehat\sigma_0).
$$
Using Stokes equations  and integrating by parts we obtain
$$
|H({\bf w})|=\bigg\vert\int_{\widehat\sigma_0}\nabla_y\tilde
p\cdot{\bf
w}\,dy\bigg\vert\leq\nu\bigg\vert\int_{\widehat\sigma_0}\nabla_y{\bf
u}\cdot\nabla_y{\bf w}\,dy\bigg\vert
$$
$$
+\bigg\vert\int_{\widehat\sigma_0}\tilde{\bf  f}\cdot{\bf
w}\,dy\bigg\vert\leq c\Big( \Vert\tilde{\bf
f}\Vert_{L^\frac32(\widehat\sigma_0)}+\Vert\nabla_y{\bf
u}\Vert_{L^\frac32(\widehat\sigma_0)} \Big)\Vert\nabla_y{\bf
w}\Vert_{L^3(\widehat\sigma_0)}.
$$
The norm of the functional $H$ is equivalent to $\Vert\tilde
p-\tilde p_{0}\Vert_{L^\frac32(\widehat\sigma_0)}$ (see, e.g.,  \cite{Pil1}, \cite{Pil2}). Hence, the last
estimate gives
$$
\Vert\tilde p-\tilde p_{0}\Vert_{L^\frac32(\widehat\sigma_0)}\leq c\Big(
\Vert\tilde{\bf f}\Vert_{L^\frac32(\widehat\sigma_0)}+\Vert\nabla_y{\bf
u}\Vert_{L^\frac32(\widehat\sigma_0)} \Big)
$$
and inequality \eqref{p15} takes the form
$$
\|{\bf u}\|_{W^{2,\frac32}(\sigma_0)} +\|\nabla_y\tilde
p\|_{L_\frac32(\sigma_0)}\leq c\Big( \|\tilde {\bf
f}\|_{L_\frac32(\widehat\sigma_0)} +\|{\bf
u}\|_{W^{1,\frac32}(\widehat\sigma_0)} \Big).
$$
In particular,
$$
\|\nabla^2_y\ue\|_{L^{\frac32}(\sigma_0)}+\|\nabla_y\tilde
p\|_{L_\frac32(\sigma_0)}\leq c\Big( \|\tilde
{f}\|_{L_\frac32(\widehat\sigma_0)} +\|{\bf
u}\|_{W^{1,\frac32}(\widehat\sigma_0)} \Big).
$$

Returning to coordinates $x$,  we get
\begin{equation}
\label{p17}
\begin{array}{lcr}
\|\nabla^2_x\ue\|_{L^{\frac32}(\omega_k)}+ \|\nabla_x p\|_{L^\frac32(\omega_k)}\leq\\
\\
\qquad\qquad \leq c\Big( \|{\bf f}\|_{L^\frac32(\widehat\omega_k)}
+2^{-2k}\|{\bf u}\|_{L^\frac32(\widehat\omega_k )}
+2^{-k}\|\nabla_x{\bf u}\|_{L^\frac32(\widehat\omega_k )} \Big).
\end{array}
\end{equation}
Estimates \eqref{p17} hold for $1\leq k\leq l-1$. For $k=l$ we
obtain, by the same argument, the following estimate
\begin{equation}
\label{p171}
\begin{array}{lcr}
\|\nabla^2_x\ue\|_{L^{\frac32}(\omega_l)}+ \|\nabla_x p\|_{L^\frac32(\omega_l)}\leq\\
\\
\leq c\Big(
\|{\bf f}\|_{L^\frac32(\omega_{l-1}\cup\omega_l)} +2^{-2l}\|{\bf
u}\|_{L^\frac32(\omega_{l-1}\cup\omega_l)} +2^{-l}\|\nabla_x{\bf
u}\|_{L^\frac32(\omega_{l-1}\cup\omega_l)} \Big).
\end{array}
\end{equation}
Furthermore, by local estimates (see~\cite{Agmon},
\cite{Solonnikov}) for the solution to (\ref{SE-b0})
we have
\begin{equation}
\label{p19}
\begin{array}{lcr}
\|\nabla^2_x\ue\|_{L^{\frac32}(\omega_0)}+\|\nabla_x p\|_{L^{\frac32}(\omega_0)}\leq\\
\\
\leq c\Big( \| {\bf
f}\|_{L^\frac32(\omega_{0}\cup\omega_1)}+\|{\bf
u}\|_{L^\frac32(\omega_{0}\cup\omega_1)}  +\|\nabla_x{\bf
u}\|_{L^\frac32(\omega_{0}\cup\omega_1)}\Big).
\end{array}
\end{equation}
(recall that ${\bf u}|_{\partial \Omega}=0$). Summing
estimates (\ref{p17})--\eqref{p19} by $k$ and taking into account
that $r\sim 2^k$ for $x\in \omega_k$,  we derive
\begin{equation}\label{p20}
\begin{array}{lcr}
\|\nabla^2_x\ue\|_{L^{\frac32}(\Omega_R)}+ \|\nabla_x
p\|_{L^\frac32(\Omega_R)}\leq
\\
\\
 \leq
c\Big( \|{\bf f}\|_{L^\frac32(\Omega_R)} +\||x|^{-2}{\bf
u}\|_{L^\frac32(\Omega_R )}+\||x|^{-1}\nabla_x{\bf
u}\|_{L^\frac32(\Omega_R )} \Big)\leq
\\
\\
\leq c'\Big( \|{\bf
f}\|_{L^\frac32(\Omega_R)} +\||x|^{-1}{\bf
u}\|_{L^2(\Omega_R )} +\|\nabla_x{\bf u}\|_{L^2(\Omega_R )}
\Big)\leq
\\
\\
\leq  c''\Big( \|{\bf
f}\|_{L^\frac32(\Omega_R)} +\|\nabla_x{\bf u}\|_{L^2(\Omega_R )}
\Big),
\end{array}
\end{equation}
where the constant $c''$ is independent of $R$.
Here we have used the H\"older inequality and the well know estimate
$$
\int_{\R^3}\frac{|g(x)|^2}{|x|^2}dx\leq c\int_{\R^3}|\nabla g(x)|^2dx\quad \forall \; g\in C_0^\infty(\R^3).
$$

The required estimate~(\ref{L-p1}) for the solution of the
nonhomogeneous boundary value problem \eqref{SE-b} follows from
(\ref{p20}) and the estimate \eqref{ExtensionEst} for the
extension ${\bf A}$ constructed in  Section 2.1.
 \qed

\begin{lemA}
\label{LemPr2}  {\sl Let $\Omega\subset\R^3 $ be an exterior
axially symmetric domain $(\ref{Omega})$, ${\bf f}\in
L^{3/2}(\Omega_R)$ and  ${\bf a}\in W^{3/2,2}(\partial\Omega)$.
Take $R_0>0$ such that $\partial\Omega\Subset\frac12B_{R_0}$ and
$R\ge R_0$. Let $\ue\in W^{1,2}(\Omega_R)$ be a solution to
the~Navier--Stokes system
\begin{equation}
\label{NS-b}
\begin{array}{r@{}l}
 -\nu \Delta{\bf u}+\big({\bf u}\cdot \nabla\big){\bf u} +\nabla p&{}= {\bf f} \quad\hbox{\rm in }
\Omega_R,\\[2pt]
\div\mathbf{u}  &{}= 0\,\quad\hbox{\rm in } \Omega_R,\\[2pt]
{\u}&{}={\a} \quad\hbox{\rm on } \partial\Omega,\\[2pt]
\ue & {}=\boldsymbol{\sigma}\quad\hbox{\rm on } \partial B_R,
\end{array}
\end{equation}
where   $\boldsymbol{\sigma}(x)$ is defined by~(\ref{SE-s}). Then
the estimates
\begin{equation}
\label{LLL-u1} \|\ue\|_{L^6(\Omega_R)}\le c \bigl(\|\nabla
\ue\|_{L^2(\Omega_R)}+\|\a\|_{W^{3/2,2}(\partial\Omega)}\bigr),
\end{equation}
\begin{equation}
\label{LLL-p1}
\begin{array}{lcr} \|\nabla^2 {\bf u}\|_{L^{3/2}(\Omega_R)}+ \|\nabla p\|_{L^{3/2}(\Omega_R)}\le c \Big(
\|\fe\|_{L^{3/2}(\Omega_R)}+\|\nabla\ue\|_{L^2(\Omega_R)}+\\
\\
\qquad+\|\a\|_{W^{3/2,2}(\partial\Omega)}+\bigl(\|\nabla\ue\|_{L^2(\Omega_R)}+\|\a\|_{W^{3/2,2}(\partial\Omega)}\bigr)^2
\Big)
\end{array}
\end{equation}
hold. Here the constant~$c$ is independent of $\a,\ue,\fe,R$.
}\end{lemA}

\pr Estimate \eqref{LLL-u1} follows by the same argument as
\eqref{L-u1}. In order to prove \eqref{LLL-p1}, we consider the
Navier--Stokes system \eqref{NS-b} as the Stokes one with the
right-hand side  ${\bf f}'={\bf f}-({\bf u}\cdot\nabla){\bf u}$.
Since
$$
\|({\bf u}\cdot\nabla){\bf u}\|_{L^{3/2}(\Omega_R)}\leq \|{\bf u}\|_{L^{6}(\Omega_R)}\|\nabla{\bf u}\|_{L^{2}(\Omega_R)}
\leq c\|\nabla{\bf u}\|_{L^{2}(\Omega_R)}^2,
$$
estimate \eqref{LLL-p1} follows from \eqref{L-p1}. \qed

\subsection{On Morse-Sard and Luzin N-properties of Sobolev
functions from $W^{2,1}$} \label{SMS}

Let us recall some classical differentiability properties of
Sobolev functions.

\begin{lem}[see Proposition~1 in \cite{Dor}]
\label{kmpThDor} {\sl Let  $\psi\in W^{2,1}(\R^2)$. Then the
function~$\psi$ is continuous and there exists a set $A_{\psi}$
such that $\mathfrak{H}^1(A_{\psi})=0$, and the function $\psi$ is
differentiable (in the classical sense) at each $x\in\R^2\setminus
A_{\psi}$. Furthermore, the classical derivative at such points
$x$ coincides with $\nabla\psi(x)=\lim\limits_{r\to 0}
\dashint_{B_r(x)}{ \nabla\psi}(z)dz$, and \ $\lim\limits_{r\to
0}\dashint\nolimits_{B_r(x)}|\nabla\psi(z)-\nabla\psi(x)|^2dz=0$.}
\end{lem}

Here and henceforth we denote by $\mathfrak{H}^1$ the
one-dimensional Hausdorff measure, i.e.,
$\mathfrak{H}^1(F)=\lim\limits_{t\to 0+}\mathfrak{H}^1_t(F)$,
where $\mathfrak{H}^1_t(F)=\inf\{\sum\limits_{i=1}^\infty {\rm
diam} F_i:\, {\rm diam} F_i\leq t, F\subset
\bigcup\limits_{i=1}^\infty F_i\}$.

The next theorem have been proved recently by J. Bourgain,
M.~Korobkov and J. Kristensen \cite{korob}.

\begin{theo}
\label{kmpTh1.1}{\sl  Let  ${\mathcal D}\subset\R^2$ be a bounded
domain with Lipschitz boundary and  $\psi\in W^{2,1}({\mathcal
D})$. Then

{\rm (i)} $\mathfrak{H}^1(\{\psi(x)\,:\,x\in\bar{\mathcal
D}\setminus A_\psi\,\,\&\,\,\nabla \psi(x)=0\})=0$;

{\rm (ii)} for every $\varepsilon>0$ there exists $\delta>0$ such
that for any set $U\subset \bar{\mathcal D}$ with
$\mathfrak{H}^1_\infty(U)<\delta$ the inequality
$\mathfrak{H}^1(\psi(U))<\varepsilon$ holds;

{\rm (iii)} for $\mathfrak{H}^1$--almost all $y\in
\psi(\bar{\mathcal D})\subset \mathbb{R}$ the preimage
$\psi^{-1}(y)$ is a finite disjoint family of $C^1$--curves $S_j$,
$ j=1, 2, \ldots, N(y)$. Each $S_j$ is either a cycle in
${\mathcal D}$ $($i.e., $S_j\subset{\mathcal D}$ is homeomorphic
to the unit circle $\mathbb{S}^1)$ or it is a simple arc with
endpoints on $\partial{\mathcal D}$ $($in this case $S_j$ is
transversal to $\partial{\mathcal D}\,)$. }
\end{theo}

\subsection{Some facts from topology}
\label{Kronrod-s}

We shall need  some topological definitions and results. By  {\it
continuum} we mean a compact connected set. We understand
connectedness  in the sense of general topology. A~subset of a
topological space is called {\it an arc} if it is homeomorphic to
the unit interval~$[0,1]$.

Let us shortly present  some results from the classical paper of
A.S.~Kron-\\rod~ \cite{Kronrod} concerning level sets of
continuous functions. Let ${Q}=[0,1]\times[0,1]$ be a square in
$\mathbb{R}^2$ and let $f$ be a continuous function  on ${Q}$.
Denote by $E_t$ a level set of the function $f$, i.e.,
$E_t=\{x\in{Q}: f(x)=t\}$. A component $K$  of the level set $E_t$
containing a point $x_0$ is a maximal connected subset of $E_t$
containing $x_0$. By $T_f$ denote a family of all connected
components of level sets of~$f$. It was established
in~\cite{Kronrod} that $T_f$ equipped by a natural
topology\footnote{The convergence in $T_f$ is defined as follows:
$T_f\ni C_i\to C$ iff $\sup\limits_{x\in C_i}\dist(x,C)\to0$.} is
a one-dimensional topological tree\footnote{A locally connected
continuum $T$ is called~{\it a topological tree}, if it does not
contain a curve homeomorphic to a circle, or, equivalently, if any
two different points of~$T$ can be joined by a unique arc. This
definition implies that $T$ has topological dimension~1.}.
Endpoints of this tree\footnote{A point of a continuum~$K$ is
called an {\it endpoint of $K$} (resp., {\it a branching point
of~$K$}) if its topological index equals~1 (more or equal to~$3$
resp.). For a~topological tree~$T$ this definition is equivalent
to the~following: a point $C\in T$ is an~endpoint of~$T$ (resp., {
a branching point of~$T$}), if the set $T\setminus\{C\}$ is
connected (resp., if $T\setminus\{C\}$ has more than two connected
components).} are the components~$C\in T_f$ which do not
separate~$Q$, i.e., $Q\setminus C$ is a connected set. Branching
points of the tree are the components $C\in T_f$ such that
$Q\setminus C$ has more than two connected components (see
\cite[Theorem 5]{Kronrod}). By results of \cite[Lemma~1]{Kronrod},
the set of all branching points of~$T_f$ is at most countable. The
main property of a tree is that any two points could be joined by
a unique arc. Therefore, the same is true for~$T_f$.

\begin{lem} [see Lemma~13 in \cite{Kronrod}]
\label{kmpLem6} If $f\in C(Q)$, then for any two different points
$A\in T_f$ and $B\in T_f$, there exists a unique arc
$J=J(A,B)\subset T_f$ joining $A$ to $B$. Moreover, for every
inner point $C$ of this arc the points $A,B$ lie in  different
connected components of the set $T_f\setminus\{C\}$.
\end{lem}

We can reformulate the above Lemma in the following equivalent
form.

\begin{lem} \label{kmpLem7}{\sl If  $f\in C(Q)$, then for
any two different points $A,B\in T_f$, there exists a~continuous
injective function $\varphi:[0,1]\to T_f$ with the properties

{\rm (i)} $\varphi(0)=A$, $\varphi(1)= B$;

{\rm (ii)} for any $t_0\in[0,1]$,
$$
\lim\limits_{[0,1]\ni t\to t_0}\sup\limits_{x\in
\varphi(t)}\dist(x,\varphi(t_0))\to0;
$$

{\rm (iii)}  for any $t\in(0,1)$ the sets $A,B$ lie in different
connected components of the set \ $Q\setminus\varphi(t)$.}
\end{lem}

\begin{lemr}
\label{kmpRem2} {\rm If in  Lemma~\ref{kmpLem7} $f\in W^{2,1}(Q)$,
then by Theorem~\ref{kmpTh1.1}~(iii), there exists a dense subset
$E$ of $(0,1)$ such that $\varphi(t)$ is a $C^1$-- curve for every
$t\in E$. Moreover, $\varphi(t)$ is either a cycle   or a simple
arc with endpoints on $\partial Q$.}
\end{lemr}

\begin{lemr}
\label{kmpRem1.2} {\rm All results of
Lemmas~\ref{kmpLem6}--\ref{kmpLem7} remain valid for level sets of
continuous functions $f:\overline\D_0\to\R$, where
$\overline\D_0\subset\R^2$ is a compact set homeomorphic to the
unit square $Q=[0,1]^2$.}
\end{lemr}
\

\section{Leray's argument ``reductio ad absurdum'' }
\label{poet} \setcounter{theo}{0} \setcounter{lem}{0}
\setcounter{lemr}{0}\setcounter{equation}{0}

Let us consider the Navier--Stokes problem (\ref{NS}) with ${\bf
f}\in W^{1,2}(\Omega)\cap L^{6/5}(\Omega)$ in the $C^2$-smooth
axially symmetric exterior domain $\Omega\subset\R^3$  defined
by~(\ref{Omega}).
 Without loss of
generality, we may assume that ${\bf f}={\rm curl}\,{\bf b}\in
W^{1,2}(\Omega)\cap L^{6/5}(\Omega)$.\footnote{By the
Helmholtz--Weyl decomposition, $\fe$ can be represented as the sum
$\fe={\rm curl}\,{\bf b}+\nabla \varphi$ with ${\rm curl}\,{\bf
b}\in W^{1,2}(\Omega)\cap L^{6/5}(\Omega)$, and the
gradient part is included then into the pressure term (see, e.g.,
\cite{Lad}, \cite{Galdibook}).}

We shall prove
Theorem~\ref{kmpTh4.2} for
$$
{\bf u}_0={\bf 0}.
$$
The proof for ${\bf u}_0\ne {\bf 0}$ follows the same steps with
minor standard modification, see Section~\ref{contrad-u0}.

By a {\it weak solution} of problem (\ref{NS}) we mean a function
${\bf u}$ such that ${\bf w}={\bf u}-{\bf A}\in H(\Omega)$ and the
integral identity
\begin{equation}\label{weaksolution}
\begin{array}{lcr}
\nu\int\limits_\Omega\nabla{\bf w}\cdot\nabla\bfeta
\,dx=-\nu\int\limits_\Omega\nabla{\bf A}\cdot\nabla\bfeta
\,dx
-\int\limits_\Omega\big({\bf A} \cdot \nabla\big){\bf
A}\cdot\bfeta\,dx\\
\\
- \int\limits_\Omega\big({\bf A} \cdot
\nabla\big){\bf w} \cdot\bfeta\,dx
- \int\limits_\Omega\big({\bf w}\cdot \nabla\big){\bf
w}\cdot\bfeta\,dx\\
\\
- \int\limits_\Omega\big({\bf w}\cdot
\nabla\big){\bf A} \cdot\bfeta\,dx
+\int\limits_\Omega{\bf f} \cdot
\bfeta\,dx
\end{array}
\end{equation}
holds for any $\bfeta\in J_0^\infty(\Omega)$, where $
J_0^\infty(\Omega)$ is a set of all infinitely smooth solenoidal
vector-fields with compact support in $\Omega$. Here ${\bf A} $ is
the extension of the boundary data constructed in Lemma
\ref{kmpLem14.1}. We shall look for the axially symmetric (axially
symmetric with no swirl) weak solution of problem (\ref{NS}). We
find this solution as a limit of weak solution to the
Navier--Stokes problem in a sequence of bounded domain $\Omega_k$
that in the limit exhaust  the unbounded domain $\Omega$. The
following result concerning the solvability of the Navier-Stokes
problem in axially symmetric bounded domains was proved in
\cite{kpr_a_ann}.

\begin{theo} \label{Th_Ex_b} {\sl
Let $\Omega'=\Omega_0\setminus
\bigl(\bigcup\limits_{j=1}^N\bar\Omega_j\bigr)$ be an axially
symmetric bounded domain in $\R^3$ with multiply connected
$C^2$-smooth boundary $\partial\Omega'$ consisting of $N+1$
disjoint components $\Gamma_j=\partial\Omega_j$, $j=0,\dots, N$.
If ${\bf f}\in W^{1,2}(\Omega')$ and ${\bf a}\in
W^{3/2,2}(\partial\Omega')$ are axially symmetric and ${\bf a}$
satisfies
$$
\int_{\partial\Omega'}{\bf a}\cdot{\bf n}\,dS=0,
$$
then $(\ref{NS}_1)$--$(\ref{NS}_3)$ with $\Omega=\Omega'$ admits
at least one weak axially symmetric solution ${\bf u}\in
W^{1,2}(\Omega')$. Moreover, if ${\bf f}$ and ${\bf a}$
 are axially symmetric with
no swirl, then the problem $(\ref{NS}_1)$--$(\ref{NS}_3)$ with
$\Omega=\Omega'$ admits at least one weak axially symmetric
solution with no swirl. }
\end{theo}

 Consider the sequence of boundary value problems
\begin{equation}
\label{NSESS}
\begin{array}{r@{}l}
 \nu \Delta\widehat{\bf w}_k - (\widehat{\bf w}_k +{\bf A} )\cdot\nabla(\widehat{\bf w}_k +{\bf A} )
 +\nu \Delta{\bf A} -\nabla \widehat p_k &{}= {\bf f} \quad\hbox{\rm in }
\Omega_k,\\[2pt]
\div\widehat{\bf w}_k   &{}= 0\,\quad\hbox{\rm in } \Omega_k,\\[2pt]
\widehat{\bf w}_k &{}={\bf 0} \quad\hbox{\rm on }
\partial\Omega_k,\end{array}
\end{equation}
where $\Omega_k=B_k\cap\Omega$ for $k\ge k_0$, $B_k=\{x:\ \ |x|<
k\}$, $\frac12B_{k_0}\supset \bigcup\limits_{i=1}^N\bar\Omega_i$.
By Theorem~\ref{Th_Ex_b},  each problem~(\ref{NSESS}) has an
axially symmetric solution $\w_k\in H(\Omega_k)$ satisfying the
integral identity
\begin{equation}\label{weaksolution_k}
\begin{array}{lcr}
\displaystyle\nu\int\limits_\Omega\!\!\nabla\widehat{\bf
w}_k\cdot\nabla\bfeta dx\!=\!\int\limits_\Omega\!\!{\bf f} \cdot \bfeta\,dx-\nu\int\limits_\Omega\nabla{\bf A}
\cdot\nabla\bfeta dx\!-\!\int\limits_\Omega\!\!\big({\bf A} \cdot
\nabla\big){\bf A}\cdot\bfeta dx\\[12pt]
\displaystyle - \int\limits_\Omega\big({\bf A}\cdot
\nabla\big)\widehat{\bf w}_k \cdot\bfeta\,dx-
\int\limits_\Omega\big(\widehat{\bf w}_k\cdot
\nabla\big)\widehat{\bf w}_k\cdot\bfeta\,dx-
\int\limits_\Omega\big(\widehat{\bf w}_k\cdot \nabla\big){\bf A}
\cdot\bfeta\,dx
\end{array}
\end{equation}
for all $\bfeta\in H(\Omega_k)$. Here we have assumed that
$\widehat{\bf w}_k$ and $\bfeta$ are extended by zero  to  the
whole domain $\Omega$.

Assume that there is a positive constant $c$ independent of $k$
such that
\begin{equation}
\label{kkkaf} \int_\Omega|\nabla\widehat{\bf w}_k|^2\le c
\end{equation}
(possibly along a subsequence of $\{\widehat{\bf
w}_k\}_{k\in{\mathbb N}}$). The
estimate~\eqref{kkkaf} implies  the existence of a solution to
problem~(\ref{NS}). Indeed, the sequence $\widehat{\bf w}_k$ is
bounded in $H(\Omega)$. Hence, $\widehat{\bf w}_k$ converges
weakly (modulo a subsequence) in $H(\Omega)$ and strongly in
$L^q_{\rm loc}(\overline\Omega)$ $(q<6)$ to a function
$\widehat{\bf w}\in H(\Omega)$. Taking any test function $\bfeta$
with compact support, we can find $k$ such that ${\rm
supp}\,\bfeta\subset \Omega_k$. Thus, we can pass to a limit as
$k\to\infty$ in \eqref{weaksolution_k} and we obtain for the limit
function $\widehat{\bf w}$ the integral identity
\eqref{weaksolution}. Then, by definition,  $\u=\widehat{\bf
w}+\h$ is a weak solution to the Navier--Stokes
problem~(\ref{NS}). Thus, to prove the assertion of
Theorem~\ref{kmpTh4.2}, it is sufficient to establish the uniform
estimate~(\ref{kkkaf}).

We shall prove (\ref{kkkaf}) following a classical {\it reductio
ad absurdum\/} argument of J. Leray \cite{Leray} and O.A.
Ladyzhenskaia \cite{Lad}. Indeed, if (\ref{kkkaf}) is not true,
then there exists a sequence $\{\widehat\w_k\}_{k\in{\mathbb N}}$
such that
$$
\lim_{k\to+\infty}J_k^2=+\infty,\quad
J_k^2=\int_\Omega|\nabla\widehat{\bf w}_k|^2.
$$
The sequence ${\bf w}_k=\widehat{\bf w}_k/J_k$ is   bounded in
$H(\Omega)$ and it  holds
\begin{equation}
\label{NSESS'}
\begin{array}{lcr}
 \dfrac{\nu}{J_k}\int\limits_{\Omega}\nabla{\bf w}_k\cdot\nabla\bfeta\,dx= -\int\limits_\Omega ({\bf w}_k\cdot\nabla){\bf w}_k\cdot\bfeta\,dx+ \dfrac{1}
{J_k}\int\limits_\Omega({\bf w}_k\cdot\nabla)\bfeta\cdot{\bf A}\,dx\\
\\
+\dfrac{1}{J_k}\int\limits_\Omega ({\bf
A}\cdot\nabla)\bfeta\cdot{\bf
w}_k\,dx+\dfrac{1}{J_k^2}\int\limits_\Omega({\bf
A}\cdot\nabla)\bfeta\cdot{\bf
A}\,dx\\
\\
-\dfrac{\nu}{J_k^2}\int\limits_\Omega\nabla{\bf
A}\cdot\nabla\bfeta\,dx+\dfrac{1}{J_k^2}\int\limits_\Omega{\bf
f}\cdot\bfeta\,dx
\end{array}
\end{equation}
for all $\bfeta\in H(\Omega_k)$. Extracting  a subsequence (if
necessary) we can assume that ${\bf w}_k$ converges weakly in
$H(\Omega)$ and strongly in $L^q_{\rm loc} (\overline\Omega)$
$(q<6)$  to a vector field ${\bf v}\in H(\Omega)$ with
\begin{equation}\label{norm1}
\int\limits_\Omega|\nabla {\bf v}|^2\le1.
\end{equation} Fixing in \eqref{NSESS'} a
solenoidal smooth  $\bfeta$ with compact support and letting
$k\to+\infty$ we get
\begin{equation}\label{Eweak}
\int\limits_\Omega({\bf v}\cdot\nabla){\bf
v}\cdot\bfeta\,dx=0\quad\forall \,\bfeta\in J_0^\infty(\Omega),
\end{equation}
Hence, ${\bf v}\in H(\Omega)$ is a weak solution to the Euler
equations, and for some $ p\in D^{1,3/2}(\Omega)$ the pair $({\bf
v},  p)$ satisfies the Euler equations almost everywhere:
\begin{equation}
\label{2.1}\left\{\begin{array}{rcl} \big({\bf
v}\cdot\nabla\big){\bf v}+\nabla p & = & 0 \qquad \ \ \
\hbox{\rm in }\;\;\Omega,\\[4pt]
\div{\bf v} & = & 0\qquad \ \ \ \hbox{\rm in }\;\;\Omega,
\\[4pt]
{\bf v} &  = & 0\ \
 \qquad\  \hbox{\rm on }\;\;\partial\Omega.
\end{array}\right.
\end{equation}
Adding some constants to $p$ (if necessary) by virtue of the
Sobolev inequality (see, e.g., \cite{Galdibook} II.6) we may
assume without loss of generality that
\begin{equation}
\label{BOIN} \|p\|_{L^3(\Omega)}\le \const.
\end{equation}

Put $\nu_k=(J_k)^{-1}\nu$. Multiplying equations~(\ref{NSESS}) by
$\frac{1}{J^2_k}=\frac{\nu^2_k}{\nu^2}$, we see that the pair
$\big(\ue_k= \frac{1}{J_k}\widehat{\w}_k+\frac{1}{J_k}{\bf A}, \
p_k=\frac{1}{J^2_k}\widehat p_k\big)$ satisfies the following
system

\begin{equation}\label{NSk}
\left\{\begin{array}{rcl}-\nu_k\Delta{\bf u}_k +\big({\bf
u}_k\cdot\nabla\big){\bf u}_k+\nabla p_k & = & \fe_k\qquad
\hbox{\rm in }\;\;\Omega_k,
\\[4pt]
\div{\bf u}_k & = & 0 \;\qquad \hbox{\rm in }\;\;\Omega_k,
\\[4pt]
{\bf u}_k &  = & {\bf a}_k
 \quad\  \hbox{\rm on }\;\;\partial\Omega_k,
\end{array}\right.
\end{equation}
where $\fe_k=\frac{\nu_k^2}{\nu^2}\,{\bf f}$, \ ${\bf
a}_k=\frac{\nu_k}\nu\,{\bf A}$, $\ue_k\in W^{3,2}_{\loc}(\Omega)$,
$p_k\in W^{2,2}_{\loc}(\Omega)$\footnote{The interior  regularity
of the solution depends on the regularity of $\fe\in
W^{1,2}(\Omega)$, but not on the regularity of the boundary value
${\bf a}$, \ see \cite{Lad}.}. From Corollary~\ref{LemPr2} we have
\begin{equation}
\label{epgr33}\|\ue_k\|_{L^6(\Omega_k)}\leq\const.
\end{equation}
For $R>R_0$ denote $\Omega_R=\Omega\cap B(0,R)$. Since by H\"older
inequality $\|\fe\|_{L^{3/2}(\Omega_R)}\le
\|\fe\|_{L^{2}(\Omega_R)}\sqrt{R}$, Corollary~\ref{LemPr2} implies
also
\begin{equation}
\label{epgr34}\|\nabla p_k\|_{L^{3/2}(\Omega_R)}\leq
C\sqrt{R}\qquad\forall R\in[R_0,R_k],
\end{equation}
where $C$ does not depend\footnote{Of course, the above
assumptions
 ${\bf
f}\in W^{1,2}(\Omega)\cap L^{6/5}(\Omega)$ imply
$\|\fe\|_{L^{3/2}(\Omega_R)}\le C$ and $\|\nabla
p_k\|_{L^{3/2}(\Omega_k)}\le C$, but we prefer to use here weaker
estimate~(\ref{epgr34}) which holds also in the general
case~$\ue_0\ne{\bf0}$ to make our arguments more universal.}
on~$k,R$. By construction, we have the weak convergences
$\ue_k\rightharpoonup \ve\mbox{ \ in \
}W^{1,2}_\loc(\overline\Omega),\quad p_k\rightharpoonup p\mbox{ \
in \ }W^{1,3/2}_\loc(\overline\Omega)$\footnote{The weak
convergence in $W_\loc^{1,2}(\overline\Omega)$ means the weak
convergence in $W^{1,2}(\Omega')$ for every bounded subdomain
$\Omega'\subset\Omega$.}.

In conclusion,  we can prove the following lemma.

\begin{lem}
\label{lem_Leray_symm} {\sl Assume that  $\Omega\subset\R^3$ is an
exterior  axially symmetric domain of type \eqref{Omega} with
$C^2$-smooth boundary $\partial\Omega$,  and ${\bf a}\in
W^{3/2,2}(\partial\Omega)$, ${\bf f}={\rm curl}\,{\bf b}\in
W^{1,2}(\Omega)$ are axially symmetric. If the assertion of
Theorem~\ref{kmpTh4.2} is false, then there exist $ \ve, p$ with
the following properties.

\medskip

(E) \ \,The axially symmetric functions $\ve\in H(\Omega)$, $p\in
D^{1,3/2}(\Omega)$ satisfy the Euler system \eqref{2.1} and
$\|p\|_{L^3(\Omega)}<\infty$.

\medskip

(E-NS) \ \,Condition (E) is satisfied and there exist a sequences
of axially symmetric functions $\ue_k \in W^{1,2}(\Omega_k)$,
$p_k\in {W^{1,3/2}(\Omega_k)}$,  $\Omega_k=\Omega\cap B_{R_k}$,
$R_k\to\infty$ as $k\to\infty$, and numbers $\nu_k\to0+$, such
that estimates (\ref{epgr33})--(\ref{epgr34}) hold, the pair
$(\ue_k,p_k)$ satisfies \eqref{NSk} with
$\fe_k=\frac{\nu_k^2}{\nu^2}\,{\bf f}$, \ ${\bf
a}_k=\frac{\nu_k}\nu\,{\bf A}$ (here ${\bf A}$ is solenoidal
extension of~${\bf a}$ from Lemma~\ref{kmpLem14.1}), and
\begin{equation}
\label{E-NS-ax} \|\nabla\ue_k\|_{L^2(\Omega_k)}\to1,\quad
\ue_k\rightharpoonup \ve\mbox{ \ in \
}W_\loc^{1,2}(\overline\Omega),\quad p_k\rightharpoonup p\mbox{ \
in \ }W_\loc^{1,3/2}(\overline\Omega),
\end{equation}
\begin{equation}\label{cont_e}
\begin{array}{l}
\displaystyle \nu= \int _\Omega ({\bf v}\cdot\nabla){\bf
v}\cdot\h\,dx
\end{array}\end{equation}
Moreover, $\ue_k\in W^{3,2}_{\loc}(\Omega)$ and $p_k\in
W^{2,2}_{\loc}(\Omega)$. }
\end{lem}

\pr We need to prove only the identity~(\ref{cont_e}), all other
properties are already established above. Choosing $\bfeta={\bf
w}_k$ in (\ref{NSESS'}) yields
 \begin{equation}
\label{NaaS}
\begin{array}{l}
\displaystyle \nu=\int _\Omega ({\bf w}_k\cdot\nabla){\bf
w}_k\cdot{\bf A} \,dx+ {1\over J_k}\int_\Omega{\bf A}\cdot\nabla{\bf w}_k\cdot{\bf A} \,dx\\[10pt]
\displaystyle -{1\over J_k}\int_\Omega\nabla{\bf A}
\cdot\nabla{\bf w}_k\,dx+\dfrac{1}{J_k}\int\limits_\Omega{\bf
f}\cdot{\bf w}_k\,dx.
\end{array}
\end{equation}
By the  H\"older inequality
$$
\begin{array}{l}
\displaystyle\left|\int_\Omega{\bf A}\cdot\nabla{\bf w}_k\cdot{\bf
A} \,dx\right| \le \|{\bf A}\|_{L^4(\Omega)}^2\|\nabla{\bf
w}_k\|_{L^2(\Omega)}
\le \|{\bf A}\|_{L^4(\Omega)}^2,\\[10pt]
\displaystyle\left|\int_\Omega\nabla{\bf A}\cdot\nabla{\bf
w}_k\,dx \right|  \le \|\nabla{\bf A}\|_{L^2(\Omega)}\|\nabla{\bf
w}_k\|_{L^2(\Omega)}\le \|\nabla{\bf A}\|_{L^2(\Omega)},\\[10pt]
\displaystyle\left|\int\limits_\Omega{\bf f}\cdot{\bf
w}_k\,dx\right|\le \|{\bf f}\|_{L^{6/5}(\Omega)}\|{\bf
w}_k\|_{L^6(\Omega)}.
\end{array}
$$
Note  that
$$
\begin{array}{r}
\displaystyle  \int_\Omega({\bf w}_k\cdot\nabla)\h\cdot {\bf
w}_k\,dx - \int_\Omega({\bf v} \cdot\nabla)\h\cdot{\bf v}\,dx =
\int_\Omega\bigl(({\bf w}_k-{\bf v})\cdot\nabla\bigr)
\h\cdot{\bf w}_k\,dx\\[12pt]
\displaystyle+\int_\Omega({\bf v}\cdot\nabla)\h\cdot({\bf
w}_k-{\bf v})\,dx ={\cal J}^{(1)}_k+{\cal J}^{(2)}_k\to 0
\quad{\mbox as}\;\; k\to \infty.
\end{array}
$$
Indeed, by Cauchy's  and Hardy's inequalities
$$
\begin{array}{l}
\displaystyle |{\cal J}^{(1)}_k|\le\left|\,\int_{\Omega_R}
\bigl(({\bf w}_k-{\bf v})\cdot\nabla\bigr)\h\cdot{\bf
w}_k\,dx\right| +{c\over R}\int_{\R^3\setminus B_R}
{|{\bf w}_k-{\bf v}||{\bf w}_k|\over |x|^2}\,dx\\[12pt]
\displaystyle\le \|{\bf w}_k-{\bf v}\|_{L^4(\Omega_R)} \|{\bf w}_k\|_{L^4(\Omega_R)}\|\nabla{\bf A}\|_{L^2(\Omega_R)}\\[12pt]
\displaystyle +{c_1\over R}\|\nabla({\bf w}_k-{\bf
v})\|_{L^2(\R^3\setminus B_R)} \|\nabla{\bf
w}_k\|_{L^2(\R^3\setminus B_R)}\le c(R)\|{\bf w}_k-{\bf
v}\|_{L^2(\Omega_R)}+{2c_1\over R},
\end{array}
$$
where $\Omega_R=B_R\cap\Omega$ and $c_1$ does not depend on $R$.
Hence, letting first $k\to+\infty$ and then $R\to+\infty$, we
conclude that
$$
\lim_{k\to+\infty}{\cal J}^{(1)}_k=0.
$$
Analogously, it can be shown that $\lim\limits_{k\to+\infty}{\cal
J}^{(2)}_k=0$. Consequently, we can let $k\to+\infty$ in
(\ref{NaaS}) to get the required identity~(\ref{cont_e}).\qed

\section{Euler equation} \label{eueq}

\setcounter{theo}{0} \setcounter{lem}{0}
\setcounter{lemr}{0}\setcounter{equation}{0}

In this section we assume that the assumptions~(E) (from
Lemma~\ref{lem_Leray_symm}) are satisfied. For definiteness, we
assume that

(SO)  $\Omega$ is the domain (\ref{Omega}) symmetric with respect
to the axis $O_{x_3}$ and
$$\Gamma_j\cap O_{x_3}\ne\emptyset,\quad j=1,\dots,M',$$
$$\Gamma_j\cap O_{x_3}=\emptyset,\quad j=M'+1,\dots,N.$$
(We allow also the cases $M'=N$ or $M'=0$, i.e.,
 when all components (resp., no components) of the
 boundary intersect the axis of symmetry.)

Denote $P_+=\{(0,{x_2},{x_3}):{x_2}>0,\ {x_3}\in\R\}$, \
${\D}=\Omega\cap P_+$, $\D_j=\Omega_j\cap P_+$. Of course, on
$P_+$ the coordinates $x_2,x_3$ coincides with coordinates $r,z$.

For a set $A\subset \R^3$ put $\breve{A}:=A\cap P_+$, and for
$B\subset P_+$ denote by $\widetilde B$ the set in $\R^3$ obtained
by rotation of $B$ around $O_z$-axis. From the conditions (SO${}$)
one can easily see that

(S${}_1$) ${\D}$ is an unbounded plane domain with Lipschitz
boundary. Moreover, $\breve\Gamma_j$ is a connected set for each
$j=1,\dots,N$. In other words, the family
$\{\breve\Gamma_j:j=1,\dots,N\}$ coincides with the family of all
connected components of the set $P_+\cap\partial{\D}$.

Then ${\bf v}$ and $p$ satisfy the following system in the plane
domain~${\mathcal D}$:
\begin{equation} \label{2.1'}
\left\{\begin{array}{rcl} \dfrac{\partial p}{\partial
z}+v_r\dfrac{\partial v_z}{\partial r}+ v_z\dfrac{\partial
v_z}{\partial z}=0,\\[6pt]
\dfrac{\partial p}{\partial
r}-\dfrac{(v_\theta)^2}r+v_r\dfrac{\partial v_r}{\partial r}+
v_z\dfrac{\partial v_r}{\partial z}=0,\\[6pt]
\dfrac{v_\theta v_r}r+v_r\dfrac{\partial v_\theta}{\partial r}+
v_z\dfrac{\partial v_\theta}{\partial z}=0,\\[6pt]
\dfrac{\partial (rv_r)}{\partial r}+\dfrac{\partial
(rv_z)}{\partial z}=0
\end{array}\right.
\end{equation}
(these equations are fulfilled for almost all $x\in{\mathcal
D}$\,).

We have the following integral estimates: ${\bf v}\in
W^{1,2}_\loc({\mathcal D})$,
\begin{equation}
\label{ax'3}\int_{{\mathcal D}}r|\nabla{\bf
v}(r,z)|^2\,drdz<\infty.
\end{equation}
\begin{equation}
\label{un_ax3}\int_{{\mathcal D}}r|{\bf v}(r,z)|^6\,drdz<\infty.
\end{equation}
 Also, the  inclusions $\nabla p\in
L^{3/2}(\Omega)$, $p\in L^3(\Omega)$ can be rewritten in the
following two-dimensional form:
\begin{equation}
\label{un_p2}\int_{{\mathcal D}}r|\nabla
p(r,z)|^{3/2}\,drdz<\infty,
\end{equation}
\begin{equation}
\label{p-bl3}\int_{{\mathcal D}}r|p(r,z)|^{3}\,drdz<\infty.
\end{equation}

The next statement  was proved in \cite[Lemma 4]{KaPi1} and in
\cite[Theorem 2.2]{Amick}.

\begin{theo}
\label{kmpTh2.3'} {\sl Let the conditions {\rm (E)} be fulfilled.
Then
\begin{equation} \label{bp2} \forall j\in\{1,\dots,N\} \ \exists\, \widehat p_j\in\R:\quad
p(x)\equiv \widehat p_j\quad\mbox{for }\Ha^2-\mbox{almost all }
x\in\Gamma_j.\end{equation} In particular, by axial symmetry,
\begin{equation}
\label{bp1}p(x)\equiv \widehat p_j\quad\mbox{for
}\Ha^1-\mbox{almost all } x\in \breve\Gamma_j.\end{equation} }
\end{theo}

\medskip

We shall prove (see Corollary~\ref{un_cor1}), that in the axially
symmetric case $\widehat p_1=\dots =\widehat p_{M'}=0$, i.e.,
these values coincide with the constant limit of the pressure
function $p$ at infinity. To formulate this result, we need some
preparation. Below without loss of generality we assume that the
functions $\ve,p$ are extended to the whole half-plane $P_+$  as
follows:
\begin{equation}
\label{axc10.10} \ve(x):=0,\quad x\in P_+\setminus\D,
\end{equation}
\begin{equation}
\label{axc110} p(x):=\widehat p_j,\ \, x\in P_+\cap\bar\D_j,\
j=1,\dots,N.
\end{equation}
Obviously,  the extended functions inherit the properties of the
previous ones. Namely, $\ve\in W^{1,2}_\loc(P_+)$, $p\in
W^{1,3/2}_\loc(P_+)$, and the Euler equations  (\ref{2.1'}) are
fulfilled almost everywhere in  $P_+$. Of course, for the
corresponding axial-symmetric functions of three variables we have
$\ve\in H(\R^3)$, $p\in D^{1,3/2}(\R^3)$,  and the Euler equations
(\ref{2.1}) are fulfilled almost everywhere in $\R^3$.

\begin{lem} \label{pressure} {\sl For almost all $r_0>0$
\begin{equation}
\label{un_p4}|p(r_0,z)|+|\ve(r_0,z)|\to 0\quad\mbox{as
}|z|\to\infty.
\end{equation}
}
\end{lem}

\pr This Lemma follows from the fact that for almost all straight
lines $L_{r_0}=\{(r,z)\in\R^2:r=r_0\}$ we have the inclusion
$|p(r_0,\cdot)|+|\ve(r_0,\cdot)|^2\in L^3(\R)\cap D^{1,3/2}(\R)$,
see (\ref{ax'3})-(\ref{p-bl3}). $\qed$

\medskip

The main result of this Section is a weak version of Bernoulli Law
for the Sobolev solution $({\bf v}, p)$ of Euler equations
\eqref{2.1'} (see Theorem~\ref{kmpTh2.2} below). To formulate and
to prove it, we need some preparation.

The last equality in (\ref{2.1'}) (which is fulfilled, after the
above extension agreement, see~(\ref{axc10.10})--(\ref{axc110}),
in the whole half-plane~$P_+$) implies the existence of a stream
function $\psi\in W^{2,2}_{\loc}(P_+)$ such that
\begin{equation}
\label{ax7}\frac{\partial\psi}{\partial
r}=-rv_z,\quad\frac{\partial\psi}{\partial z}=rv_r.
\end{equation}
By (\ref{ax7}), formula~(\ref{un_ax3}) can be rewritten in the
following form:
\begin{equation}
\label{ax3'}\int_{P_+}\frac{|\nabla\psi(r,z)|^6}{r^5}\,drdz<\infty.
\end{equation}
By Sobolev Embedding Theorem, $\psi\in C(P_+)$
(recall, that $P_+$ is an open half-plane, so here we
do not assert the continuity at the points of singularity
axis~$O_z$). By virtue of~(\ref{axc10.10}), we have $\nabla
\psi(x)=0$ for almost all $x\in{\mathcal D}_j$. Then
\begin{equation}
\label{axc10} \forall j\in\{1,\dots,N\}\ \ \exists\, \xi_j\in\R:
 \quad \psi(x)\equiv\xi_j\qquad \forall x\in P_+\cap\bar\D_j.
\end{equation}

Denote by $\Phi= p+\dfrac{|{\bf v}|^2}{2}$ the total head pressure
corresponding to the solution $({\bf v}, p)$.
From~(\ref{ax'3})--(\ref{p-bl3}) we get
\begin{equation}
\label{axc2}\int_{P_+}r|\Phi(r,z)|^3\,drdz+\int_{P_+}r|\nabla
\Phi(r,z)|^{3/2}\,drdz<\infty.
\end{equation}
By direct calculations one easily gets the identity
\begin{equation}
\label{2.2} v_r \Phi_r+ v_z\Phi_z=0
\end{equation}
for almost all $x\in P_+$. Identities
(\ref{axc10.10})--(\ref{axc110}) mean that
\begin{equation}
\label{axc11} \Phi(x)\equiv\widehat p_j\qquad \forall x\in
P_+\cap\bar\D_j,\ \,j=1,\dots,N.
\end{equation}

\medskip

\begin{theo}[Bernoulli Law for Sobolev solutions]
\label{kmpTh2.2} {\sl Let the conditions~{\rm (E)} be valid. Then
there exists a set $A_\ve\subset P_+$ with $\Ha^1(A_\ve)=0$, such
that for any compact connected\footnote{We understand the
connectedness in the sense of general topology.} set $K\subset
P_+$ the following property holds : if
\begin{equation}
\label{2.4} \psi\big|_{K}=\const,
\end{equation}
then
\begin{equation}
\label{2.5'}  \Phi(x_1)=\Phi(x_2) \quad\mbox{for
 all \,}x_1,x_2\in K\setminus A_{\bf v}.
\end{equation}
}
\end{theo}

Theorem~\ref{kmpTh2.2} was obtained for bounded plane domains
in~\cite[Theorem~1]{korob1} (see also~\cite{kpr} for detailed
proof). For the axially symmetric bounded domains the result was
proved in~\cite[Theorem~3.3]{kpr_a_arx}. The proof for exterior
axially symmetric domains is similar. To prove
Theorem~\ref{kmpTh2.2}, we have to overcome two obstacles. First
difficulty is the lack of the classical regularity, and here the
results of~\cite{korob} have a decisive role (according to these
results, almost all level sets of plane $W^{2,1}$-functions are
$C^1$-curves, see Section~\ref{SMS}). The second obstacle is the
set where $\nabla\psi(x)=0\ne\nabla\Phi(x)$, i.e., where
$v_r(x)=v_z(x)=0$, but $v_\theta(x)\ne 0$. Namely, if we do not
assume the boundary conditions~(\ref{2.1}${}_3$), then in general
even in smooth case (\ref{2.4}) does not imply~(\ref{2.5'}).  For
example, if $v_r=v_z=0$ in the whole domain, $v_\theta=r$, then
$\psi\equiv\const$ on the whole domain, while $\Phi=r^2\ne\const$.
Without boundary assumptions one can prove only the assertion
similar to Lemma~\ref{lemTh2.2unb} (see below). But
Lemma~\ref{lemTh2.2unb} together with boundary
conditions~(\ref{2.1}${}_3$) imply Theorem~\ref{kmpTh2.2}.

First, we prove some auxiliary results.

\begin{lem}
\label{pp} {\sl Let the conditions {\rm (E)} be fulfilled. Then
\begin{equation} \label{co-1}p\in
W^{2,1}_\loc(P_+).
\end{equation}
Moreover, for any $\e>0$  the inclusion
\begin{equation} \label{co1}p\in
D^{2,1}(P_{\e})
\end{equation}
holds, where $P_{\e}=\{(r,z)\in P_+:r>\e\}$.}
\end{lem}

\pr  Clearly, $p\in D^{1,3/2}(\R^3)$ is the weak solution to the
Poisson equation
\begin{equation} \label{r2.2}
\Delta p =-\nabla{\mathbf v}\cdot\nabla{\mathbf
v}^\top\quad\mbox{in }\R^3
\end{equation}
(recall, that after our agreement about extension of ${\mathbf
v}$ and $p$, see (\ref{axc10.10})--(\ref{axc110}), the Euler
equations (\ref{2.1}) are fulfilled in the whole $\R^3$). Let
$$
G(x)={1\over 4\pi}\int_{\Omega}{(\nabla{\mathbf
v}\cdot\nabla{\mathbf v}^\top)(y)\over |x-y|} dy.
$$

By the results of \cite{CLMS}, $\nabla{\mathbf
v}\cdot\nabla{\mathbf v}^\top$  belongs to the Hardy space
$\H^1(\R^3)$. Hence by Calder\'on--Zygmund theorem for Hardy's
spaces \cite{Stein} $G\in D^{2,1}(\R^3)\cap D^{1,3/2}(\R^3)$.
Consider the function $p_*=p-G$. By construction $p_*\in
D^{1,3/2}(\R^3)$ and $\Delta p_*=0$ in $\R^3$. In particular,
$\nabla p_*\in L^{3/2}(\R^3)$ is a harmonic (in the sense of
distributions) function. From the mean-value property it follows
that $p_*\equiv\const$. Consequently, $p\in D^{2,1}(\R^3)$. $\qed$

\bigskip

From inclusion~(\ref{co-1}) it follows that ~$\dfrac{\partial^2
p}{\partial r\partial z}\equiv \dfrac{\partial^2 p}{\partial
z\partial r}$ for almost all $x\in P_+$. Denote $Z=\{x\in
P_+:v_r(x)=v_z(x)=0\}$. Equations~(\ref{2.1'}) yield the equality
$$
\frac{\partial p}{\partial z}=0,\quad \frac{\partial p}{\partial
r}=\frac{(v_\theta)^2}{r}\quad\mbox{for almost all }x\in Z,
$$
and, using identity $\dfrac{\partial^2 p}{\partial r\partial
z}\equiv \dfrac{\partial^2 p}{\partial z\partial r}$ and equations
~(\ref{2.1'}), it is easy to deduce that
\begin{equation}
\label{2.2a}\frac{\partial \Phi}{\partial z}(x)=0\quad\mbox{for
almost all }x\in P_+\mbox{ such that }v_r(x)=v_z(x)=0.
\end{equation}

In the assertion below we collect the basic properties of Sobolev
functions applied to ${\bf v}$ and $\Phi$. Here and in the next
statements we assume that the conditions {\rm (E)} are fulfilled
and that all functions are extended to the whole half-plane $P_+$
(see \eqref{axc10.10}--\eqref{axc110},
\eqref{axc10}--\eqref{axc11}\,).

\begin{theo}[see, e.g., \cite{kpr_a_arx}]
\label{kmpTh2.1} {\sl There exists a set $A_{\bf v}\subset P_+$
such that:

 {\rm (i)}\quad $ \mathfrak{H}^1(A_{\bf v})=0$;

{\rm (ii)} For  all  $x\in P_+\setminus A_{\bf v}$
\begin{displaymath}
\lim\limits_{R\to 0}\dashint\nolimits_{B_R(x)}|{\bf v}(y)-{\bf
v}(x)|^2dy=\lim\limits_{R\to
0}\dashint\nolimits_{B_R(x)}|{\Phi}(y)-{\Phi}(x)|^{3/2}dy=0,
\end{displaymath}
\begin{displaymath}
\lim\limits_{R\to
0}\frac1R\int\nolimits_{B_R(x)}|\nabla\Phi(y)|^{3/2}dy=0.
\end{displaymath}
Moreover, the function $\psi$ is differentiable at $x\in
P_+\setminus A_{\bf v}$ and $\nabla\psi(x)=(-rv_z(x), rv_r(x))$;

{\rm  (iii) } For all $\varepsilon >0$ there exists an open set
$U\subset \mathbb{R}^2$ such that
$\mathfrak{H}^1_\infty(U)<\varepsilon$, $A_{\bf v}\subset U$, and
the functions ${\bf v}, \Phi$ are continuous in $P_+\setminus U$;

{\rm  (iv) } For each $x_0=(r_0,z_0)\in P_+\setminus A_{\bf v}$
and for any $\varepsilon>0$ the convergence
\begin{equation}
\label{z2}\lim\limits_{\rho\to0+}\frac1{2\rho}\Ha^1(E(x_0,\varepsilon,\rho))\to1
\end{equation}
holds, where
$$
E(x_0,\varepsilon,\rho):=\Big\{t\in(-\rho,\rho):
\int\limits_{r_0-\rho}^{r_0+\rho}\biggl|\frac{\partial
\Phi}{\partial r}(r,z_0+t)\biggr|\,dr+
\int\limits_{z_0-\rho}^{z_0+\rho}\biggl|\frac{\partial
\Phi}{\partial z}(r_0+t,z)\biggr|\,dz$$
$$+\sup\limits_{r\in[r_0-\rho,r_0+\rho]}|\Phi(r,z_0+t)-\Phi(x_0)|+
\sup\limits_{z\in[z_0-\rho,z_0+\rho]}|\Phi(r_0+t,z)-\Phi(x_0)|<\varepsilon\Big\}.$$

{\rm  (v) } Take any function $g\in C^1(\R^2)$ and a closed set
$F\subset P_+$ such that $\nabla g\ne 0$ on $F$. Then for almost
all $y\in g(F)$ and for all the connected components $K$ of the
set $F\cap g^{-1}(y)$ the equality $K\cap A_{\ve}=\emptyset$
holds, the restriction $\Phi|_K$ is an absolutely continuous
function, and formulas~(\ref{2.1'}),\,(\ref{2.2}) are fulfilled
$\Ha^1$-almost everywhere on~$K$.}
\end{theo}

Most of these properties are from \cite{evans}.  The property~(iv)
follows directly from the second convergence formula in~(ii). The
last property~(v) follows (by coordinate transformation, cf.
\cite[\S1.1.7]{maz'ya}) from the well-known fact that any function
$f\in W^{1,1}$ is absolutely continuous along almost all
coordinate lines. The same fact together with (\ref{2.2a}),
(\ref{2.2}) imply

\begin{lem}
\label{kmpLem2}{\sl Denote $L_{r_0}=\{ (r,z)\in
P_+:\, r=r_0 \}$.  Then for almost all $r_0>0$ the equality
$L_{r_0}\cap A_\ve=\emptyset$ holds. Moreover, $p(r_0,\cdot)$,
$\ve(r_0,\cdot)$ are absolutely continuous functions (locally) and
\begin{equation}
\label{z1}\frac{\partial \Phi}{\partial z}(r_0,z)=0\quad\mbox{for
almost all }z\in\R\mbox{ such that }v_r(r_0,z)=0.
\end{equation}}
\end{lem}

We need also the following assertion which was proved in
\cite{kpr_a_arx} (see Lemma 3.6).

\begin{lem} \label{lemTh5} {\sl Suppose for $r_0>0$
the assertion of Lemma~\ref{kmpLem2} is fulfilled, i.e., the
equality $L_{r_0}\cap A_\ve=\emptyset$ holds, $p(r_0,\cdot)$,
$\ve(r_0,\cdot)$ are absolutely continuous functions, and
formula~(\ref{z1}) is valid. Let $F\subset\R$ be a compact set
such that
\begin{equation}
\label{T1.1}\psi(r_0,z)\equiv\const\quad\mbox{for all }z\in F
\end{equation}
and
\begin{equation}
\label{T1.2}\Phi(r_0,\alpha)=\Phi(r_0,\beta) \quad\mbox{ for any
interval }(\alpha,\beta)\mbox{ adjoining }F
\end{equation}
(recall that $(\alpha,\beta)$ is called an interval adjoining~$F$
if $\alpha,\beta\in F$ and $(\alpha,\beta)\cap F=\emptyset$\,).
Then
\begin{equation}
\label{T1.3}\Phi(r_0,z)\equiv\const\quad\mbox{for all }z\in F.
\end{equation}
}
\end{lem}

The next lemma plays the key role in the proof of the~Bernoulli
Law.

\begin{lem}
\label{lemTh2.2unb} {\sl Let $P\subset P_+$ be a rectangle
$P:=\{(r,z):r\in[r_1,r_2],\ z\in[z_1,z_2]\}$, $r_1>0$, and
$K\subset P$ be a connected component of the set $\{x\in
P:\psi(x)=y_0\}$, where $y_0\in\R$.
 Then there
exists an absolute continuous function~$f:[r_1,r_2]\to\R$ such
that
\begin{equation}
\label{z1'u}\Phi(r,z)=f(r)\quad\mbox{for all }(r,z)\in K\setminus
A_\ve.
\end{equation}
Moreover, for each $r_0\in[r_1,r_2]$ if $f(r)\ne \const$ locally
in any neighborhood of $r_0$, then
\begin{equation}
\label{z2'u} (r_0,z)\in K\quad \forall z\in[z_1,z_2].
\end{equation}
In other words, if for $r_0\in[r_1,r_2]$ the inclusion
(\ref{z2'u}) is not valid, then there exists a neighborhood
$U(r_0)$ such that $f|_{U(r_0)\cap[r_1,r_2]}=\const$. }
\end{lem}

\begin{lemr}
\label{rem__int}{\rm Notice that the above lemma is valid without
assumption~(\ref{2.1}${}_3$) (that
$\ve|_{\partial\Omega}=0$\,). It is enough to suppose that
$\ve\in W^{1,2}(P)$, $p\in W^{1,3/2}(P)$ satisfy Euler
system~(\ref{2.1'}) almost everywhere in~$P$. }
\end{lemr}

{\bf Proof of Lemma~\ref{lemTh2.2unb}}  splits into six steps (see
below). Steps 1--4 coincides (in essential) with the corresponding
steps (having the same numbers) from the proof of Theorem~3.3 in
\cite{kpr_a_arx}, where we assumed additionally that
\begin{equation}\label{bc4} \ve(x)=0\quad \forall
x\in\partial P\setminus \{(r_1,z):z\in(z_1,z_2)\}.
\end{equation}
The new arguments appear in Steps~5 and 6, where we cannot simply repeat the arguments from \cite{kpr_a_arx} because of absence of boundary conditions~(\ref{bc4}).

\textsc{Step 1.} Put $P^\circ=\Int P=(r_1,r_2)\times(z_1,z_2)$.
Let $U_i$ be the connected components of the open set
$U=P^\circ\setminus K$. In \cite[Step~1 of the proof of
Theorem~3.3]{kpr_a_arx} it was proved that for any $U_i$ there
exists a number $\beta_i$ such that
 \begin{equation}\label{T4}
\Phi(x)=\beta_i \quad\mbox{for
 all \,}x\in K\cap\partial U_i\setminus A_\ve.
\end{equation}

\textsc{Step 2.} Using (\ref{T4}) and Lemma~\ref{lemTh5}, in
\cite[Step~2 of the proof of Theorem~3.3]{kpr_a_arx} it was proved
that for almost all $r \in(r_1,r_2)$ the identities
\begin{equation} \label{T5} \Phi(r ,z')=\Phi(r
,z'')\quad\forall\,(r ,z'),(r ,z'')\in K
\end{equation} hold.

\textsc{Step 3.} Denote by $\Pr_r E$ the projection of the set $E$
onto the $r$-axis. From \eqref{T4}, \eqref{T5} it follows that
\begin{equation}
\label{T5.1}
\Pr_rU_i\cap\Pr_rU_j\ne\emptyset\Rightarrow\beta_i=\beta_j.
\end{equation}
Denote $V=\Pr_r(P^\circ\setminus K)$,\footnote{By construction,
the set $[r_1,r_2]\setminus V$ coincides with the set of values
$r\in[r_1,r_2]$ such that the whole segment $\{r\}\times[z_1,z_2]$
lies in $K$. The set $V$ can be empty, but it can also coincide
with the whole $[r_1,r_2]$. For example, if $K$ is a circle, then
$V=[r_1,r_2]$. Further, $V=\emptyset$ iff $K=P$.}
 and let $(u_k,v_k)$ be the family of intervals adjoining
the compact set $[r_1,r_2]\setminus V$ (in other words,
$(u_k,v_k)$ are the maximal intervals of the open
set~$V=\bigcup_i\Pr_rU_i$\,). Then, repeating the
assertion~(\ref{T5.1}), we have
\begin{equation}
\label{T5.2}\Pr_rU_i\subset (u_k,v_k)\supset
\Pr_rU_j\Rightarrow\beta_i=\beta_j.
\end{equation}
In other words, for any interval $(u_k,v_k)$  adjoining the set
$[r_1,r_2]\setminus V$ there exists a constant $\beta(k)$ such
that
\begin{equation}
\label{T6''}  \Pr_rU_i\subset (u_k,v_k)\Rightarrow \bigl(
\Phi(x)=\beta(k) \quad\mbox{for
 all \,}x\in K\cap\partial U_i\setminus A_\ve\bigr).
\end{equation}
By construction, for any $r\in(u_k,v_k)$ there exists $z\in
(z_1,z_2)$ such that $(r,z)\in U_i$ for some component $U_i$ with
$\Pr_rU_i\subset (u_k,v_k)$. From this fact and from the
identities~(\ref{T6''}) and (\ref{T5}) it follows that
\begin{equation}
\label{T5'v} \Phi(r,z)=\beta(k)\quad\mbox{for
 almost all }r\in(u_k,v_k)\mbox{ and for any }(r,z)\in K.
\end{equation}

\textsc{Step 4.} It was proved in \cite[Step~4 of the proof of
Theorem~3.3]{kpr_a_arx} that the assertions (\ref{T6''}) and
(\ref{T5'v})
 imply
\begin{equation}
\label{T6'}  r \in[u_k,v_k]\Rightarrow \Phi(r ,z)=\beta(k)\quad
\forall(r ,z)\in K\setminus A_\ve .
\end{equation}
Note, that formulas (\ref{T6''}) and (\ref{T5'v}) in our proof
correspond to the formulas (3.38) and (3.39) from
\cite{kpr_a_arx}. Further, the role of the rectangle $P$ from
\cite[Step~4 of the proof of Theorem~3.3]{kpr_a_arx} is played in
our proof by the rectangle $[u_k,v_k]\times[z_1,z_2]$.

\textsc{Step 5.} By construction, for each $r
\in[r_1,r_2]\setminus V$ the whole segment $\{r \}\times[z_1,z_2]$
is contained in $K$. On this step we shall prove that $\Phi(x)$ is
constant along each of these segments, i.e.,
\begin{equation}
\label{T7}  \forall r \in[r_1,r_2]\setminus V\quad \exists \beta(r
) \quad \Phi(r ,z)=\beta(r )\quad \forall (r ,z)\in P\setminus
A_\ve.
\end{equation}
Indeed, if $r $ is an endpoint of an adjoining interval, i.e., if
$r =u_k$ or $r =v_k$, then  (\ref{T7}) immediately follows
from~(\ref{T6'}). Further, by Step~2 there exists a set $\tilde
E\subset [r_1,r_2]\setminus V$ such that $\Ha^1([r_1,r_2]\setminus
(V\cup \tilde E))=0$ and the assertion (\ref{T7}) is true for
any~$r \in\tilde E$. But for any $r \in [r_1,r_2]\setminus V$
there exists a sequence of points $r_\mu\to r$, $\mu=3,4,\dots$,
such that each $r_\mu$ is an endpoint of some adjoining interval
$(u_k,v_k)$ or $r_\mu\in\tilde E$, i.e., the assertion (\ref{T7})
is true for each~$r_\mu$. Finally, the assertion~(\ref{T7})
follows for $r =\lim\limits_{\mu\to\infty} r_\mu$ from the
continuity properties of~$\Phi(\cdot)$ (see
Theorem~\ref{kmpTh2.1} (ii)--(iv)).

\textsc{Step 6.} Define the target function~$f:[r_1,r_2]\to\R$ as
follows: $f(r)=\beta(k)$ for $r\in(u_k,v_k)$ (see Step~4) and
$f(r)=\beta(r)$ for $r\in[r_1,r_2]\setminus V$ (see Step~5). By
construction~(see (\ref{T6'}) and (\ref{T7})\,) we have the
identity~(\ref{z1'u}). Also by construction,
\begin{equation}
\label{T8}  f(r)\equiv \const\quad\mbox{on each adjoining interval
}(u_k,v_k).
\end{equation}
Take an arbitrary $z_0\in (z_1,z_2)$ such that the segment
$[r_1,r_2]\times\{z_0\}$ does not contain points from~$A_\ve$ and
$\Phi(\cdot,z_0)$ in an absolutely continuous function, i.e.,
$\Phi(\cdot,z_0)\in W^{1,1}([r_1,r_2])$. Then by construction,
$\Phi(r,z_0)=f(r)$ for each $r\in[r_1,r_2]\setminus V$, in
particular, $f(r)$ coincides with an absolute continuous function
on the set $[r_1,r_2]\setminus V$. The last fact together
with~(\ref{T8}) implies the absolute continuity of $f(\cdot)$ on
the whole interval~$[r_1,r_2]$. The Lemma is proved completely.
$\qed$

{\bf Proof of Theorem~\ref{kmpTh2.2}.} Suppose the conditions (E)
are fulfilled and $K\subset P_+$ is a compact connected set,
$\psi|_K\equiv\const$. Take the set~$A_\ve$ from
Theorem~\ref{kmpTh2.1}. Let $P$ be a rectangle
$P:=\{(r,z):r\in[r_1,r_2],\ z\in[z_1,z_2]\}$, $r_1>0$, such that
$K\subset P$, and
\begin{equation}
\label{T8'}  \Pr_rK=[r_1,r_2].
\end{equation}
We may assume without loss of generality that $K$ is a connected
component of the set $\{x\in P:\psi(x)=y_0\}$, where $y_0\in\R$.
Apply Lemma~\ref{lemTh2.2unb} to this situation, and take the
corresponding function~$f(r)$. Then the target identity
(\ref{2.5'}) is equivalent to the identity $f(r)\equiv\const$. We
prove this fact getting a contradiction.  Suppose the last
identity is false. Consider the nonempty compact set $\mathcal R
=\{r_0\in[r_1,r_2]: f(r)\ne\const$ in any neighborhood of $r_0\}$.
By assumption, $\mathcal R\ne\emptyset$, thus $f(\mathcal
R)=f([r_1,r_2])$ is an interval of positive length\footnote{Notice
that the set $\mathcal R$ itself may have empty interior, for
example, if $f$ is a Cantor staircase type  function,
i.e., if $f$ is constant on each interval $I_j$, and the union of these intervals is everywhere dense set, then $\mathcal R$ coincides with corresponding Cantor set and has empty interior.} Since $f$ is an absolute
continuous function, it has Luzin $N$-property, i.e., it maps sets
of measure zero into the sets of measure zero. In particular, the
measure of $\mathcal R$ must be positive, moreover, there exists a
set of positive measure $\mathcal R'\subset \mathcal R$ such that
$f(r)\ne0$ for each $r\in \mathcal R'$. So, by
Lemma~\ref{pressure} there exists $r_0\in \mathcal R'$
such that $L_{r_0}\cap A_\ve=\emptyset$, $$
|p(r_0,z)|+|\ve(r_0,z)|\to 0\quad\mbox{as }|z|\to\infty,
$$
and
\begin{equation}
\label{T9}  f(r_0)\ne 0.
\end{equation}

Take a sequence of rectangles
$P_m=[r_1,r_2]\times[z^m_{1},z^m_{2}]$ such that
$z^m_1\to-\infty$, $z^m_2\to+\infty$. Let $K_m$ be a connected
component of the level set $\{x\in P_m:\psi(x)=y_0\}$ containing
$K$. From (\ref{T8'}) it follows that
\begin{equation}
\label{T10}  \Pr_rK_m=[r_1,r_2].\end{equation} Apply
Lemma~\ref{lemTh2.2unb} to $K_m, P_m$, and take the corresponding
function~$f_m(r)$. By construction (see, e.g., \eqref{T8'},
(\ref{T10})\,) functions $f_m$ do not depend on $m$, i.e.,
$f(r)=f_m(r)$ for each $r\in[r_1,r_2]$ and for all $m=1,2,\dots$.
Then by the second assertion of Lemma~\ref{lemTh2.2unb}, the whole
segment $\{r_0\}\times[z^m_{1},z^m_{2}]$ is contained in $K_m$ for
each $m$, and $\Phi(r_0,z)\equiv f(r_0)$ for all $z\in
[z^m_{1},z^m_{2}]$. Passing to a limit as $m\to\infty$, we get
$\Phi(r_0,z)\equiv f(r_0)\ne0$ for all $z\in \R$. The last
identity contradicts convergence (\ref{un_p4}). The Theorem is
proved. $\qed$

\vspace{0.3cm}

For $\varepsilon>0$ and $R>0$ denote by $S_{\varepsilon,R}$ the
set $ S_{\varepsilon,R}=\{(r,z)\in P_+:r\ge\varepsilon,\
r^2+z^2=R^2\}$.

\begin{lem} \label{pressure-spheres} {\sl For any $\varepsilon>0$ there exists a sequence $\rho_j>0$,
$\rho_j\to+\infty$, such that
$S_{\varepsilon,\rho_j}\cap A_\ve=\emptyset$ and
\begin{equation}
\label{un_p_sp}\sup_{x\in S_{\varepsilon,\rho_j}}|\Phi(x)|\to
0\quad\mbox{as }j\to\infty.
\end{equation}
}
\end{lem}

\pr Fix $\varepsilon>0$. From the inclusion $\nabla\Phi\in
L^{3/2}(\R^3)$ it follows that there exists a sequence
$\rho_j\to\infty$ such that $S_{\varepsilon,\rho_j}\cap
A_\ve=\emptyset$ and
$$\int_{S_{\varepsilon,\rho_j}}|\nabla\Phi(x)|^{3/2}\,d\Ha^1\le\frac1{\rho_j}\quad\mbox{as
}j\to\infty.
$$
By H\"older inequality,
$$
\int_{S_{\varepsilon,\rho_j}}|\nabla\Phi(x)|\,d\Ha^1\le
\biggl(\int_{S_{\varepsilon,\rho_j}}|\nabla
\Phi|^{3/2}\,d\Ha^1\biggr)^{2/3} (\pi \rho_j)^{1/3}\le
\bigl(\frac{\pi}{\rho_j}\bigr)^{1/3},
$$
consequently, $$\diam\Phi(S_{\varepsilon,\rho_j})\to
0\quad\mbox{ as }j\to\infty,$$ that in virtue of
Lemma~\ref{pressure} implies the assertion of
Lemma~\ref{pressure-spheres}.  $\qed$

\medskip

One of the main results of this Section is the following.

\begin{theo}
\label{kmpTh2.3_un} {\sl Assume that conditions {\rm (E)} are
satisfied. Let $K_j$ be a sequence of continuums  with the
following properties: $K_j\subset\bar{\mathcal D}\setminus O_z$,
$\psi|_{K_j}=\const$, and
$\lim\limits_{j\to\infty}\inf\limits_{(r,z)\in K_j}r=0$, \
$\varliminf\limits_{j\to\infty}\sup\limits_{(r,z)\in K_j}r>0$.
Then $\Phi(K_j)\to 0$ as $j\to\infty$. Here we denote by $\Phi(K_j)$ the corresponding constant
$c_j\in\R$ such that $\Phi(x)=c_j$ for all $x\in K_j\setminus
A_\ve$ (see Theorem~\ref{kmpTh2.2}).}
\end{theo}

 \pr
 Let the assumptions of the Theorem be fulfilled. We shall use the Bernoulli law and the fact
that the axis $O_z$ is "almost" a stream line. More precisely,
$O_z$ is a singularity line for $\bf v$, $\psi$, $p$, but it can
be accurately approximated by usual stream lines (on which
$\Phi=\const$). Recall, that the functions $\ve,p,\Phi,\psi$ are
extended to the whole half-plane $P_+$ (see \eqref{axc10.10},
\eqref{axc110}, \eqref{axc11}, \eqref{axc10}), and the assertion
of the Bernoulli Law (Theorem~\ref{kmpTh2.2}) is true for these
extended functions.

 By Lemma~\ref{pressure}, there exists a constant $r_0>0$
such that the  straight line $L_{r_0}$ satisfies the assertion of Lemma~\ref{kmpLem2} and
\begin{equation}\label{un_bx}
\begin{array}{lcr}
L_{r_0}\cap A_\ve=\emptyset, \quad p(r_0,z)\to 0,\quad |{\bf
v}(r_0, z)|\to 0 \;\;{\rm as}\;\; |z|\to\infty,\\
\\
{\bf v}(r_0, \cdot)\in L^6(\R), \quad {\bf v}(r_0, \cdot)\in
D^{1,2}(\R)\subset C(\R),\\
\\
p(r_0, \cdot)\in D^{1,3/2}(\R)\subset C(\R),\\
\\
r_0<\varliminf\limits_{j\to\infty}\sup\limits_{(r,z)\in K_j}r.
\end{array}
\end{equation}

In particular, the function $\Phi(r_0, \cdot)$ is continuous on
$\R$. In virtue of the last inequality, we can assume without loss of
generality that \begin{equation} \label{nun1} K_j\cap
L_{r_0}\ne\emptyset\quad\forall j\in\N.\end{equation} Suppose that
the assertion of Theorem~\ref{kmpTh2.3_un} is false, i.e.,
\begin{equation} \label{nun2} \Phi(K_j)\to p_*\ne0\quad\mbox{ as }j\to\infty.\end{equation}
Then by (\ref{un_bx}${}_1$) there exists a constant $C_1>0$ such
that
\begin{equation} \label{nun3} \sup\limits_{j\in\N,\ (r_0,z)\in K_j}|z|\le C_1,\end{equation}
consequently,
\begin{equation} \label{nun4} K_j\cap
\{(r_0,z):z\in[-C_1,C_1]\}\ne\emptyset\quad\forall
j\in\N.
\end{equation}

Take a sequence of numbers $\tilde z_i\in(C_1,\infty)$ with
$\tilde z_i<\tilde z_{i+1}\to
+\infty\quad\mbox{as }i\to\infty.
$
Since, by estimate~(\ref{ax3'}),
$$\int\limits_{-\tilde z_i}^{\tilde z_i}\int\limits_{0}^1\frac{|\nabla\psi(r,z)|^6}{r^5}dr\,dz<\infty,$$
we conclude that there exists  a~sequence of numbers $r_{ik}\to0+$ such that
$$\int\limits_{-\tilde z_i}^{\tilde z_i}\frac{|\nabla\psi(r_{ik},z)|^6}{r_{ik}^5}dz<\frac1{r_{ik}},$$ that is,
$$\int\limits_{-\tilde z_i}^{\tilde z_i}|\nabla\psi(r_{ik},z)|^6\,dz<r^4_{ik}\to0\qquad\mbox{ as }k\to\infty.$$
In particular, by H\"older inequality
$$\int\limits_{-\tilde z_i}^{\tilde z_i}|\nabla\psi(r_{ik},z)|\,dz<\frac1i$$
for sufficiently large~$k$. This implies  the existence
of sequences
$r_i=r_{ik_i}\in (0,r_0)$ with the following
properties
\begin{equation}
\label{un_b7}r_i>r_{i+1}\to 0\quad\mbox{as }i\to\infty,
\end{equation}
\begin{equation}
\label{un_b8}\int_{-\tilde z_i}^{\tilde
z_i}|\nabla\psi(r_i,z)|\,dz\to 0\quad\mbox{as }i\to\infty.
\end{equation}
Further, by the same estimate~(\ref{ax3'}) and similar arguments,
there exists a sequence
$z_i\in[\tilde z_i-1,\tilde z_i]$ such that
\begin{equation}
\label{un_b8'}
\int_{r_i}^{r_0}|\nabla\psi(r,z_i)|\,dr+\int_{r_i}^{r_0}|\nabla\psi(r,-z_i)|\,dr\to
0\quad\mbox{as }i\to\infty.
\end{equation}
Thus,
\begin{equation}
\label{un_b8''}\int_{-z_i}^{z_i}|\nabla\psi(r_i,z)|\,dz+
\int_{r_i}^{r_0}|\nabla\psi(r,z_i)|\,dr+\int_{r_i}^{r_0}|\nabla\psi(r,-z_i)|\,dr\to
0
\end{equation}
as $i\to\infty$. Let $P_i$ be the~rectangle
$P_i=[r_i,r_0]\times[-z_i,z_i]$. Denote by $A_i, B_i$ the points
$A_i=(r_0,-z_i)\in L_{r_0}\cap\partial P_i$, $B_i=(r_0,z_i)\in
L_{r_0}\cap\partial P_i$. Denote by
$[A_i,B_i]=\{(r_0,z):z\in[-z_i,z_i]\}$ the closed segment and by
$]A_i,B_i[=[A_i,B_i]\setminus\{A_i,B_i\}$ the corresponding open
one. Put $T_i=(\partial P_i)\setminus\,]A_i,B_i[$. By construction
(see~(\ref{un_b8''})\,)
\begin{equation} \label{un_ax11'}
\diam\psi(T_i)\to 0\quad\mbox{as }i\to \infty.
\end{equation}

Take $R_*>\sqrt{r_0^2+C^2_1}$ such that
\begin{equation} \label{nun12}S_{r_0,R_*}\cap
K_j=\emptyset\quad\forall j\in\N.
\end{equation}
The existence of such $R_*$ follows from
Lemma~\ref{pressure-spheres} and (\ref{nun2})\,).
Indeed, Lemma~\ref{pressure-spheres}  gives us $\sup|\Phi(S_{r_0,R_k})|\to 0$ for some $R_k\to\infty$,
and (\ref{nun2}) means that~$\Phi(K_j)\to p_*\ne0$.

For $x\in P_i$ denote by $K^i_x$ the connected component of the
level set $\{y\in P_i:\psi(y)=\psi(x)\}$ containing~$x$. Put
$F_i=\{z\in [-z_i,z_i]:K^i_{(r_0,z)}\cap T_i\ne\emptyset\}$. Then
for each $i\in\N$ there exists an index $j(i)\ge i$ such that
\begin{equation} \label{nun13}\forall j\ge j(i) \quad
\{(r_0,z)\in K_j: z\in F_i\}\ne\emptyset.
\end{equation}
Indeed, the connected set $K_j$ intersect $L_{r_0}$ and $L_{r_i}$
for sufficiently large~$j$, moreover, by~(\ref{nun3}) \ $K_j\cap
L_{r_0}\subset \{r_0\}\times[-C_1,C_1]\subset[A_i,B_i]\subset
P_i$. From the last assertions and~(\ref{nun12}), by obvious
topological reasons, we derive the existence of
a~connected set $K^i_j\subset
K_j\cap P_i$ which intersect both lines~$L_{r_0}$ and $L_{r_i}$.
This means validity of~(\ref{nun13}).

Now take a point   $z^i_{j(i)}\in F_i$ such that
$(r_0,z^i_{j(i)})\in K_j$. Since the sequence  of points
$z^i_{j(i)}$ is bounded (see~(\ref{nun3})\,), we may assume
without loss of generality that
$$z^i_{j(i)}\to z_*\quad\mbox{ as }i\to\infty.$$ Then
$\psi(K_{j(i)})\to\psi(r_0,z_*)$ as $i\to\infty$. Denote
$\xi_*=\psi(r_0,z_*)$. Since by construction $K_{j(i)}\cap
T_i\ne\emptyset$ and the convergence~(\ref{un_ax11'}) holds, we
have
$$\sup\limits_{x\in T_i}|\psi(x)-\xi_*|\to 0\quad\mbox{as }i\to
\infty.
$$

By construction we have also the following properties of
sets~$F_i$.

($\mbox{I}_\sim$) \ $F_i$ is a compact set, $\pm z_i\in F_i$.

Indeed, the set $F_i\subset[z_i,z_0]$ is closed because of the
following reason: if $F_i\ni z_\mu\to z$, then there exists a
subsequence $z_{\mu_k}$ such that the components
$K^i_{(r_0,z_{\mu_k})}$ converge with respect to the Hausdorff
distance\footnote{The Hausdorff
distance $d_H$ between two compact sets $A,B\subset\R^n$ is
defined as follows: $d_H(A,B) = \max\bigl(\sup\limits_{a\in A}
\dist(a,B), \sup\limits_{b\in B} \dist(b,A)\bigr)$ (see, e.g., \S~7.3.1 in \cite{metric}). By Blaschke selection theorem [ibid], for any uniformly bounded
sequence of compact sets $A_i\subset\R^n$ there exists a
subsequence $A_{i_j}$ which converges to some compact set $A_0$
with respect to the Hausdorff distance. Of course, if all $A_i$
are compact connected sets, then the limit set $A_0$ is also connected.} to some set $K$. Of course, $K\ni
(r_0,z)$ is a compact connected set, $\psi|_K=\const$, and $K\cap
T_i\ne\emptyset$ (since by construction $K^i_{(r_0,z_{\mu})}\cap
T_i\ne\emptyset$). Therefore, $K\subset K^i_{(r_0,z)}$ (see the definition of the sets~$K^i_x$),
hence $z\in F_i$.

Put $U_i=]-z_i,z_i[\,\setminus F_i$. For $z\in U_i$ \,let\,
$(\alpha_i(z),\beta_i(z))$ be the maximal open interval from $U_i$
containing~$z$. Of course, $\alpha_i(z),\beta_i(z)\in F_i$. The
next two properties are evident.

($\mbox{II}_\sim$) \ $U_i$ is an open set, \ \ $U_i\subset
U_{i+1}$.

($\mbox{III}_\sim$) \ $\sup\limits_{x\in
T_i}|\psi(x)-\xi_*|\ge\sup\limits_{z\in
F_i}|\psi(r_0,z)-\xi_*|\to0\quad\mbox{as }i\to\infty$.

($\mbox{IV}_\sim$) \ $\forall z\in U_i$ $\exists$ a compact
connected set $K\subset P_i$ such that $(r_0,\alpha_i(z))\in K$,
$(r_0,\beta_i(z))\in K$ and $\psi|_K=\const$.

\medskip

Indeed, if the components $K^i_{(r_0,\alpha_i(z))}$,
$K^i_{(r_0,\beta_i(z))}$ do not coincide, then, by results
of~\cite{Kronrod}, there exists a compact connected set $K'\subset
P_i$, $\psi|_{K'}=\const$, which separates them, i.e., points
$(r_0,\alpha_i(z))$, $(r_0,\beta_i(z))$ lie in the different
connected components of the set $P_i\setminus K'$. Then, by
topological obviousness, $K'\cap T_i\ne\emptyset$ and
$K'\cap\{(r_0,z):z\in (\alpha_i(z),\beta_i(z))\}\ne\emptyset$. But
the last formulas contradict the condition
$(\alpha_i(z),\beta_i(z))\cap F_i=\emptyset$ and the definition of
$F_i$. Thus the property ($\mbox{IV}_\sim$) is proved.

The property~($\mbox{IV}_\sim$) together with the Bernoulli Law
(Theorem~\ref{kmpTh2.2}) imply the following identity:

\medskip
($\mbox{V}_\sim$) \ $\forall z\in U_i$ \
$\Phi(r_0,\alpha_i(z))=\Phi(r_0,\beta_i(z))$.

\medskip
Put $U=\bigcup\limits_iU_i$, $F=\R\setminus U$. Then we have

\medskip
($\mbox{VI}_\sim$) \ $F$ is a closed set, $z_*\in F$, $U$ is an
open set.

\medskip
For $z\in U$ put $\alpha(z)=\lim\limits_{i\to\infty}\alpha_i(z)$,
$\beta(z)=\lim\limits_{i\to\infty}\beta_i(z)$. Notice that the
limits exist since  the functions $\alpha_i(z),\beta_i(z)$ are
monotone in virtue of ($\mbox{II}_\sim$). Moreover, if $z<z_*$,
then $\beta(z)$ is finite because of inequalities $\beta_i(z)\le
z_*$. Analogously, if $U \ni z>z_*$, then $\alpha(z)\in[z_*,z)$.
By construction, $(\alpha(z),\beta(z))\subset U$, and, if
$\alpha(z)$ or $\beta(z)$ is finite, then it belongs to $F$.

From ($\mbox{III}_\sim$), ($\mbox{V}_\sim$) and continuity of
$\psi$ and $\Phi(r_0,\cdot)$ we have

\medskip
($\mbox{VII}_\sim$) \ $\forall z\in F$ \ $\psi(r_0,z)=\xi_*$.

($\mbox{VIII}_\sim$) \ $\forall z\in U$ \ if both values
$\alpha(z)$ and $\beta(z)$ are finite, then
$\Phi(r_0,\alpha(z))=\Phi(r_0,\beta(z))$.

\medskip
Then Lemma~\ref{lemTh5} yields
\begin{equation}
\label{un_ax13} \forall z\in F\quad \Phi(r_0,z)=\Phi(r_0,z_*),
\end{equation}
and from (\ref{nun2}), (\ref{un_ax13}) and the choice of
$(r_0,z_*)$ we deduce that
\begin{equation}
\label{un_ax14} \forall z\in F\quad \Phi(r_0,z)=p_*.
\end{equation}
Now to finish the proof of the Theorem, i.e., to receive a~desired
contradiction, we need to deduce the identity
\begin{equation}
\label{un_ax15} p_*=0
\end{equation}
(it will contradict the assumption~(\ref{nun2})\,). For this
purpose, consider two possible cases.

(i) Suppose that $F$ is an unbounded set. Then by \eqref{un_bx} \
$\lim\limits_{F\ni z\to\pm\infty}\Phi(r_0,z)=0$, and the target
equality~(\ref{un_ax15}) follows from~(\ref{un_ax14}).

(ii) Suppose that $F$ is bounded. Then there exists $z\in U$,
$z<z_*$, such that $\alpha(z)=-\infty$. By definition,
$\alpha_i(z)\to-\infty$ as $i\to\infty$. By \eqref{un_bx}, \
$\lim\limits_{i\to\infty}\Phi(r_0,\alpha_i(z))=0$. On the other
hand, by ($\mbox{V}_\sim$) \ we have
$$\lim\limits_{i\to\infty}\Phi(r_0,\alpha_i(z))=\lim\limits_{i\to\infty}\Phi(r_0,\beta_i(z))=
\Phi(r_0,\beta(z))=p_*$$ (the last two equalities follow
from~(\ref{un_ax14}) and from the finiteness of $\beta(z)\in
F\cap(-\infty,z_*]$). Thus, the equality (\ref{un_ax15}) is
proved, but it contradicts the assumption~(\ref{nun2}). The
obtained contradiction finishes the proof of the Theorem. $\qed$

\medskip

\begin{lemA}
\label{un_cor1} {\sl Assume that conditions {\rm (E)} are
satisfied. Then $\Phi|_{\Gamma_j}\equiv0$ whenever $\Gamma_j\cap
O_z\ne\emptyset$, i.e., $$\widehat p_1 =\dots= \widehat p_{M'} =
0,$$ where $\widehat p_j$ are the constants from
Theorem~\ref{kmpTh2.3'}. }
\end{lemA}

\section{Obtaining a contradiction}
\label{contrad_s}
\setcounter{theo}{0} \setcounter{lem}{0}
\setcounter{lemr}{0}\setcounter{equation}{0}

From now on we assume that  assumptions (E-NS) (see
Lemma~\ref{lem_Leray_symm}) are satisfied. Our goal is to prove
that they   lead to a contradiction. This implies the validity of
Theorem~\ref{kmpTh4.2}.

First, we introduce the main idea of the proof (which is taken
from~\cite{kpr_a_ann}) in a heuristic way. It is well known that
every $\Phi_k=p_k+\frac12|\ue_k|^2$ satisfies the linear elliptic
equation
\begin{equation}\label{cle_laps}
\Delta\Phi_k=\omega_k^2+\frac1{\nu_k}\div(\Phi_k\ue_k)-\frac{1}{\nu_k}{\bf
f}_k\cdot\ue_k,
\end{equation}
where $\ov_k=\curl \ue_k$ and $\omega_k(x)=|\ov_k(x)|$.
 If $\fe_k=0$,  then by Hopf's maximum principle, in
a~subdomain $\Omega'\Subset\Omega$ with $C^2$-- smooth
boundary~$\partial\Omega'$ the maximum of $\Phi_k$ is attained at
the boundary $\partial\Omega'$, and if $x_*\in\partial\Omega'$ is
a maximum point, then  the normal derivative of $\Phi_k$ at $x_*$
is strictly positive. It is not sufficient to apply this property
directly. Instead  we will use  some ''integral analogs'' that
lead to a contradiction by using the Coarea formula (see
Lemmas~\ref{ax-lkr11}--\ref{real_an}). For sufficiently large $k$
we construct a~set $E_k\subset \Omega$ (see the~proof of
Lemma~\ref{ax-lkr12}) consisting of level sets of $\Phi_k$ such
that $E_k$ separates the boundary components $\Gamma_j$ where
$\Phi\ne0$ from the boundary components $\Gamma_i$ where $\Phi=0$
and from infinity. On the one hand, the area of each of these
level sets is bounded from below (since they separate the boundary
components), and by the Coarea formula this implies  the estimate
from below for $\int_{E_k}|\nabla \Phi_k|$ (see the proof of
Lemma~\ref{ax-lkr12}). On the other hand, elliptic equation
(\ref{cle_laps}) for $\Phi_k$, the convergence $\fe_k\to 0$, and
boundary conditions~(\ref{NSk}${}_3$) allow us
 to estimate $\int_{E_k}|\nabla \Phi_k|^2$ from above (see Lemma~\ref{ax-lkr11}),
 and this   asymptotically contradicts
the previous estimate.  (We use also isoperimetric inequality, see the
proof of Lemma~\ref{ax-lkr12}, and some elementary Lemmas from
real analysis, see Appendix).

Recall, that by assumptions~(SO)
\begin{equation}
\label{gamj}
\begin{array}{l}
\Gamma_j\cap O_z\ne\emptyset, \quad j=1,\dots, M',
\\
\Gamma_j\cap O_z=\emptyset, \quad j=M'+1,\dots, N.
\end{array}
\end{equation}
Consider the constants $\widehat p_j$ from Theorem~\ref{kmpTh2.3'}
(see also Corollary~\ref{un_cor1}). We need the following fact.
\begin{lem}
\label{nui} {\sl The identity
 \begin{equation}
\label{pnui} -\nu =\sum_{j=M'+1}^N \widehat p_j\F_j
\end{equation}
holds. }
\end{lem}

\pr We calculate the integral $\int\limits_\Omega ({\bf v}\cdot\nabla){\bf v}\cdot {\bf A}\, dx$ in the equality~(\ref{cont_e}) by
using the Euler equations (\ref{2.1}${}_1$). In virtue
of assumptions~(E),
\begin{equation}
\label{pomt} \|p\|_{L^3(\Omega)}+\|\nabla
p\|_{L^{3/2}(\Omega)}<\infty,
\end{equation}
and by Theorem~\ref{kmpTh2.3'} and Corollary~\ref{un_cor1} we have
$$ p(x)|_{\Gamma_j}\equiv 0,\;\;  j=1,\ldots, M',
\quad\;p(x)|_{\Gamma_j}\equiv \widehat p_j, \quad j=M'+1,\ldots,
N.
$$
From the inclusions (\ref{pomt}) it is easy to deduce
\begin{equation}\label{1pub}\sup\limits_{m\in\N}\int\limits_{|x|=m} | p|^2 \,dS<\infty.
\end{equation}
Indeed, for balls $B_m=B(0,m)$ we have the uniform estimates
\begin{equation}\label{1s3}
\int\limits_{B_m\setminus \frac12B_m}{p}^2\,dx\le
\biggl(\int\limits_{B_m\setminus
\frac12B_m}{p}^3\,dx\biggr)^{\frac23}\cdot\biggl(\meas\bigl({B_m\setminus
\frac12B_m}\bigr)\biggr)^{\frac13}\le C_1m,
\end{equation}
\begin{equation}\label{1s4}
\int\limits_{B_m\setminus \frac12B_m}|{p}\nabla{p}|\,dx\le
\biggl(\int\limits_{B_m\setminus
\frac12B_m}{p}^3\,dS\biggr)^{\frac13}\cdot\biggl(\int\limits_{B_m\setminus
\frac12B_m}|\nabla{p}|^{3/2}\,dS\biggr)^{\frac23}\le C_2.
\end{equation}
Denote
$\sigma_1=\min\limits_{R\in[\frac12m,m]}
\int_{S_R}{p}^2\,dS$, \ $\sigma_2=
\max\limits_{R\in[\frac12m,m]}\int_{S_R}{p}^2\,dS-\sigma_1$.
Then
\begin{equation}\label{1sp1}
\sigma_1\le\frac2m
\int\limits_{B_m\setminus
\frac12B_m}{p}^2\,dx\le
2C_1.
\end{equation}
Analogously, since
\begin{equation}\label{1s1}
\biggl(\int\limits_{S_R}{p}^2\,dS\biggr)'_R=\frac2R\int\limits_{S_{R}}{p}^2\,dS
+2\int\limits_{S_R}{p}\nabla{p}\cdot\n\,dS,
\end{equation}
by (\ref{1s3})--(\ref{1s4}) we have
\begin{equation}\label{1s2}
\sigma_2\le\frac4{m}\int\limits_{B_m\setminus
\frac12B_m}{p}^2\,dx+2\int\limits_{B_m\setminus
\frac12B_m}|{p}\nabla{p}|\,dx\le 4C_1+2C_2.
\end{equation}
Because of inequality
$\int\limits_{|x|=m} | p|^2 \,dS\le \sigma_1+\sigma_2$, we have proved the required uniform boundedness of these integrals.

Hence
$$\int\limits_{|x|=m}
| p||\h|\,dS\le\biggl(\int\limits_{|x|=m} |
p|^2dS\biggr)^{\frac12}\biggl(\int\limits_{|x|=m}
|\h|^{2}dS\biggr)^{\frac12}\to0\quad\mbox{as }m\to\infty.$$
Multiplying equations (\ref{2.1}${}_1$) by $\h$, integrating by
parts in $\Omega_m=\{x\in\Omega:|x|<m\}$ and passing to a limit as
$m\to\infty$, we obtain
\begin{equation}
\label{pcf}
\begin{array}{l}
\int\limits_{\Omega}({\bf v}\cdot\nabla){\bf v}\cdot\h\,dx=
-\int\limits_{\Omega}\nabla
 p\cdot\h\,dx=-\lim\limits_{m\to\infty}\int\limits_{\Omega_m}\div(p\h)\,dx\\
= -\sum\limits_{j=M+1}^N\widehat
p_j\F_j+\lim\limits_{m\to\infty}\int\limits_{|x|=m} p\h\cdot{\bf
n} \,dS
\\
= -\sum\limits_{j=M+1}^N\widehat p_j\F_j.
\end{array}
\end{equation}
The required
equality~(\ref{pnui}) follows from the last identity and (\ref{cont_e}).\qed

\medskip

If $\widehat p_{M+1}=\ldots=\widehat p_{N}=0$, we get $\int\limits_{\Omega}({\bf v}\cdot\nabla){\bf v}\cdot\h\,dx=0$. However, this contradicts the equality ~(\ref{pnui}). Thus, there is $\widehat p_j\neq 0$ for some
$j\in\{M+1,\dots,N\}$.

\medskip

Further we consider separately three possible cases.

(a) The maximum of $\Phi$ is attained at infinity, i.e., it is
zero:
\begin{equation}\label{as-prev1}
0=\esssup\limits_{x\in\Omega}\Phi(x).
\end{equation}

(b) The maximum of $\Phi$ is attained on a boundary component
which does not intersect the symmetry axis:
\begin{equation}\label{as1-axxx}
0<\widehat p_N=\max\limits_{j=M'+1,\dots,N}\widehat
p_j=\esssup\limits_{x\in\Omega}\Phi(x),
\end{equation}

(c) The  maximum of $\Phi$ is not zero and it is not
attained\footnote{The case
$\esssup\limits_{x\in\Omega}\Phi(x)=+\infty$ is not excluded. } on
$\partial\Omega$:
\begin{equation}\label{as-prev-id-ax} \max\limits_{j=M'+1,\dots,N}\widehat
p_j<\esssup\limits_{x\in\Omega}\Phi(x)>0.
\end{equation}

\subsection{The case $\esssup\limits_{x\in\Omega}\Phi(x)=0$.}
\label{EPcontr-axx}

Let us consider  case (\ref{as-prev1}). By
Corollary~\ref{un_cor1},
\begin{equation}\label{as-prev1'}
\widehat p_1=\dots=\widehat
p_{M'}=\esssup\limits_{x\in\Omega}\Phi(x)=0.
\end{equation}Since the identity
$\widehat p_{M'+1}=\dots=\widehat p_N=0$ is impossible, we have that $\widehat p_j<0$ for some
$j\in\{M'+1,\dots,N\}$. Change (if necessary) the numbering of the
boundary components $\Gamma_{M'+1}$, \dots, $\Gamma_{N}$ so that
\begin{equation}\label{as01-ax}
\widehat p_j=0, \quad j=0,\dots,M,\quad M\ge M', \end{equation}
\begin{equation}\label{as1-ax'}
\widehat p_{M+1}=\dots=\widehat p_N<0.
\end{equation}

Recall, that in our notation $P_+=\{(0,{x_2},{x_3}):{x_2}>0,\
{x_3}\in\R\}$, \ ${\D}=\Omega\cap P_+$. Of course, on $P_+$ the
coordinates $x_2,x_3$ coincides with coordinates $r,z$, and
$O_z=O_{x_3}$ is a~symmetry axis of $\Omega$.
 For a set $A\subset
\R^3$ put $\breve{A}:=A\cap P_+$.

We receive
a~contradiction following the arguments of \cite{kpr_a_arx},
\cite{kpr_a_ann}. Take the positive constant
$\delta_{p}=-\sup\limits_{j=M+1,\dots,N}\Phi(\Gamma_j)$. Our first
goal is to separate the boundary components where $\Phi<0$ from
infinity and from the singularity axis~$O_z$ by level sets
of~$\Phi$ compactly supported in $\D$. More precisely, for any
$t\in(0,\delta_{p})$ and $j=M+1,\dots,N$ we construct
a~continuum $A_j(t)\Subset P_+$ with the following properties:

(i) The set $\breve\Gamma_j$ lies in a~bounded connected component
of the open set~$P_+\setminus A_j(t)$;

(ii) $\psi|_{A_j(t)}\equiv\const$, \ $\Phi(A_j(t))=-t$;

(iii) (monotonicity) If $0<t_1<t_2<\delta_{p}$, then $A_j(t_1)$
lies in the unbounded connected component of the set~$P_+\setminus
A_j(t_2)$ (in other words, the set $A_j(t_2)\cup\breve\Gamma_j$
lies in the bounded connected component of the set~$P_+\setminus
A_j(t_1)$, see Fig.1).
\begin{center}
\includegraphics[scale=0.4]{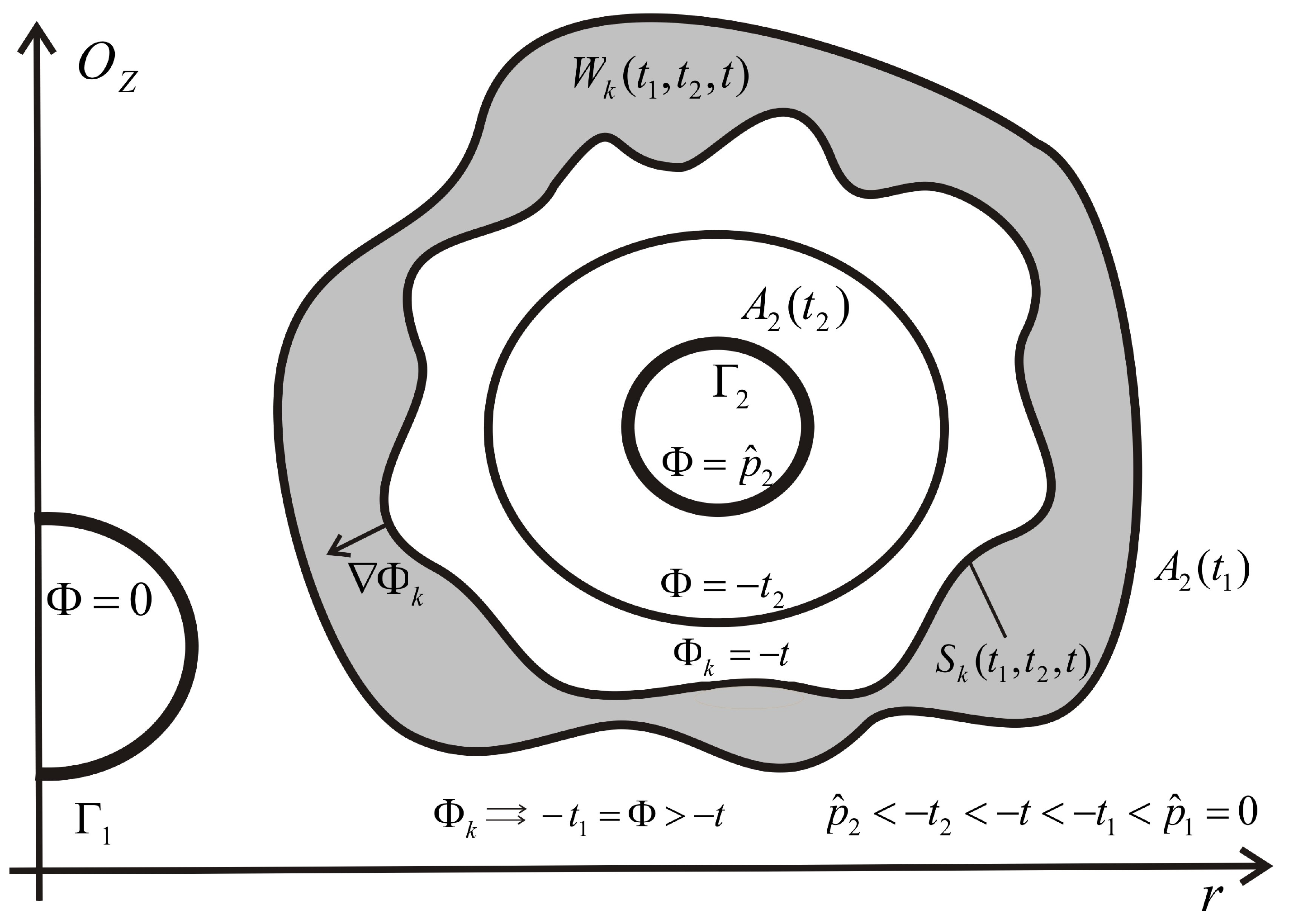}
\end{center}
\begin{center}
Fig. 1. {\sl The surface $S_k(t_1,t_2,t)$ for the case of $M=1$, $N=2$.}
\end{center}

For this construction, we shall use the results of
Subsection~\ref{Kronrod-s}. To do it, we have to consider the restrictions
of the stream function~$\psi$ to suitable compact subdomains of $P_+$.

Fix~$j\in\{M+1,\dots,N\}$. Take $r_j>0$ such that $L_{r_j}\cap
\Gamma_j\ne\emptyset$ and the conditions (\ref{un_bx}) are
satisfied with~$r_j$ instead of~$r_0$. In particular,
\begin{equation}\label{as-m1}
L_{r_j}\cap A_\ve=\emptyset,\quad \Phi(r_j,\cdot)\in
C(\R),\quad\mbox{ and }\quad\Phi(r_j,z)\to 0\mbox{ as
}z\to\infty.\end{equation} Take monotone sequences of positive
numbers~$\e_i\to 0+$, $\rho_i\to+\infty$ such that
\begin{eqnarray}\label{as-m2}
 \e_i<\inf\limits_{(r,z)\in\breve\Gamma_j}r\quad\mbox{ and }\quad S_{\e_i,\rho_i}\cap A_\ve=\emptyset\quad\forall i\in\N,\\
\label{as-m3}
\lim\limits_{i\to\infty}\sup|\Phi(S_{\e_i,\rho_i})|=0,
\end{eqnarray}
where $S_{\e,\rho}=\{(r,z)\in P_+:r\ge\e,\ \sqrt{r^2+z^2}=\rho\}$
(the existence of such sequences follow from
Lemma~\ref{pressure-spheres}).

By Remark~\ref{kmpRem1.2}, we can apply Kronrod's results to the
stream function~$\psi|_{\bar\D_{\e_i,\rho_i}}$, where
$\D_{\e,\rho}=\{(r,z)\in P_+:r>\e,\ \sqrt{r^2+z^2}<\rho\}$. Notice that $\breve\Gamma_j\subset D_{\e_i,\rho_i}$.
Accordingly, $T^i_{\psi}$ means the corresponding Kronrod tree for
the restriction $\psi|_{\bar\D_{\e_i,\rho_i}}$.

For any element~$C\in T^i_\psi$ with $C\setminus
A_\ve\ne\emptyset$  we can define the value $\Phi(C)$ as
$\Phi(C)=\Phi(x)$, where $x\in C\setminus A_\ve$. This definition
is correct because of the Bernoulli Law. (In particular, $\Phi(C)$ is well
defined if $\diam C>0$.)

Take $x_i=(r_j,z_i)\in S_{\e_i,\rho_i}$. Denote by $B_{x_i}$ and
$B^i_j$ the elements of $T^i_\psi$ with $x_i\in B_{x_i}$ and
$\breve\Gamma_j\subset B^i_j$. Consider the
arc~$[B^i_j,B_{x_i}]\subset T^i_\psi$. Recall, that, by
definition, a connected component $C$ of a level set
of~$\psi|_{\bar\D_{\e_i,\rho_i}}$ belongs to the
arc~$[B^i_j,B_{x_i}]$ iff $C=B_{x_i}$, or $C=B^i_j$, or $C$
separates $B_{x_i}$ from $B^i_j$ in $\bar\D_{\e_i,\rho_i}$, i.e.,
if $B_{x_i}$ and $B^i_j$ lie in different connected components
of~$\bar\D_{\e_i,\rho_i}\setminus C$. In particular, since
$B_{x_i}\cap L_{r_j}\ne\emptyset\ne B^i_j\cap L_{r_j}$, we have
\begin{equation}\label{f-m1}
C\cap L_{r_j}\ne\emptyset\quad\forall C\in
[B^i_j,B_{x_i}].\end{equation} Therefore, in view of equality
\begin{equation}\label{f-m2}
L_{r_j}\cap A_\ve=\emptyset\end{equation} the value $\Phi(C)$ is
well defined for all~$C\in [B^i_j,B_{x_i}]$. Moreover, we have
\begin{lem}
\label{l-m1} {\sl The restriction $\Phi|_{[B^i_j,B_{x_i}]}$ is a
continuous function. }
\end{lem}

\pr The assertion follows immediately\footnote{See also the proof of Lemma 3.5 in \cite{kpr_a_ann}.} from the assumptions
(\ref{f-m1})--(\ref{f-m2}), from the continuity of
$\Phi(r_j,\cdot)$, and from the definition of convergence
in~$T^i_\psi$ (see Subsection~\ref{Kronrod-s}\,). \qed

Define the natural order on~$[B^i_j,B_{x_i}]$. Namely, we say,
that $A< C$ for some different elements~$A,C\in[B^i_j,B_{x_i}]$
iff $C$ closer to~$B_{x_i}$ than $A$, i.e., if the sets $B_{x_i}$
and $C$ lie in the same connected component of the
set~$\bar\D_{\e_i,\rho_i}\setminus A$.

Put $$K_{\e_i}=\min\{C\in [B^i_j,B_{x_i}]: C\cap \partial
\D_{\e_i,\rho_i}\ne\emptyset\}.$$  The next assertion is an analog
of Lemma~4.6 from~\cite{kpr_a_ann}.

\begin{lem}
\label{l-m2} {\sl $\Phi(K_{\e_i})\to 0$ as $i\to\infty$. In
particular, $|\Phi(K_{\e_i})|<|\widehat p_j|=|\Phi(\Gamma_j)|$,
and, consequently, $B^i_j<K_{\e_i}$ for sufficiently large~$i$.}
\end{lem}

\pr By definition, $K_{\e_i}\cap \partial
\D_{\e_i,\rho_i}\ne\emptyset$. By construction, $\partial
\D_{\e_i,\rho_i}\subset S_{\e_i,\rho_i}\cup L_{\e_i}$. If
$K_{\e_i}\cap S_{\e_i,\rho_i}\ne\emptyset$, then the smallness of
$\Phi(K_{\e_i})$ follows immediately from the
assumption~(\ref{as-m3}). Now let
\begin{equation}\label{f-m4}
K_{\e_i}\cap L_{\e_i}\ne\emptyset.\end{equation} Recall, that,
by~(\ref{f-m1}), we have also
\begin{equation}\label{f-m5}
K_{\e_i}\cap L_{r_j}\ne\emptyset.\end{equation} Now the smallness
of $\Phi(K_{\e_i})$ follows from
Theorem~\ref{kmpTh2.3_un}. $\qed$

\medskip

In view of above Lemma we may assume without loss of generality
that \begin{equation}\label{f-m6}|\Phi(K_{\e_i})|<|\widehat
p_j|=|\Phi(\Gamma_j)|\quad\mbox{ and }\quad
B^i_j<K_{\e_i}\quad\forall i\in\N.
\end{equation}

By construction,
\begin{equation}\label{f-m7}C\cap\partial \D_{\e_i,\rho_i}=\emptyset
\quad\forall C\in[B^i_j,K_{\e_i}).
\end{equation}
In particular,
\begin{equation}\label{f-m8}B^i_j\cap\partial \D_{\e_i,\rho_i}=\emptyset
\quad\forall i\in\N.
\end{equation}
Therefore, really $B^i_j$ does not depend on~$i$, so we have
\begin{equation}\label{f-m9}B^i_j\equiv B_j
\quad\forall i\in\N
\end{equation}
for some continuum~$B_j$. Also, by equality
$$K_{\e_i}=\sup\{C\in [B_j,B_{x_i}]: C\cap \partial
\D_{\e_i,\rho_i}=\emptyset\}$$ and by inclusions
$\D_{\e_i,\rho_i}\Subset \D_{\e_{i+1},\rho_{i+1}}$ we have
\begin{equation}\label{f-m10}[B_j,K_{\e_i})\subset
[B_j,K_{\e_{i+1}})\quad\forall i\in\N,
\end{equation}
where, as usual,
$[B_j,K_{\e_i})=[K_{\e_i},B_j]\setminus\{K_{\e_i}\}$.
 Denote
$$[B_j,\infty)=\bigcup\limits_{i\in\N}[B_j,K_{\e_i}).$$
The set $[B_j,\infty)$ inherits the order and the topology from
the arcs~$[B_j,K_{\e_i})$. Obviously, $[B_j,\infty)$ is
homeomorphic to the ray~$[0,\infty)\subset\R$.

By construction, we have the following properties of the set
$[B_j,\infty)$.

($*1$) each element $C\in [B_j,\infty)$ is a continuum, $C\Subset
P_+$, and the set $\breve\Gamma_j$ lies in a~bounded connected
component of the open set~$P_+\setminus C$ for $C\ne B_j$.

($*2$) \ $\psi|_{C}\equiv\const$ for all $C\in [B_j,\infty)$.

($*3$) (monotonicity) \ If $C', C''\in [B_j,\infty)$ and $C'<C''$,
then $C''$ lies in the unbounded connected component of the
set~$P_+\setminus C'$; \ i.e.,  the set $C'\cup B_j$ lies
in the bounded connected component of the set~$P_+\setminus
C''$.

($*4$) (continuity) \ If $C_m\to C_0\in[B_j,\infty)$, then
$\sup\limits_{x\in C_m}\dist(x,C_0)\to0$ and
$\Phi(C_m)\to\Phi(C_0)$ as $m\to\infty$. In particular,
$\Phi|_{[B_j,\infty)}$ is a continuous function.

($*5$) (range of values) \ $\Phi(C)<0$ for every
$C\in[B_j,\infty)$. Moreover, if $[B_j,\infty)\ni C_m\to \infty$
as $m\to\infty$, then $\Phi(C_m)\to0$ as $m\to\infty$.

In the last property we use the following natural definition: for a
sequence $C_m\in [B_j,\infty)$ we say that $C_m\to\infty$ as
$m\to\infty$ if for any $i\in\N$ there exists $M(i)$ such that
$C_m\notin [B_j,K_{\e_{i}})$ for all $m>M(i)$.

We say that a set $\mathcal Z\subset [B_j,\infty)$ has $T$-measure
zero if $\Ha^1(\{\psi(C):C\in \mathcal Z\})=0$.

\begin{lem}
\label{regc-ax} {\sl For every $j=M+1,\dots,N$,  $T$-almost all
$C\in[B_j,\infty)$ are $C^1$-curves homeomorphic to the circle and $C\cap A_\ve=\emptyset$. Moreover,
there exists
a~subsequence~$\Phi_{k_l}$ such that $\Phi_{k_l}|_C$ converges to
$\Phi|_C$ uniformly $\Phi_{k_l}|_C\rightrightarrows\Phi|_C$ on $T$-almost all
$C\in[B_j,\infty)$.}
\end{lem}

Below we assume (without loss of generality) that the subsequence $\Phi_{k_l}$ coincides with the whole sequence $\Phi_{k}$.
\\

\pr  The first assertion of the lemma follows from
Theorem~\ref{kmpTh1.1}~(iii) and ~(\ref{f-m7}). The validity of
the second one for $T$-almost all
$C\in[B_j,\infty)$ was proved in~\cite[Lemma 3.3]{kpr}. \qed \\

Below we will call {\it regular } the cycles $C$ which satisfy
the~assertion of Lemma~\ref{regc-ax}.

Since $\diam C>0$ for every $C\in [B_j,\infty)$, by
\cite[Lemma~3.6]{kpr_a_ann} we obtain that the function
$\Phi|_{[B_j,\infty)}$ has the following analog of Luzin's
$N$-property.

\begin{lem}
\label{lkr7} {\sl For every $j=M+1,\dots,N$, if $\mathcal Z\subset
[B_j,\infty)$ has $T$-measure zero, then $\Ha^1(\{\Phi(C):C\in
\mathcal Z\})=0$.}
\end{lem}

From the last two assertions we get

\begin{lemA}
\label{regPhi-ax} {\sl For every $j=M+1,\dots,N$ and for almost
all $t\in(0,|\widehat p_j|)$ we have
$$\bigl(\,C\in[B_j,\infty)\mbox{\rm\ and }\Phi(C)=-t\,\bigr)\Rightarrow
C\mbox{\rm\ is a regular cycle}.$$}
\end{lemA}

Below we will say that a value~$t\in(0,\delta_p)$ is {\it
regular} if it satisfies the assertion of
Corollary~\ref{regPhi-ax}. Denote by~$\Ti$  the set of all
regular values. Then the set $(0,\delta_p)\setminus
\Ti$ has zero measure.

For $t\in (0,\delta_p)$ and $j\in\{M+1,\dots,N\}$ denote
$$A_j(t)=\max\{C\in[B_j,\infty):\Phi(C)=-t\}.$$
By construction, the function $A_j(t)$ is nonincreasing and
satisfies the properties (i)--(iii) from the beginning of this
subsection. Moreover, by definition of regular values we have the
following additional property:

(iv) If $t\in \Ti$, then $A_j(t)$ is a regular cycle\footnote{Some
of these cycles $A_j(t)$ could coincide, i.e., equalities of type
$A_{j_1}(t)=A_{j_2}(t)$ are possible (if Kronrod arcs
$[B_{j_1},\infty)$ and $[B_{j_2},\infty)$ have nontrivial
intersection), but this  a~priori possibility has no influence on
our arguments.}.

For $t\in\Ti$ denote by ${V}(t)$ the unbounded connected component
of the open set $\D\setminus\bigl(\cup_{j=M+1}^N A_j(t)\bigr)$.
Since $A_{j_1}(t)$ can not separate $A_{j_2}(t)$ from
infinity\footnote{Indeed, if $A_{j_2}(t)$ lies in a~bounded
component of $\D\setminus A_{j_1}(t)$, then by construction
$A_{j_1}(t)\in [B_{j_2},\infty)$ and $A_{j_1}(t)>A_{j_2}(t)$ with
respect to the~above defined order on~$[B_{j_2},\infty)$, but it
contradicts the definition of
$A_{j_1}(t)=\max\{C\in[B_{j_2},\infty):\Phi(C)=-t\}$.} for
$A_{j_1}(t)\ne A_{j_2}(t)$, we have
$$\D\cap\partial V(t)=A_{M+1}(t)\cup\dots\cup
A_N(t).$$  By construction, the sequence of domains ${V}(t)$ is
increasing, i.e., ${V}(t_1)\subset{V}(t_2)$ for $t_1<t_2$. Hence,
the sequence of sets $(\partial\D)\cap(\partial{V}(t))$ is
nondecreasing:
\begin{equation}\label{ref1}
(\partial\D)\cap \partial{V}(t_1)
\subseteqq(\partial\D)\cap\partial{V}(t_2)\quad\mbox{ if }t_1<t_2.
\end{equation}
Every set $(\partial\D)\cap\partial{V}(t)\setminus O_z$ consists
of several components $\breve\Gamma_l$ with $l\le M$ (since cycles
$\cup_{j=M+1}^N A_j(t)$ separate infinity from
$\breve\Gamma_{M+1},\dots,\breve\Gamma_N$, but not necessary from
other $\breve\Gamma_l$\,). Since there are only finitely many
components $\breve\Gamma_l$, using monotonicity~property
(\ref{ref1}) we conclude that for sufficiently small~$t$ the set
$(\partial\D)\cap(\partial{V}(t))$  is independent of $t$. So we
may assume, without loss of generality, that
$(\partial\D)\cap(\partial{V}(t))\setminus
O_z=\breve\Gamma_{1}\cup\dots\cup\breve \Gamma_K$ for $t\in\Ti$,
where $M'\le K\le M$. Therefore,
\begin{equation}\label{boundary1}\partial{V}(t)\setminus
O_z=\breve\Gamma_{1}\cup\dots\cup\breve\Gamma_K\cup
A_{M+1}(t)\cup\dots\cup A_N(t),\quad\ t\in\Ti.\end{equation}

Let $t_1,t_2\in\Ti$ and $t_1<t_2$. The next geometrical objects
plays an~important role in the estimates below: for
$t\in(t_1,t_2)$ we define the level set $S_k(t,t_1,t_2)\subset
\{x\in\D:\Phi_k(x)=-t\}$ separating cycles $\cup_{j=M+1}^N
A_j(t_1)$ from $\cup_{j=M+1}^NA_j(t_2)$ as follows. Namely, take
arbitrary $t',t''\in\Ti$ such that $t_1<t'<t''<t_2$. From
Properties~(ii),(iv) we have the uniform convergence
$\Phi_k|_{A_j(t_1)}\rightrightarrows -t_1$,
$\Phi_k|_{A_j(t_2)}\rightrightarrows -t_2$ as $k\to\infty$ for
every~$j=M+1,\dots,N$. Thus there exists
$k_\circ=k_\circ(t_1,t_2,t',t'')\in\N$ such that for all $k\ge
k_\circ$
\begin{equation}\label{boundary0}\Phi_k|_{A_j(t_1)}>-t',\quad
\Phi_k|_{A_j(t_2)}< -t''\quad\forall j=M+1,\dots,N.\end{equation}
In particular,
\begin{equation}\label{boundary2}
\begin{array}{lcr}
\forall t\in[t',t'']\ \forall k\ge k_\circ\quad
\Phi_k|_{A_j(t_1)}> -t,\quad\Phi_k|_{A_j(t_2)}< -t,\quad
\\
\\
\forall
j=M+1,\dots,N.
\end{array}
\end{equation}

For $k\ge k_\circ$, $j=M+1,\dots,N$, and $t\in[t',t'']$ denote by
$W^j_k(t_1,t_2;t)$ the connected component of the open set $\{x\in
V(t_2)\setminus \overline V(t_1):\Phi_k(x)>-t\}$ such that
$\partial W^j_k(t_1,t_2;t)\supset A_j(t_1)$ (see Fig.1) and put
$$W_k(t_1,t_2;t)=\bigcup\limits_{j=M+1}^N W^j_k(t_1,t_2;t),\qquad
S_k(t_1,t_2;t)=(\partial W_k(t_1,t_2;t))\cap V(t_2)
\setminus\overline V(t_1).$$ Clearly, $\Phi_k\equiv -t$ on
$S_k(t_1,t_2;t)$. By construction (see Fig.1),
\begin{equation}\label{boundary3}\partial W_k(t_1,t_2;t)=S_k(t_1,t_2;t)\cup A_{M+1}(t_1)\cup\dots\cup
A_N(t_1).\end{equation} (Note that $W_k(t_1,t_2;t))$ and
$S_k(t_1,t_2;t)$ are well defined for all $t\in[t',t'']$ and $k\ge
k_\circ=k_\circ(t_1,t_2,t',t'')$.)

Since by (E--NS) each $\Phi_k$ belongs to   $W^{2,2}_{\loc}(\D)$,
by the Morse-Sard theorem for Sobolev functions (see
Theorem~\ref{kmpTh1.1}) we have that for almost all $t\in[t',t'']$
the level set $S_k(t_1,t_2;t)$ consists of finitely many
$C^1$-cycles and $\Phi_k$ is differentiable (in classical sense)
at every point~$x\in S_k(t_1,t_2;t)$ with $\nabla\Phi_k(x)\ne0$.
The values $t\in[t',t'']$ having the above property will be called
$k$-{\it regular}.

Recall, that for a set $A\subset P_+$ we denote by $\widetilde A$
the set in $\R^3$ obtained by rotation of $A$ around $O_z$-axis.
By construction, for every regular value~$t\in[t',t'']$ the set
$\widetilde S_{k}(t',t'';t)$ is a finite union of smooth surfaces
(tori), and
\begin{equation}\label{lac-2-ax}
\int_{\widetilde S_k(t_1,t_2;t)}\nabla\Phi_k\cdot{\bf
n}\,dS=-\int_{\widetilde S_k(t_1,t_2;t)}|\nabla\Phi_k|\,dS<0,
\end{equation}
where $\n$ is the unit outward normal vector to
$\partial\widetilde W_k(t_1,t_2;t)$.

For $h>0$ denote
$\Gamma_h=\{x\in\Omega:\dist(x,\Gamma_1\cup\dots\cup\Gamma_K)=h)\}$,
$\Omega_h=\{x\in\Omega:\dist(x,\Gamma_1\cup\dots\cup\Gamma_K)<h)\}$.
Since the distance function $\dist(x,\partial\Omega)$ is
$C^1$--regular and the norm of its gradient is equal to one in the
neighborhood of $\partial\Omega$, there is a constant $\delta_0>0$
such that for every $h\le\delta_0$ the set $\Gamma_h$ is a union
of $K$ \ $C^1$-smooth  surfaces homeomorphic to balls or torus,
and
\begin{equation}\label{lac0.1-ax}
\Ha^2(\Gamma_h)\le c_0\quad\forall h\in(0,\delta_0],
\end{equation}
where the constant $c_0=3\Ha^2(\Gamma_1\cup\dots\cup\Gamma_K)$ is
independent of~$h$.

By direct calculations, (\ref{2.1'}) implies
\begin{equation}\label{grthpax}\nabla\Phi=\ve\times\ov\quad\mbox{ in }\Omega,
\end{equation}
where $\ov=\curl\ve$, i.e.,
$$\ov=(\omega_r,\omega_\theta,\omega_z)=\bigl(-\frac{\partial v_\theta}{\partial z},\ \frac{\partial v_r}{\partial z}-\frac{\partial v_z}{\partial r},\ \frac{v_\theta}r+\frac{\partial v_\theta}{\partial r}\bigr).$$
Set $\omega(x)=|\ov(x)|$. Since $\Phi\neq \const$ on $V(t)$,
\eqref{grthpax} implies $\int_{{\widetilde V}(t)}\omega^2\,dx>0$
for every $t\in\Ti$. Hence, from the weak convergence
$\ov_k\rightharpoonup\ov$ in $L^2(\Omega)$ (recall, that $\ov_k=\curl \ue_k$,
$\omega_k(x)=|\ov_k(x)|$\,) it follows

\begin{lem}
\label{lkr8-ax} {\sl For any ${t}\in\Ti$ there exist constants
$\e_{t}>0$, $\delta_{t}\in(0,\delta_0)$ and $k_{t}\in\N$ such that
for all $k\ge k_t$
$$ A_j({t})\Subset \frac12B_k,\quad\overline\Omega_{\delta_{t}}\cap
A_j({t})=\emptyset,\quad j=M+1,\dots,N, $$
$$\Gamma_j\Subset \frac12B_k,\quad j=1,\dots,N,$$
 and
\begin{equation}
\label{omeg} \int\limits_{{\widetilde V}(t)\cap
\frac12B_k\setminus
\Omega_{\delta_t}}\omega_k^2\,dx>\varepsilon_t.
\end{equation}
Here $B_k=\{x\in\R^3:|x|<R_k\}$ are the balls where the
solutions $\ue_k\in W^{1,2}(\Omega\cap B_k)$ are defined.}
\end{lem}

Now we are ready to prove the key estimate.

\begin{lem}
\label{ax-lkr11}{\sl For  any $t_1,t_2,t',t''\in\Ti$ with
$t_1<t'<t''<t_2$ there exists $k_*=k_*(t_1,t_2,t',t'')$ such that
for every $k\ge k_*$ and for almost all $t\in[t_1,t_2]$
the inequality
\begin{equation}\label{mec}
\int\limits_{\widetilde S_k(t_1,t_2;t)}|\nabla\Phi_k|\,dS<\F t,
\end{equation}
holds with the constant $\F$ independent of  $t,t_1,t_2,t',t''$
and $k$. }
\end{lem}

\pr
Fix $t_1,t_2,t',t''\in\Ti$ with $t_1<t'<t''<t_2$. Below we always
assume that $k\ge k_\circ(t_1,t_2,t',t'')$ (see
(\ref{boundary0})--(\ref{boundary2})\,), in particular, the set
$S_k(t_1,t_2;t)$ is well defined for all $t\in[t',t'']$. We
assume also that $R_k>2$ and $k\ge k_{t_1}$ (see Lemma~\ref{lkr8-ax}).

The main idea of the proof of~(\ref{mec}) is quite simple: we will
integrate the equation
\begin{equation}\label{cle_lap**}\Delta\Phi_k=\omega_k^2+\frac1{\nu_k}\div(\Phi_k\ue_k)-
\frac{1}{\nu_k}\fe_k\cdot\ue_k\end{equation} over the suitable
domain $\Omega_k(t)$ with $\partial\Omega_k(t)\supset \widetilde
S_k(t_1,t_2;t)$ such that the corresponding boundary integrals
\begin{equation}\label{r**1}
\biggr|\int_{\bigl(\partial\Omega_k(t)\bigr)\setminus \widetilde
S_k(t_1,t_2;t)}\nabla\Phi_k\cdot{\bf n}\,dS\biggr|
\end{equation}
\begin{equation}\label{r**2}
\frac1{\nu_k}\biggr|\int_{\bigl(\partial\Omega_k(t)\bigr)\setminus
\widetilde S_k(t_1,t_2;t)}\Phi_k\ue_k\cdot{\bf n}\,dS\biggr|
\end{equation}
are negligible. We split the construction of the domain
$\Omega_k(t)$ into several steps.

\medskip

{\sc Step 1.} We claim that for any $\e>0$ the estimate
\begin{equation}\label{ax-key1}
\frac1{{R}^2}\biggl|\int\limits_{S_{R}}\Phi_k\,dS\biggr|<\e
\end{equation}
holds for sufficiently large~$R$ uniformly with respect to~$k$. Indeed, the weak convergence
$\Phi_k\rightharpoonup \Phi\mbox{ \ in \
}W^{1,3/2}_\loc(\overline\Omega)$ implies
\begin{equation}\label{nax-key1}
\frac1{R_*^2}\biggl|\int\limits_{S_{R_*}}\Phi_k\,dS\biggr|\to
\frac1{R_*^2}\biggl|\int\limits_{S_{R_*}}\Phi\,dS\biggr|
\end{equation}
for any fixed $R_*>R_0$. Take $R_\e>R_0$ sufficiently large such that
\begin{equation}\label{nth1}
\frac1{R_\e^2}\biggl|\int\limits_{S_{R_\e}}\Phi\,dS\biggr|<\frac\e2
\end{equation}
(this inequality holds for sufficiently large~$R_\e$ because of the inclusion $\Phi\in
L^3(\Omega)\cap D^{1,3/2}(\Omega)$, see the proof of~(\ref{1pub})\,). Take arbitrary
$R\in[R_\e,R_k]$. Then $R\le 2^l R_\e$ for some
$l\in\N$. We have
\begin{equation}\label{nth2}
\begin{array}{lcr}\displaystyle \biggl|\frac1{R^2}\int\limits_{S_{R}}\Phi_k\,dS-\frac1{R_\e^2}\int\limits_{S_{R_\e}}\Phi_k\,dS
\biggr|=
\biggl|\int\limits_{|x|\in[R_\e,R]}\frac1{|x|^3}x\cdot\nabla\Phi_k\,dx\biggr|\\
\displaystyle \le\biggl|\sum\limits_{m=0}^{l-1}
\int\limits_{|x|\in[2^mR_\e,2^{m+1}R_\e]}\frac1{|x|^3}x\cdot\nabla\Phi_k\,dx\biggr|
\\ \displaystyle \le \sum\limits_{m=0}^{l-1}
\biggl(\int\limits_{|x|\in[2^mR_\e,2^{m+1}R_\e]}|\nabla\Phi_k|^{\frac32}\,dx\biggr)^{\frac23}
\biggl(\int\limits_{|x|\in[2^mR_\e,2^{m+1}R_\e]}\frac1{|x|^6}\,dx\biggr)^{\frac13}
\\ \le\displaystyle C \sum\limits_{m=0}^{l-1}\sqrt{2^mR_\e}\frac1{2^mR_\e}\le
C' \frac1{\sqrt{R_\e}},
\end{array}
\end{equation}
where the constants $C$, $C'$  do not depend on~$l$ and $k$ (here we have used the
estimate~(\ref{epgr34})\,). Consequently, if we take
a~sufficiently large~$R_\e$, then
$$\biggl|\frac1{R^2}\int\limits_{S_{R}}\Phi_k\,dS-\frac1{R_\e^2}\int\limits_{S_{R_\e}}\Phi_k\,dS
\biggr|<\frac\e2$$ for all $k\in\N$ and $R>R_\e$. Now the
required estimate~(\ref{ax-key1}) follows from the last inequality
and formulas~(\ref{nth1}), (\ref{nax-key1}) (with $R_*=R_\e$).

\medskip
{\sc Step 2.}  By  direct calculations, (\ref{NSk}) implies
$$\nabla\Phi_k=-\nu_k\curl\,\ov_k+\ov_k\times\ue_k+\fe_k=-\nu_k\curl\,\ov_k+\ov_k\times\ue_k+\frac{\nu_k^2}{\nu^2}\,\curl\,{\bf b}.$$
By the Stokes theorem, for any $C^1$-smooth closed surface
$S\subset\Omega$ and ${\bf g}\in W^{2,2}(\Omega)$ we have
$$\int_S\curl{\bf g}\cdot{\bf n}\,dS=0.$$  So, in particular,
\begin{equation}\label{uesk0}\int_S\nabla\Phi_k\cdot{\bf
n}\,dS=\int_S(\ov_k\times\ue_k)\cdot{\bf n}\,dS.
\end{equation}
Since by construction for every $ x\in
S_{R_k}=\{y\in\R^3:|y|=R_k\}$ there holds the equality
\begin{equation}\label{uesk}
\ue_k(x)\equiv \frac{\nu_k}{\nu}\h(x)\equiv
-\nu_k\frac{\sum_{i=1}^N F_i}{4\pi\nu}\,\frac{x}{|x|^3},
\end{equation}
we see that
\begin{equation}\label{lac0.5}
\int_{S_{R_k}}\nabla\Phi_k\cdot{\bf n}\,dS=0.
\end{equation}
Indeed, by (\ref{uesk}) the vector $\ue_k(x)$ is parallel to $x$ for
every $x\in S_{R_k}$,
consequently,
$\ov_k(x)\times\ue_k(x)$ is perpendicular to~$x$ for $x\in S_{R_k}$, and by virtue of~(\ref{uesk0}) we obtain~(\ref{lac0.5}).

Furthermore, using~(\ref{uesk}) we get
\begin{equation}\label{r**7}
\frac1{\nu_k}\biggr|\int_{S_{R_k}}\Phi_k\ue_k\cdot{\bf
n}\,dS\biggr|=
\frac{C}{R_k^2}\biggr|\int_{S_{R_k}}\Phi_k\,dS\biggr|.
\end{equation}
 Thus applying (\ref{ax-key1}) for
sufficiently large~$k$ we have
\begin{equation}\label{remop2'}
\frac1{\nu_k}\biggr|\int_{S_{R_k}}\Phi_k\ue_k\cdot{\bf
n}\,dS\biggr|<\e.
\end{equation}

\medskip

{\sc Step 3.} Denote~$\Gamma_0=\Gamma_1\cup\dots\cup\Gamma_K$.
Recall, that by the pressure normalization condition,
\begin{equation}\label{remo1}
\Phi|_{\Gamma_0}=0.\end{equation} From uniform boundedness
$\|\Phi_k\|_{L^3(\Omega_{\delta_0})}+\|\nabla\Phi_k\|_{L^{3/2}(\Omega_{\delta_0})}\le
C$ we easily have
\begin{equation}\label{remob0.4}
\int\limits_{\Gamma_{h}}\Phi_k^2\,dS<\sigma_0 \quad\forall
h\in(0,\delta_0]\ \ \forall k\in\N,
\end{equation}
where $\sigma_0>0$ is independent of~$k,h$ (recall, that
$\Omega_h=\{x\in\Omega:\dist(x,\Gamma_0)\le h\}$,\ \ $\Gamma_h=\{x\in\Omega:\dist(x,\Gamma_0)=h\}$\,). For $r_0>0$
denote $$T_{r_0}=\{(x_1,x_2,x_3)\in\R^3:x_1^2+x_2^2<r_0^2\},$$
$$C_{r_0}=\partial T_{r_0}=\{(x_1,x_2,x_3)\in\R^3:x_1^2+x_2^2=r_0^2\}.$$
We claim, that for every~$r_0>0$ and for any~$\sigma>0$ there
exist constants $\delta_{r_0}(\sigma)\in(0,\delta_0)$,
$k_{r_0}(\sigma)\in\N$ such that
\begin{equation}\label{belac0.4}
\int\limits_{\Gamma_{h}\setminus T_{r_0}}\Phi_k^2\,dS<\sigma^2
\quad\forall h\in(0,\delta_{r_0}(\sigma)]\ \ \forall k\ge
k_{r_0}(\sigma).
\end{equation}
Indeed,  by the classical formula of changing variables in integral,
there exists a constant~$\delta_1\in(0,\delta_0)$ such that for
$h'<h''<\delta_1$ the following formula
\begin{equation}\label{us1}\begin{array}{lcr}
\int\limits_{\Gamma_{h''}}g^2\,dS\leq C_1(h',h'')
\int\limits_{\Gamma_{h'}}g^2\,dS+
\\
\\
\qquad\qquad+C_2(h',h'')\int\limits_{\Omega_{h''}\setminus\Omega_{h'}}|g\nabla
g|\,dx, \qquad \forall g\in W^{1,3/2}(\Omega_h),
\end{array}
\end{equation}
holds. Here
$C_1(h',h'')\to1$ and $C_2(h',h'')\to2$ as $|h'-h''|\to0$. Then
the boundary conditions~(\ref{remo1}) and
$\|\Phi\|_{L^{3}}+\|\nabla\Phi\|_{L^{3/2}}<\infty$ imply that for
every~$\sigma>0$ there exists $h_\sigma\in(0,\delta_1)$ such that
\begin{equation}\label{belac0.4'}
\int\limits_{\Gamma_{h}}\Phi^2\,dS<\frac16\sigma^2 \quad\forall
h\in(0,h_\sigma],
\end{equation}
\begin{equation}\label{belac0.5}
\int\limits_{\Omega_{h_\sigma}}|\Phi\nabla\Phi|\,dx<\frac16\sigma^2.
\end{equation}
From the weak convergence $\nabla\Phi_k\rightharpoonup \nabla\Phi$
in $L^{3/2}(\Omega_{h_\sigma})$, the axial symmetry and from the
Sobolev embedding theorems for plane domains it follows that $\Phi_k\to\Phi$ strongly in
$L^{3}(\Omega_{h_\sigma}\setminus T_{r_0})$ for any fixed $r_0>0$.
Thus  $\Phi_k\nabla\Phi_k\to \Phi\nabla\Phi$ strongly in
$L^1(\Omega_{h_\sigma}\setminus T_{r_0})$ and, consequently,
\begin{equation}\label{belac0.6}
\int\limits_{\Omega_{h_\sigma}\setminus
T_{r_0}}|\Phi_k\nabla\Phi_k|\,dx<\frac16\sigma^2\quad\mbox{ for
sufficiently large }k.
\end{equation}
Moreover, the uniform convergence $\Phi_{k}|_{\Gamma_h\setminus
T_{r_0}}\rightrightarrows\Phi|_{\Gamma_h\setminus T_{r_0}}$ holds\footnote{The convergence holds possibly only over a subsequence, which we denote again by $\Phi_k$.}
as $k\to\infty$ for almost all $h\in (0,h_\sigma)$ (see
\cite{Amick}, \cite{kpr}\footnote{In \cite{Amick} Amick proved the
uniform convergence $\Phi_k\rightrightarrows\Phi$ on almost all
circles. However, his method can be easily modified to    prove
the uniform convergence on almost all level lines of every
$C^1$-smooth function with nonzero gradient. Such modification was
done in the proof of Lemma~3.3 of \cite{kpr}.}; cf. also with Lemma~\ref{regc-ax} above). The last
assertion, together with~(\ref{belac0.6}), (\ref{belac0.4'}), and
(\ref{us1}) implies the required claim~(\ref{belac0.4}).

Our purpose on the next several steps is as follows: for arbitrary
$\e>0$ and for sufficiently large~$k$ to construct a~surface
$S_k(\e)\subset\Omega_{h}$, $h\ll1$, separating $\Gamma_0$ from other boundary components and from
infinity such that
\begin{equation}\label{remop1}
\biggr|\int_{S_k(\e)}\nabla\Phi_k\cdot{\bf
n}\,dS\biggr|=\biggl|\int_{S_k(\e)}(\ov_k\times\ue_k)\cdot{\bf
n}\,dS\biggr| <3\e,
\end{equation}
\begin{equation}\label{remop2}
\frac1{\nu_k}\biggr|\int_{S_k(\e)}\Phi_k\ue_k\,dS\biggr|<3\e.
\end{equation}
This construction is more complicated because of the singularity
axis~$O_z$: we do not have the uniform smallness
of~$\int\limits_{\Gamma_{h}}\Phi^2\,dS$ near the singularity line;
these integrals are only uniformly bounded (see~(\ref{remob0.4}), (\ref{belac0.4})\,).

Recall, that in our notation $\ue_k=\frac{\nu_k}{\nu}\h+\w_k$,
where $\w_k\in H(\Omega)$, $\|\w_k\|_{H(\Omega)}=1$, and $\h$ is a
solenoidal extension of~$\a$ from Lemma~\ref{kmpLem14.1}.

\medskip

{\sc Step 4.} We claim that for arbitrary $\e>0$ there exists
a~constant $r_0=r_0(\e)>0$ such that
\begin{equation}\label{remo5}
\frac{\nu_k}{\nu}\!\!\!\int_{\Gamma_h\cap
T_{2r_0}}\!\!\!\bigl(|\h|\cdot|\nabla\h|+|\h|\cdot|\nabla\w_k|+|\w_k|\cdot|\nabla\h|\bigr)\,dS<\e\quad\forall
h\in(0,\delta_0],
\end{equation}
\begin{equation}\label{remp6}
\frac{\nu^2_k}{\nu^2}\int_{\Gamma_h\cap
T_{2r_0}}|\h|^2\,dS<\frac13\e\sigma^2_0\qquad\forall
h\in(0,\delta_0],
\end{equation}
where~$\sigma_0$ is a constant from~(\ref{remob0.4}). This claim
is easily  deduced from the uniform estimates
\begin{equation}\label{est_i1}
\|\h\|_{L^\infty(\Omega)}\le C \|\h\|_{W^{2,2}(\Omega)}<\infty,
\end{equation}
\begin{equation}\label{est_i2}
\|\nu_k\ue_k\|_{W^{2,3/2}(\Omega_{\delta_0})}\le C,
\end{equation}
where $C$ is
independent of~$k$ (the last inequality follows from the
construction (see~(\ref{NSk})\,)  by well-known
estimates~\cite{Agmon}, \cite{Solonnikov} for the solutions to the
Stokes system, see also formula~(\ref{p15}) of the present paper). For example, prove that
for sufficiently small~$r_0$ the estimate
\begin{equation}\label{est_i3}
\frac{\nu_k}{\nu}\!\!\!\int_{\Gamma_h\cap
T_{2r_0}}\!\!\!|\h|\cdot|\nabla\w_k|\,dS<\e\quad\forall
h\in(0,\delta_0]
\end{equation}
holds for all~$k$. Indeed, (\ref{est_i2}) implies
$\|\nu_k\nabla\w_k\|_{W^{1,3/2}(\Omega_{\delta_0})}\le C_1$, consequently, by Sobolev Imbedding Theorem,
\begin{equation}\label{est_i4}
\int_{\Gamma_h}\nu_k^2|\nabla\w_k|^2\,dS\leq C_2\qquad\forall h\in(0,\delta_0],
\end{equation}
where $C_2$ does not depend on~$h,k$. Thus
\begin{equation}\label{est_i5}
\begin{array}{lcr}
\dfrac{\nu_k}{\nu}\int\limits_{\Gamma_h\cap
T_{2r_0}}|\h|\cdot|\nabla\w_k|\,dS\\
\\
\leq\dfrac1\nu\|\h\|_{L^\infty(\Omega)}
\biggl(\int\limits_{\Gamma_h}\nu_k^2|\nabla\w_k|^2\,dS\biggr)^{\frac12}\Ha^2(T_{2r_0}\cap\Gamma_h)^{\frac12}\le
C_3\sqrt{r_0},
\end{array}
\end{equation}
that gives the required estimate~(\ref{est_i3}) (here $\Ha^2$ means the two-dimensional
Hausdorff measure, i.e., area of the set). Other estimates in (\ref{remo5})--(\ref{remp6})
are proved analogously.

\medskip

{\sc Step 5.} Take arbitrary $\e>0$ and $\lambda>0$. We claim that
for sufficiently large $k$ there exists a distance $\bar h_k=\bar
h_k(\e,\lambda)$ such that
\begin{equation}\label{remo6}
\bar h_k\in (\lambda \nu_k^2,\lambda \nu_k^2j_\e],
\end{equation}
\begin{equation}\label{remo8}
\int_{\Gamma_{\bar h_k}}|\w_k|^2<c_1\lambda\nu_k^2 j_\e,
\end{equation}
\begin{equation}\label{remo7}
\int_{\Gamma_{\bar h_k}}|\w_k|\cdot|\nabla\w_k|\,dS<\e,
\end{equation}
\begin{equation}\label{remo9}
\frac{\nu_k}{\nu}\int_{\Gamma_{\bar
h_k}}\bigl(|\h|\cdot|\nabla\h|+|\h|\cdot|\nabla\w_k|+|\w_k|\cdot|\nabla\h|\bigr)\,dS<\frac{c_2}{\sqrt{\lambda}},
\end{equation}
where $j_\e\in\N$ depends on $\e$ only, and $c_1,c_2$ are independent
of~$k,\e,\lambda$.

To prove the above claim, fix $\lambda>0$ and $\e>0$ and put
$h_j=\lambda \nu_k^2j$,
$U_j=\{x\in\Omega:\dist(x,\Gamma_0)\in(h_{j-1},h_j)\}$, \ $\xi_j=\sup\limits_{h\le h_j}\int\limits_{\Gamma_{ h}}\w_k^2\,dS$.
Then
$$\xi_j\le c\int\limits_{U_1\cup\dots\cup
U_j}|\w_k|\cdot|\nabla\w_k|\,dx\le c\biggl(\int\limits_{U_1\cup\dots\cup
U_j}|\w_k|^2\,dx\biggr)^{\frac12}\biggl(\int\limits_{U_1\cup\dots\cup
U_j}|\nabla\w_k|^2\,dx\biggr)^{\frac12}
$$
$$
\le c\bigl(\xi_jh_j\bigr)^{\frac12}$$
(here we use the identities $\|\nabla\w_k\|_{L^2(\Omega)}=1$ and
$\int\limits_{U_1\cup\dots\cup
U_j}|\w_k|^2\,dx=\int_0^{h_j}dh\int_{\Gamma_h}\w_k^2\,dS$\,).
Hence $\xi_j\le c_1 h_j$ with $c_1=c^2$, i.e.,
\begin{equation}\label{remo10}
\int_{\Gamma_{ h}}\w_k^2\,dS\le
c_1\lambda\nu_k^2j\qquad\forall h\le h_j.
\end{equation}
Consequently,
\begin{equation}\label{remo11}\begin{array}{rcl}
\int\limits_{U_j}|\w_k|\cdot|\nabla\w_k|\,dx\le
\biggl(\int\limits\limits_{U_j}|\w_k|^2\,dx\biggr)^{\frac12}\cdot
\biggl(\int\limits\limits_{U_j}|\nabla\w_k|^2\,dx\biggr)^{\frac12}\\
\\
\le c_3\biggl(\lambda\nu_k^2j[h_j-h_{j-1}]\biggr)^{\frac12}\cdot
\biggl(\int\limits_{U_j}|\nabla\w_k|^2\,dx\biggr)^{\frac12}\\
\\
 = c_3
[h_j-h_{j-1}]
\biggl(j\int\limits\limits_{U_j}|\nabla\w_k|^2\,dx\biggr)^{\frac12}.\end{array}
\end{equation}
We need to find an index $j\in\N$ such that the right-hand side of the last inequality is sufficiently small. To this end, denote $j_{\max}=\min\{j\in
\N:c_3^2j\int\limits_{U_j}|\nabla\w_k|^2\,dx<\frac14\e^2\}$, i.e.,
$$\int\limits_{U_j}|\nabla\w_k|^2\,dx\ge\frac{\e^2}{4jc_3^2}\qquad\ \forall j<j_{\max}.$$
Therefore,
by virtue of $\int\limits_{\Omega}|\nabla\w_k|^2\,dx=1$ and
$1+\frac12+\dots+\frac{1}{j_{\max}}\sim \ln j_{\max}$ we have $j_{\max}<
c_4\exp\bigl(\frac{4c^2_3}{\e^2}\bigr)+1$.  Denote the right-hand
side of the last estimate by $j_\e$, i.e., $j_\e=c_4\exp\bigl(\frac{4c^2_3}{\e^2}\bigr)+1$. By construction, $j_\e$ is
independent of~$k,\lambda$. We assume here that $k$ is
sufficiently large so that $h_{j_{\max}}\le
h_{j_{\e}}<\delta_0$. By choice of $j_{\max}$ we get
\begin{equation}\label{est_p2}
\frac1{h_{j_{\max}}-h_{j_{\max}-1}}
\int\limits_{U_j}|\w_k|\cdot|\nabla\w_k|\,dx<\frac\e2.
\end{equation}
We have also for arbitrary nonnegative
function $f\in L^2(\Omega)$
\begin{equation}\label{remo12}
\begin{array}{lcr}
\nu_k\int\limits_{U_j}f\,dx\le c_5\nu_k
\biggl(\int\limits_{U_j}f^2\,dx\biggr)^{\frac12}\cdot
[h_j-h_{j-1}]^{\frac12}
\\
\\
\qquad\qquad =c_5 \dfrac1{\sqrt{\lambda}}[h_j-h_{j-1}]
\biggl(\int\limits_{U_j}f^2\,dx\biggr)^{\frac12}.
\end{array}
\end{equation}
Applying this fact to the function in the left-hand side of~(\ref{remo9}), we obtain
\begin{equation}\label{remo_pe}
\begin{array}{lcr}
\dfrac{\nu_k}{\nu[h_{j_{\max}}-h_{j_{\max}-1}]}\int\limits_{U_{j_{\max}}}
\bigl(|\h|\cdot|\nabla\h|+|\h|\cdot|\nabla\w_k|+|\w_k|\cdot|\nabla\h|\bigr)\,dS
\\
\\
<\dfrac{c_2}{2\sqrt{\lambda}}.
\end{array}
\end{equation}
Estimates (\ref{est_p2}) and (\ref{remo_pe}) together with identity
$$\int\limits_{U_{j_{\max}}}f\,dx=\int\limits_{h_{j_{\max}-1}}^{h_{j_{\max}}}\,dh\int_{\Gamma_h}f\,dS$$
imply the existence of $\bar
h_k\in(h_{j_{\max}-1},h_{j_{\max}})$ with the required properties
(\ref{remo7})--(\ref{remo9}). By construction, the property~(\ref{remo6}) is fulfilled as well, and~(\ref{remo8}) follows immediately from~(\ref{remo10}). This finishes the proof of the claim of Step~5.

\medskip

{\sc Step 6.} Now let us define the target surface $S_k(\e)$ with
properties (\ref{remop1})--(\ref{remop2}). Take arbitrary $\e>0$
and fix it. We apply the last Step~5 two times. First, take
$\lambda'>0$ sufficiently large in order to satisfy the condition
\begin{equation}\label{remo13}
\frac{c_2}{\sqrt{\lambda'}}<\e,
\end{equation}
where $c_2$ is a constant from the estimate~(\ref{remo9}). Let $\bar
h'_k$ be the corresponding distance from Step~5. Take the
parameter~$r_0>0$ from Step~4. Since $\bar h'_k\to0$ as
$k\to\infty$, from~(\ref{belac0.4}) and (\ref{remo8}) we  have for sufficiently large~$k$ the inequality
\begin{equation}\label{remo14}
\frac1{\nu_k}\biggr|\int_{\Gamma_{\bar h'_k}\setminus
T_{r_0}}\Phi_k\ue_k\,dS\biggr|<\e.
\end{equation}
Furthermore,
estimates~(\ref{remo7})--(\ref{remo9}) yield
\begin{equation}\label{remop15}
\int_{\Gamma_{\bar h'_k}}|\ue_k|\cdot|\nabla\ue_k|\,dS<2\e.
\end{equation}
Second, take $\lambda''>0$ sufficiently small such that
\begin{equation}\label{remo16}
\lambda'' j_\e\sigma_0<\frac13\e^2,
\end{equation}
where $j_\e,\sigma_0$ are parameters from formulas~(\ref{remo8})
and (\ref{remob0.4}). Let ${\bar h''_k}$ be the corresponding
distance from Step~4. Then, by formulas (\ref{remo8}),
(\ref{remp6}) and (\ref{remob0.4}), we have
\begin{equation}\label{remo17}
\frac1{\nu_k}\biggr|\int_{\Gamma_{\bar h''_k}\cap
T_{2r_0}}\Phi_k\ue_k\,dS\biggr|<\e
\end{equation}
for sufficiently large~$k$. Moreover, from (\ref{remo5}) and
(\ref{remo7}) it follows that
\begin{equation}\label{remo18}
\int_{\Gamma_{\bar h''_k}\cap
T_{2r_0}}|\ue_k|\cdot|\nabla\ue_k|\,dS<2\e.
\end{equation}
The required surface $S_k(\e)$ will be the union of three parts (see Fig.2):
$\Gamma_{\bar h'_k}\setminus T_{r_1}$ with $\Gamma_{\bar
h''_k}\cap T_{r_1}$ and $C_{r_1}\cap\Omega_{\bar
h'_k}\setminus\Omega_{\bar h''_k}$, where~$r_1\in[r_0,2r_0]$ is such that the following conditions \begin{equation}\label{remo19}
\int_{C_{r_1}\cap\Omega_{\bar
h'_k}}|\ue_k|\cdot|\nabla\ue_k|\,dS<\e,
\end{equation}
\begin{equation}\label{remo20}
\frac1{\nu_k}\biggr|\int_{C_{r_1}\cap\Omega_{\bar
h'_k}}\Phi_k\ue_k\,dS\biggr|<\e
\end{equation}
hold (we want to achieve the required
properties~(\ref{remop1})--(\ref{remop2})).

\begin{center}
\includegraphics[scale=0.4]{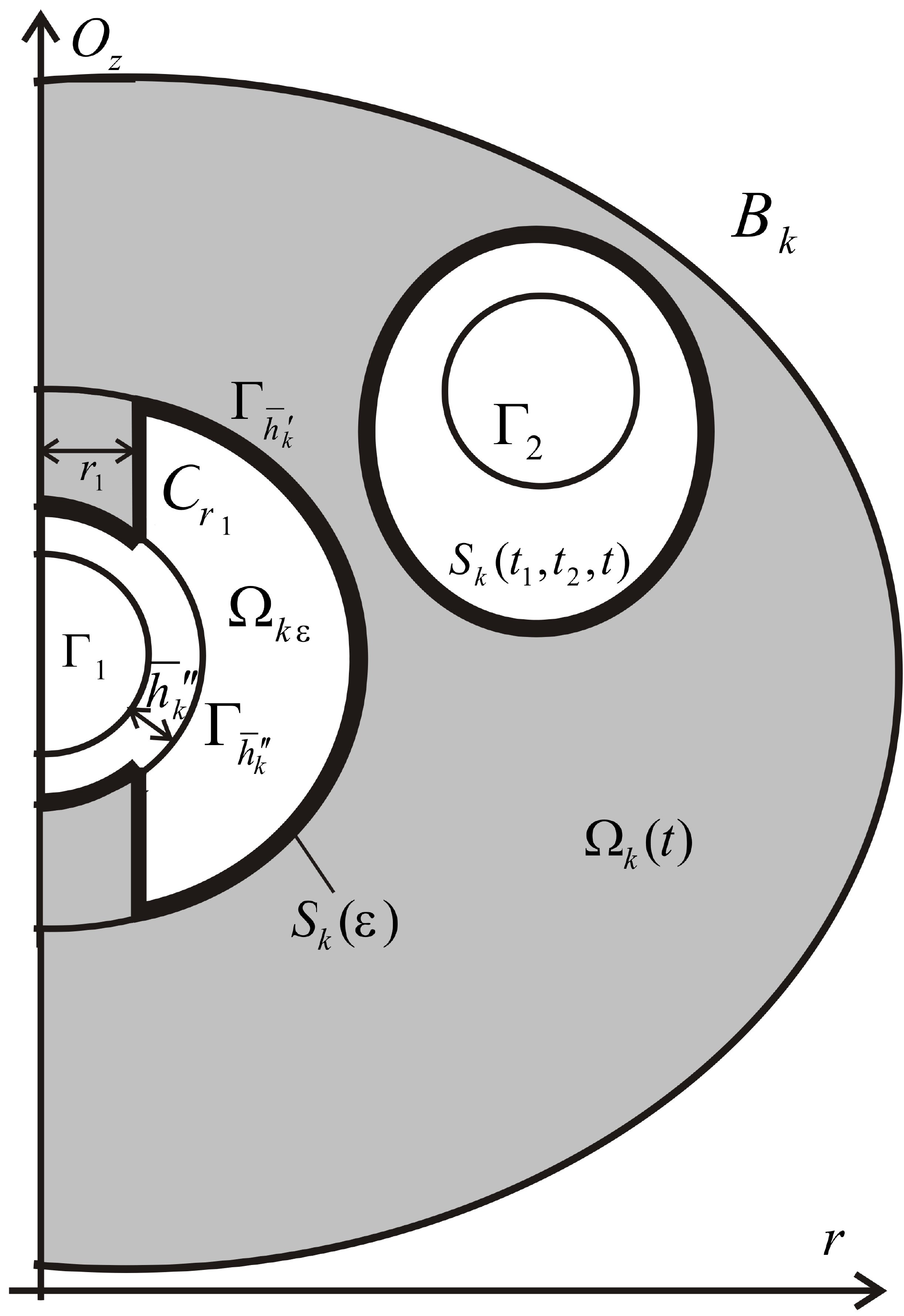}
\end{center}
\begin{center}
Fig. 2. {\sl The domain $\Omega_k(t)$ for the case $M=K=1$, $N=2$, $\Gamma_0=\Gamma_1$. }
\end{center}

So, to finish the construction $S_k(\e)$, we need to prove the
existence of~$r_1\in[r_0,2r_0]$
with properties~(\ref{remo19})--(\ref{remo20}). By one dimensional Poincar\`{e} inequality (applied on segments parallel to $O_z$-axis during the integration by Fubini theorem) we have
\begin{equation}\label{remo21}
\int_{T_{2r_0}\cap\Omega_{\bar
h'_k}}\w_k^2\,dx<C(r_0)(h'_k)^2\int_{T_{2r_0}\cap\Omega_{\bar
h'_k}}|\nabla \w_k|^2\,dx\le C(\e,r_0)\nu_k^4
\end{equation}
(second inequality follows from the estimate~$h'_k\le c(\e)\nu_k^2$,
see~(\ref{remo6}) and the beginning of Step~6\,).
Therefore,
\begin{equation}\label{remo22}
\int_{T_{2r_0}\cap\Omega_{\bar h'_k}}\ue_k^2\,dx\le
C'(\e,r_0)\nu_k^4,
\end{equation}
\begin{equation}\label{remo23}
\int_{T_{2r_0}\cap\Omega_{\bar
h'_k}}|\ue_k|\cdot|\nabla\ue_k|\,dx\le C''(\e,r_0)\nu^2_k.
\end{equation}
Here we have used the identity
$\ue_k=\w_k+\frac{\nu_k}{\nu}\h$ and the estimate~$\|\h\|_{L^\infty}<\infty$\,.
By Fubini Theorem, for an~integrable function~$f$ the identity
$$\int_{T_{2r_0}\cap\Omega_{\bar
h'_k}}f=\int\limits_0^{2r_0}dr\int\limits_{C_r\cap\Omega_{\bar
h'_k}}f\,dS$$ holds. Therefore we can choose~$r_1\in[r_0,2r_0]$
satisfying
\begin{equation}\label{remo22ll}
\int_{C_{r_1}\cap\Omega_{\bar h'_k}}\ue_k^2\,dS\le
\frac1{r_0}C'(\e,r_0)\nu_k^4,
\end{equation}
\begin{equation}\label{remo23ll}
\int_{C_{r_1}\cap\Omega_{\bar
h'_k}}|\ue_k|\cdot|\nabla\ue_k|\,dS\le
\frac1{r_0}C''(\e,r_0)\nu^2_k.
\end{equation}
Since $\nu_k\to0$ as $k\to\infty$, but $r_0$ and $\e$ are fixed
(they do not depend on~$k$), the
estimates~(\ref{remo22})--(\ref{remo23})
imply~(\ref{remo22ll})--(\ref{remo23ll}) for sufficiently
large~$k$.

So, we have constructed the required Lipschitz
surface~$S_k(\e)\subset\Omega_{2\bar h'_k}$ with properties
(\ref{remop1})--(\ref{remop2}) which separates
$\Gamma_0=\Gamma_1\cup\dots\cup\Gamma_K$ from other boundary
components and from infinity. Denote by~$\Omega_{k\e}$ the bounded
domain with the boundary $\partial\Omega_{k\e}=S_k(\e)$. By
construction, $\Gamma_0\Subset\Omega_{k\e}$ and
\begin{equation}\label{rem0}
\int_{S_k(\e)}|\ue_k|^2\,dS\to0
\end{equation}
as $k\to\infty$ (see~(\ref{remo22ll}), (\ref{remo8}),
(\ref{est_i1})\,).

\medskip

{\sc Step 7.} Take arbitrary $\e>0$ and fix it (the precise value
of~$\e$ will be specified below). Now, for $t\in\Ti\cap[t',t'']$
and sufficiently large~$k$ (in particular, such that the claims of previous
Steps are fulfilled) consider the domain (see Fig.2)
$$
\Omega_{k}(t)=\bigl[\overline{\widetilde{V}}(t_1)\cup\widetilde
 W_k(t_1,t_2;t)\bigr]\cap B_k
\setminus\overline \Omega_{k\e}.
$$
Recall,
$B_k=\{x\in\R^3:|x|<R_k\}$
are the balls where the solutions $\ue_k\in
W^{1,2}(\Omega\cap B_k)$ are defined,
the sets ${V}(t_1)$ and $W_k(t_1,t_2;t)$ were defined previously (see text after
Corollary~\ref{regPhi-ax}), and for a set
$A\subset P_+$ we denote by $\widetilde A$ the set in $\R^3$ obtained
by rotation of $A$ around $O_z$-axis.

By construction (see Fig.2),
$\partial\Omega_{k}(t) =\widetilde S_k(t_1,t_2;t)\cup S_{R_k }\cup
S_k(\e)$. Integrating the equation
\begin{equation}\label{cle_lap}\Delta\Phi_k=\omega_k^2+\frac1{\nu_k}\div(\Phi_k\ue_k)-\frac{1}{\nu_k}\fe_k\cdot\ue_k\end{equation}
over the domain $\Omega_{k}(t)$, we have
$$
\int_{ \widetilde S_k(t_1,t_2;t)}\nabla\Phi_k\cdot{\bf
n}\,dS+\int_{ S_{R_k }\cup S_k(\e)}\nabla\Phi_k\cdot{\bf
n}\,dS=\int_{\Omega_{k}(t)}\omega_k^2\,dx-\frac1{\nu_k}\int_{\Omega_{k}(t)}\fe_k\cdot\ue_k\,dx
$$
$$
+\frac{1}{\nu_k}\int_{ \widetilde
S_k(t_1,t_2;t)}\Phi_k\ue_k\cdot{\bf n}\,dS+\frac{1}{\nu_k}\int_{
S_{R_k }\cup S_k(\e)}\Phi_k\ue_k\cdot{\bf n}\,dS
$$
\begin{equation}\label{cle_lap1}
=\int_{\Omega_{k}(t)}\omega_k^2\,dx-\frac1{\nu_k}\int_{\Omega_{k}(t)}\fe_k\cdot\ue_k\,dx-t{}\bar\F+\frac{1}{\nu_k}\int_{
S_{R_k }\cup S_k(\e)}\Phi_k\ue_k\cdot{\bf n}\,dS,
\end{equation}
where $\bar\F=\frac1\nu(\F_{K+1}+\dots+\F_N)$ (here we use the
identity~$\Phi_k\equiv-t$ on $\widetilde S_k(t_1,t_2;t)$\,). In
view of (\ref{lac-2-ax}), (\ref{lac0.5}), (\ref{remop2'}) and
(\ref{remop1})--(\ref{remop2})
 we can estimate
\begin{equation}\label{cle_lap2}
\int_{ \widetilde S_k(t_1,t_2;t)}|\nabla\Phi_k|\,dS \le
t\F+8\e+\frac1{\nu_k}\int_{\Omega_{k}(t)}\fe_k\cdot\ue_k\,dx-\int_{\Omega_{k}(t)}\omega_k^2\,dx
\end{equation} with $\F=|\bar\F|$. By definition, $\frac{1}{\nu_k}\|\fe_k\|_{L^{6/5}(\Omega)}=
\frac{\nu_k}{\nu^2}\|\fe\|_{L^{6/5}(\Omega)}\to 0$ as
$k\to\infty$. Therefore, using the uniform
estimate~$\|\ue_k\|_{L^6(\Omega_k)}\le\const$, we have
$$
\Big|\frac1{\nu_k}\int_{\Omega_{k}(t)}\fe_k\cdot\ue_k\,dx\Big|<\varepsilon
$$
for sufficiently large $k$. Then (\ref{cle_lap2}) yields
\begin{equation}\label{cle_lap3}
\int_{ \widetilde S_k(t_1,t_2;t)}|\nabla\Phi_k|\,dS <
t\F+9\e-\int_{\Omega_{k}(t)}\omega_k^2\,dx.
\end{equation}
Choosing $\e$ sufficiently small so that~$9\e<\e_t$
(see~(\ref{omeg})) and a sufficiently large $k$ (with ${\bar
h'_k}<\delta_t$), we deduce from~(\ref{omeg}) that
$$9\e-\int_{\Omega_{k}(t)}\omega_k^2\,dx\le 9\e- \int_{\widetilde
 V(t_1)\cap B_k \setminus\overline \Omega_{\bar
h'_k}}\omega_k^2\,dx<0.$$ Estimate~(\ref{mec}) is proved.  \qed

\bigskip

We need the following technical fact from the one-dimensional real analysis.

\begin{lem}
\label{real-an2}{\sl Let $f:\Si\to\R$ be a positive decreasing
function defined on a~measurable set~$\Si\subset(0,\delta)$ with
$\meas[(0,\delta)\setminus\Si]=0$. Then
\begin{equation}\label{isoper2}
\sup\limits_{t_1,t_2\in\Si}\frac{[f(t_2)]^{\frac43}(t_2-t_1)}{(t_2+t_1)(f(t_1)-f(t_2))}=\infty.
\end{equation}
}
\end{lem}

The proof of this fact is  elementary, see Appendix.\\

Below we will use the key estimate~(\ref{mec}) to prove some geometrical
relations that contradict Lemma~\ref{real-an2}.
\\

For $t\in\Ti$ denote by $U(t)$ the union of bounded connected
components (tori) of the set $\R^3\setminus
\bigl(\cup_{j=M+1}^N\widetilde A_j(t)\bigr)$. By construction,
$U(t_2)\Subset U(t_1)$ for $t_1<t_2$.

\begin{lem}
\label{ax-lkr12}{\sl For  any $t_1,t_2\in\Ti$ with $t_1<t_2$ the
estimate
\begin{equation}\label{isoper1}
\meas U(t_2)^{\frac43}\le C\frac{t_2+t_1}{t_2-t_1}\bigl[\meas
U(t_1)-\meas U(t_2)\bigr]
\end{equation}
holds with the constant $C$ independent of  $t_1,t_2$. }
\end{lem}

\pr Fix $t_1,t_2\in\Ti$ with $t_1<t_2$. Take a~pair $t',t''\in\Ti$
such that $t_1<t'<t''<t_2$. For $k\ge k_*(t_1,t_2,t',t'')$ (see
Lemma~\ref{ax-lkr11}) put
$$
E_{k}=\bigcup\limits_{t\in [t',t'']} \widetilde S_k(t_1,t_2;t).
$$
 By the Coarea formula (see, e.g, \cite{Maly}), for any integrable
 function $g:E_{k}\to\R$ the equality
\begin{equation}\label{Coarea_Phi}\int\limits_{E_k}g|\nabla\Phi_k|\,dx=
\int\limits_{t'}^{t''}\int_{\widetilde
S_k(t_1,t_2;t)}g(x)\,d\Ha^2(x)\,dt
\end{equation}
holds. In particular, taking $g=|\nabla\Phi_k|$ and
using~(\ref{mec}), we obtain
\begin{equation}\label{Coarea_Phi2}
\begin{array}{lcr}
\int\limits_{E_k}|\nabla\Phi_k|^2\,dx=
\int\limits_{t'}^{t''}\int_{ \widetilde
S_{k}(t_1,t_2;t)}|\nabla\Phi_k|(x)\,d\Ha^2(x)\,dt
\\
\displaystyle
\le
\int\limits_{t'}^{t''}\F t\,dt=
\frac\F2\bigl((t'')^2-(t')^2\bigr).
\end{array}
\end{equation}
Now, taking  $g=1$ in (\ref{Coarea_Phi})  and using the H\"older
inequality we have
\begin{equation}\label{Coarea_Phi3}
\begin{array}{lcr}
\displaystyle \int\limits_{t'}^{t''}\Ha^2\bigl( \widetilde
S_{k}(t_1,t_2;t)\bigr)\,dt= \int\limits_{E_k}|\nabla\Phi_k|\,dx
\\
\displaystyle \le
\biggl(\int\limits_{E_k}|\nabla\Phi_k|^2\,dx\biggr)^{\frac12}
\bigl(\meas
(E_k)\bigr)^{\frac12}\le\sqrt{\frac\F2\bigl((t'')^2-(t')^2\bigr)\meas
(E_k)}.
\end{array}
\end{equation}
 By construction (see property (iii) in the beginning of Subsection~\ref{EPcontr-axx}), each of the sets $A_j(t_1)$ and $A_j(t_2)$ is
 a~smooth cycle surrounding the component~$\breve\Gamma_j$ for $j=M+1,\dots,N$,
 moreover, the cycle $A_j(t_1)$ lies
 in the unbounded connected component of the open
 set~$P_+\setminus A_j(t_2)$. Furthermore,
 for almost all $t\in [t',t'']$ the set $S_{k}(t_1,t_2;t)$ is a finite union of smooth cycles in $P_+$ and
$S_{k}(t_1,t_2;t)$ separates $A_j(t_1)$ from $A_j(t_2)$ for all
$j=M+1,\dots,N$. In particular, the set $U(t_2)$ is contained in
the union of bounded connected components of~$\R^3\setminus
\widetilde S_{k}(t_1,t_2;t)$. Then by the isoperimetric inequality
(see,  e.g., \cite{evans}),  $\Ha^2\bigl(\widetilde
S_{k}(t_1,t_2;t)\bigr)\ge C_*\bigl(\meas  U(t_2)\bigr)^{\frac23}$
for $t\in[t',t'']$. Therefore, (\ref{Coarea_Phi3}) implies
\begin{equation}\label{Coarea_Phi3'}
\bigl(\meas  U(t_2)\bigr)^{\frac43}(t''-t')^2\le
C\bigl((t'')^2-(t')^2\bigr)\meas (E_k).
\end{equation}
On the other hand, by definition, $S_{k}(t_1,t_2;t)\subset
V(t_2)\setminus V(t_1)$. Consequently, $\widetilde
S_{k}(t_1,t_2;t)\subset U(t_1)\setminus U(t_2)$ for all
$t\in[t',t'']$, hence from (\ref{Coarea_Phi3'}) we get
\begin{equation}\label{Coarea_Phi4}
\begin{array}{lcr}
\bigl(\meas  U(t_2)\bigr)^{\frac43}\le
C\dfrac{t''+t'}{t''-t'}\meas\bigl( U(t_1)\setminus
U(t_2)\bigr)
\\
\displaystyle
=C\dfrac{t''+t'}{t''-t'}\bigl[\meas U(t_1)-\meas
U(t_2)\bigr].
\end{array}
\end{equation}
The last estimate is valid for every pair $t'', t'\in(t_1,t_2)$.
Taking a limit as $t''\to t_2, t'\to t_1$, we obtain the required
estimate~(\ref{isoper1}).
$\qed$\\

The last estimate leads us to the~main result of this subsection.

\begin{lem}
\label{lem_Leray_fc} {\sl Assume that  $\Omega\subset\R^3$ is an
exterior  axially symmetric domain of type \eqref{Omega} with
$C^2$-smooth boundary $\partial\Omega$,  and~$\fe\in
W^{1,2}(\Omega)$, ${\bf a}\in
W^{3/2,2}(\partial\Omega)$ are axially symmetric. Then assumptions
(E-NS) and \eqref{as-prev1} lead to a contradiction.}
\end{lem}

\pr By construction, $U(t_1)\supset U(t_2)$ for $t_1,t_2\in\Ti$,
$t_1<t_2$. Thus the just obtained estimate~(\ref{isoper1})
contradicts Lemma~\ref{real-an2}. This
contradiction finishes the proof of Lemma~\ref{lem_Leray_fc}.
 \qed

\subsection{The case $0<\widehat
p_N=\esssup\limits_{x\in\bar\Omega}\Phi(x)$.}\label{Euler-contr-2}

Suppose now that (\ref{as1-axxx}) holds, i.e., the maximum of
$\Phi$ is attained on the boundary component~$\Gamma_N$ which does
not intersect the symmetry axis. Repeating the arguments from the
first part of the previous Subsection~\ref{EPcontr-axx}, we can construct
a~$C^1$~-smooth cycle $A_N\subset\D$ such that
$\psi|_{A_N}=\const$, $0<\Phi(A_N)<\widehat p_N$ and $\breve
\Gamma_N=P_+\cap\Gamma_N$ lies in the bounded connected component
of the set~$P_+\setminus A_N$. Denote this component by~$\D_b$.  The cycle $A_N$ separates
$\breve \Gamma_N$ from infinity and from the singularity
line~$O_z$. This means, that we can reduce our case to the situation
with bounded domain~$\D_b$ (surrounded by~$A_N$), considered in~\cite{kpr_a_ann}.
Describe some details of this reduction.

Let $$\D_b\cap \Gamma_j=\emptyset,\qquad j=1,\dots M_1-1,$$
$$\D_b\supset\breve \Gamma_j,\qquad j=M_1,\dots,N$$
(the case $M_1=N$ is not excluded). Making a renumeration (if necessary), we may assume without loss of generality that
$$\Phi(\breve \Gamma_j)<\widehat p_N,\qquad j=M_1,\dots,M_2,$$
$$\Phi(\breve \Gamma_j)=\widehat p_N,\qquad j=M_2+1,\dots,N$$
(the case $M_2=M_1-1$, i.e., when $\Phi$ attains maximum value at every boundary component
inside the domain~$\D_b$, is not excluded). Apply Kronrod results from Subsection~\ref{Kronrod-s}
to the restriction~$\psi|_{\D_b}$ of the stream function~$\psi$ to the domain~$\D_b$ (this is possible
because of Remark~\ref{kmpRem1.2}). Let $T_\psi$ denote
the corresponding Kronrod tree for this restriction.
Denote by $B_0$ the element of~$T_\psi$ with $B_0\supset A_N$.
Similarly, denote by~$B_j$, $j=M_1,\dots,N$, the elements of~$T_\psi$ with $B_j\supset\breve \Gamma_j$.
Adding a constant to the pressure~$p(x)$, we can assume from this moment that
$$\Phi(B_0)=\Phi(A_N)<0,$$
$$\Phi(B_j)=\Phi(\breve\Gamma_j)<0\qquad j=M_1,\dots,M_2,$$
$$\Phi(B_j)=\Phi(\breve\Gamma_j)=0\qquad j=M_2+1,\dots,N.$$ Now in order to
receive the required contradiction, one need to~consider the behavior of $\Phi$ on the Kronrod arcs $[B_j,B_N]$ and to repeat word
by word the corresponding arguments
of Subsection~4.2.2 in~\cite{kpr_a_ann} starting from Lemma~4.7. The only modifications are as follows:
now our
sets~$B_0$ and $B_{M_1},\dots, B_{M_2}$ play the role of the sets~$C_{M'}$
and $C_{M'+1},\dots,C_{M}$ from~\cite[Subsection~4.2.2]{kpr_a_ann} respectively. Also,
the domain $\D_b\cap\D$ from the present case plays the role of the domain~$D_{r_*}$
from~\cite[Subsection~4.2.2]{kpr_a_ann}.

\subsection{The case $\max\limits_{j=1,\dots,N}\widehat
p_j<\esssup\limits_{x\in\bar\Omega}\Phi(x)>0$.}\label{Euler-contr-3}

Suppose  that (\ref{as-prev-id-ax}) holds, i.e., the  maximum
of $\Phi$ is not zero and it is not attained on $\partial\Omega$
(the case $\esssup\limits_{x\in\Omega}\Phi(x)=+\infty$ is not
excluded). We start from  the following simple fact.

\begin{lem}
\label{lem_PR} {\sl Under assumptions~(\ref{as-prev-id-ax}) there exists a compact connected set~$F\subset \D$ such that
$\diam F>0$, $\psi|_F\equiv\const$, and
$$0<\Phi(F)>\max\limits_{j=1,\dots,N}\widehat p_j.$$}
\end{lem}

{\bf The proof} of this Lemma is quite similar to the proof of \cite[Lemma 3.10]{kpr_a_ann} which was done for the case of bounded plane domain. But since the present situation has some specific differences, for reader's convenience we reproduce the proof with the corresponding modifications. Denote $\sigma=\max\limits_{j=1,\dots,N}\widehat p_j$. By the assumptions,
$\Phi(x)\le\sigma$ for every $x\in P_+\cap\partial\D\setminus A_\ve$ and there is a set of a positive plane measure $E\subset \D\setminus A_\ve$ such that $\Phi(x)>\sigma$ at each~$x\in E$.
In virtue of Theorem~\ref{kmpTh2.1} (iii), there exists a straight-line segment $I=[x_0,y_0]\subset\overline\D\cap P_+$ with
$I\cap A_\ve=\emptyset$, $x_0\in\partial\D$, $y_0\in E$, such that $\Phi|_I$ is a continuous function.
By construction, $\Phi(x_0)\le\sigma$, $\Phi(y_0)\ge\sigma+\delta_0$ with some $\delta_0>0$.
Take a subinterval $I_1=[x_1,y_0]\subset\D$ such that $\Phi(x_1)=\sigma+\frac12\delta_0$ and $\Phi(x)\ge\sigma+\frac12\delta_0$ for each
$x\in[x_1,y_0]$. Then by Bernoulli's Law (see Theorem~\ref{kmpTh2.2}) $\psi\ne \const$ on $I_1$. Take a~closed rectangle $Q\subset\D$ such that $I_1\subset Q$. By Theorem~\ref{kmpTh1.1}~(iii) applied to~$\psi|_{Q}$ we can take $x\in I_1$ such that the preimage $\{y\in Q:\psi(y)=\psi(x)\}$ consists of a finite union of $C^1$-curves. Denote by $F$ the~curve containing~$x$. Then by construction
$\Phi(F)\ge\sigma+\frac12\delta_0$, $F\subset Q\subset\D$ and $\diam F>0$. $\qed$

\medskip
Fix a compact set $F$ from Lemma~\ref{lem_PR}.
Using the arguments from the first part of
Subsection~\ref{EPcontr-axx}, we can construct a~$C^1$~-smooth
cycle $A_F\subset\D$ such that $\psi|_{A_F}=\const$,
$0<\Phi(A_F)<\Phi(F)$ and $F$ lies in the bounded connected
component of the set~$P_+\setminus A_F$, denote this component by~$\D_b$ (in this procedure
the set~$F$ plays the role of the set~$\breve\Gamma_j$ from the beginning of Subsection~\ref{EPcontr-axx},
where the cycles~$A_j(t)$ were constructed with properties~(i)--(iii)).
 As before, $A_F$
separates $F$ from infinity and from the singularity line~$O_z$.
This means, that we can reduce our case to the situation with bounded
domain $D_b$ (surrounded by~$A_F$), considered in~\cite{kpr_a_ann}. Describe some details of this reduction.

Let $$\D_b\cap \Gamma_j=\emptyset,\qquad j=1,\dots M_1-1,$$
$$\D_b\supset\breve \Gamma_j,\qquad j=M_1,\dots,N$$
(the case $M_1=N+1$, i.e., when $\D_b\cap\partial\D=\emptyset$, is not excluded).
Apply the Kronrod results from Subsection~\ref{Kronrod-s}
to the restriction~$\psi|_{\D_b}$ of the stream function~$\psi$ to the domain~$\D_b$.
Let $T_\psi$ denote
the corresponding Kronrod tree for this restriction.
Denote by $B_0,B_F$ the element of~$T_\psi$ with $B_0\supset A_F$, $B_F\supset F$.
Similarly, denote by~$B_j$, $j=M_1,\dots,N$, the elements of~$T_\psi$ with $B_j\supset\breve \Gamma_j$.
Adding a constant to the pressure~$p(x)$, we can assume from this moment that
$$\Phi(B_0)=\Phi(A_F)<0,$$
$$\Phi(B_j)=\Phi(\breve\Gamma_j)<0\qquad j=M_1,\dots,N,$$
$$\Phi(B_F)=\Phi(F)=0.$$
 Now in order to
receive the required contradiction, one need to consider the behavior of $\Phi$ on the Kronrod arcs $[B_j,B_F]$ and to repeat almost word
by word the corresponding arguments
of Subsection~3.3.2 in~\cite{kpr_a_ann} after Lemma~3.10. The only modifications are as follows:
now our
sets~$B_0$, $B_{M_1},\dots, B_{N}$, and $B_F$ play the role of the sets~$B_0,\dots,B_N$ and $F$ from~\cite[Subsection~3.3.2]{kpr_a_ann} respectively. Also,
the domain $\D_b\cap\D$ from the present case plays the role of the domain~$\Omega$
from~\cite[Subsection~3.3.2]{kpr_a_ann}, and
on the final stage we have to integrate  identity (\ref{cle_lap**}) of the present paper over the three--dimensional domains $\Omega_{ik}(t)$ with $\partial \Omega_{ik}(t)=\widetilde S_{ik}(t)$. $\qed$

\bigskip We can summarize the results of Subsections~\ref{Euler-contr-2}--\ref{Euler-contr-3} in the following statement.

\begin{lem}
\label{lem_Leray_fc_3} {\sl Assume that $\Omega\subset\R^3$ is an
exterior  axially symmetric domain of type~\eqref{Omega} with
$C^2$-smooth boundary $\partial\Omega$ and ~$\fe\in
W^{1,2}(\Omega)$, ${\bf a}\in
W^{3/2,2}(\partial\Omega)$ are axially symmetric. Let (E-NS) be fulfilled. Then each assumptions
(\ref{as1-axxx}), (\ref{as-prev-id-ax}) lead to a contradiction.
}
\end{lem}

\bigskip {\bf Proof of Theorem \ref{kmpTh4.2}.} Let the hypotheses
of Theorem \ref{kmpTh4.2} be satisfied. Suppose that its assertion
fails. Then, by Lemma~\ref{lem_Leray_symm}, there exist $ \ve, p$
and a sequence $(\ue_k,p_k)$ satisfying (E-NS), and by Lemmas~\ref{lem_Leray_fc}, \ref{lem_Leray_fc_3}
these
assumptions lead to a contradiction.    \qed

\begin{lemr}
\label{rem_rot}{\rm Let in Lemma~\ref{lem_Leray_symm}  the data
$\fe$ and~$\bf a$ be  axially symmetric with no swirl. If the
corresponding assertion of Theorem~\ref{kmpTh4.2} fails, then that
conditions~(E-NS) are satisfied with $\ue_k$ axially symmetric
with no swirl as well (see Theorem~\ref{Th_Ex_b}). But since we
have proved that assumptions~(E-NS) lead to a contradiction in
the~more general case (with possible swirl), we get also the validity of
second assertions of Theorem~\ref{kmpTh4.2}. }
\end{lemr}

\section{The
case $\ue_0\ne0$} \label{contrad-u0} \setcounter{theo}{0}
\setcounter{lem}{0} \setcounter{lemr}{0}\setcounter{equation}{0}
\setcounter{lemA}{0}

In the proof of Theorem~\ref{kmpTh4.2} we have assumed that the
assigned value of the velocity at infinity is zero: $\ue_0=0$.  If
$\ue_0\ne0$, we can use the same arguments with some
modifications. First, we need some additional facts on Euler
equations.

\subsection{Some identities for solutions to the Euler system}
\label{E-u0} Let the conditions~(E) be fulfilled, i.e., axially
symmetric functions $(\ve,p)$ satisfy to Euler
equations~(\ref{2.1}) and
$$\ve\in L^6(\R^3),\quad p\in L^3(\R^3),$$
$$\nabla\ve\in L^2(\R^3),\quad \nabla p\in L^{3/2}(\R^3), \quad \nabla^2p\in L^1(\R^3)$$
(these properties were discussed in Section~\ref{eueq}).

For a $C^1$-cycle $\mathcal S\subset P_+$ (i.e., $\mathcal S$ is a
curve homeomorphic to the circle) denote by $\Omega_{\mathcal S}$
the bounded domain in $\R^3$ such that $\partial\Omega_{\mathcal
S}=\widetilde{\mathcal S}$, where, recall, $\widetilde{\mathcal
S}$ means the surface obtained by rotation of the curve~$\mathcal
S$ around the symmetry axis.

\medskip
\begin{lem}
\label{lemu-Eu-01} {\sl If conditions (E) are satisfied, then for
any $C^1$-cycle $\mathcal S\subset P_+$ with $\psi|_{\mathcal
S}\equiv\const$ we have
$$\int_{\Omega_{\mathcal S}}\ve\cdot\partial_z\ve\,dx=0.$$}
\end{lem}

\pr By Bernoulli Law (see Theorem~\ref{kmpTh2.2}\,) we have
$\Phi\equiv\const$ on $\mathcal S$, therefore,
$$\int_{\Omega_{\mathcal S}}\partial_z\Phi\,dx=\int_{\Omega_{\mathcal S}}\bigl[\partial_zp+\ve\cdot\partial_z\ve\bigr]\,dx=0.$$
Thus, to finish the proof of the Lemma, we need to check the
equality
\begin{equation}\label{p-uo}
\int_{\Omega_{\mathcal S}}\partial_zp\,\,dx=0.
\end{equation}
Denote by $\D_{\mathcal S}$ the open bounded domain in the
half-plane $P_+$ such that $\partial\D_{\mathcal S}=\mathcal S$.
Of course, $\Omega_{\mathcal S}=\widetilde\D_{\mathcal S}$. Then
the required assertion~(\ref{p-uo}) could be rewritten in the
following form
\begin{equation}\label{p-u01}
\int_{\D_{\mathcal S}}r\partial_zp\,\,drdz=\int_{\D_{\mathcal
S}}r\bigl[v_r\partial_rv_z+ v_z\partial_zv_z\bigr]\,drdz=0,
\end{equation}
where we have used the Euler equation~(\ref{2.1'}${}_1$)
for~$\partial_zp$. Since the gradient of the stream
function~$\psi$ satisfies $\nabla\psi\equiv(-rv_z,rv_r)$,
 we could
rewrite~(\ref{p-u01}), using the Coarea formula, in the following
equivalent form
\begin{equation}\label{p2-u01}
\int_\R dt\int_{\psi^{-1}(t)\cap\D_{S}}\nabla v_z\cdot{\bf l}\,ds
=\int_\R dt\int_{\psi^{-1}(t)\cap\D_{S}}\frac{\partial
v_z}{\partial s}\,ds=0,
\end{equation}
where ${\bf l}=\dfrac1{|\nabla\psi|}(rv_r,rv_z)$ is the tangent
vector to the stream lines~$\psi^{-1}(t)$. The last equality~in
(\ref{p2-u01}) is evident because almost all level lines of~$\psi$
in $\D_{\mathcal S}$ are $C^1$-curves homeomorphic to the circle
(see the Morse-Sard Theorem~\ref{kmpTh1.1}\,(iii)\,). The Lemma is
proved. \qed

\medskip

We need also the following simple technical fact.
\begin{lem}
\label{lemu-Eu2-01} {\sl If $\ue=(u_\theta,u_r,u_z)$ is
$C^1$-smooth axially-symmetric vector field in~$\Omega$ with
$\div\ue\equiv0$, then for any Lipschitz curve $\mathcal S\subset
P_+\cap\Omega$ such that $\widetilde{\mathcal S}$ is a compact
closed Lipschitz surface the identity
$$\int_{\widetilde{\mathcal S}}{\bf n}\cdot\partial_z\ue\,dS=0$$
holds.}
\end{lem}

\pr There are two possibilities: $\mathcal S$ is homeomorphic to
the circle, or~$\mathcal S$ is homeomorphic to the straight
segments with endpoints on symmetry axis. Consider the last case
(the first case could be done analogously). By identity
$\div\ue=0$, we have
$$r\partial_zu_z\equiv-\partial_r(ru_r).$$
Then by direct calculation we have
$$\int_{\widetilde{\mathcal S}}{\bf n}\cdot\partial_z\ue\,dS=\int_{\mathcal S}r{\bf n}\cdot\partial_z\ue\,ds=\int_{\mathcal S}r(-dz,dr)\cdot\partial_z\ue$$
$$=
\int_{\mathcal S} \bigl[-\frac{\partial(ru_r)}{\partial
z}dz+\frac{\partial(ru_z)}{\partial z}dr\bigr]= \int_{\mathcal S}
\bigl[-\frac{\partial(ru_r)}{\partial
z}dz-\frac{\partial(ru_r)}{\partial r}dr\bigr]=\int_{\mathcal
S}d(ru_r)=0.$$ The Lemma is proved. \qed

\subsection{The existence theorem}
\label{contrad_ue0} Let the hypotheses of Theorem \ref{kmpTh4.2}
be satisfied. Suppose that its assertion fails. Then, as in the
first part of Section~\ref{poet}, we obtain the sequence
 of solutions~$(\tilde\ue_k,p_k)$ to problems~(\ref{NSk})
with $\tilde\ue_k=\w_k+\frac{\nu_k}{\nu}(\h+\ue_0)$, where ${\bf
A}={\bf a}-{\bf u}_0$ on $\partial\Omega$,  ${\bf
A}={\boldsymbol{\sigma}}$  for sufficiently large $|x|$,
$\|\w_k\|_{H(\Omega_k)}\equiv1$, and $\nu_k\to0$ as $k\to\infty$.
Take $\ue_k=\w_k+\frac{\nu_k}{\nu}\h$ and note that $\ue_k$ is a
solution to the Navier-Stokes system
\begin{equation}\label{NSk-u0}
\left\{\begin{array}{rcl}-\nu_k\Delta{\bf u}_k +\big({\bf
u}_k\cdot\nabla\big){\bf u}_k+\nabla p_k & = &
\fe_k+\tilde\fe_k\qquad \hbox{\rm in }\;\;\Omega_k,
\\[4pt]
\div{\bf u}_k & = & 0 \;\qquad \hbox{\rm in }\;\;\Omega_k,
\\[4pt]
{\bf u}_k &  = & {\bf a}_k
 \quad\  \hbox{\rm on }\;\;\partial\Omega_k,
\end{array}\right.
\end{equation}
 with $\fe_k=\frac{\nu_k^2}{\nu^2}\,{\bf f}$, \ ${\bf
a}_k=\frac{\nu_k}\nu\,{\bf A}$, and
$$\tilde\fe_k=-\frac{\nu_k}{\nu}(\ue_0\cdot \nabla)\ue_k=\nu_k\alpha\partial_z\ue_k,$$
where $\alpha\in\R$ is a constant. Since by H\"older inequality
$\|\fe\|_{L^{3/2}(\Omega_R)}\le
\|\fe\|_{L^{2}(\Omega_R)}\sqrt{R}$, and, consequently, $\|
\tilde\fe_k\|_{L^{2}(\Omega_R)}\le C\sqrt{R}$, we conclude that
Corollary~\ref{LemPr2} implies the uniform estimate
\begin{equation}\label{pru01-u0}
\|\nabla p_k\|_{L^{3/2}(\Omega_R)}\leq C\sqrt{R}\qquad\mbox{ for
all }R\in[R_0,R_k],
\end{equation}
i.e., the estimate~(\ref{epgr34}) holds. The another needed
estimate~$\|\ue_k\|_{L^6(\Omega_k)}\le C$ follows from the Sobolev
Imbedding Theorem. So, for the total head pressure $\Phi_k
=p_k+\frac12|\ue_k|^2$ we have
\begin{equation}\label{pru01-Phi-u0}
\|\nabla \Phi_k\|_{L^{3/2}(\Omega_R)}\leq C\sqrt{R}\qquad\mbox{
for all }R\in[R_0,R_k].
\end{equation}
By the same reasons as before, \begin{equation} \label{E-NS-u0}
\ue_k\rightharpoonup \ve\mbox{ \ in \
}W_\loc^{1,2}(\overline\Omega),\quad p_k\rightharpoonup p\mbox{ \
in \ }W_\loc^{1,3/2}(\overline\Omega),
\end{equation}
\begin{equation}\label{cont_u0}
\begin{array}{l}
\displaystyle \nu= \int _\Omega ({\bf v}\cdot\nabla){\bf
v}\cdot\h\,dx,
\end{array}
\end{equation}
where the limit functions~$(\ve,p)$ satisfy the Euler equation
with zero boundary conditions, i.e., condition (E) is fulfilled.
Our goal is to receive a contradiction.

We need to discuss only the case~(\ref{as-prev1}), since other
cases~(\ref{as1-axxx})--(\ref{as-prev-id-ax}) are reduced to the
consideration of bounded domains (see
Subsections~\ref{Euler-contr-2}--\ref{Euler-contr-3}\,) and zero
or nonzero  condition at infinity has no influence on the proof.

The arguments in Section~5 up to the 7-th Step of the proof of
Lemma~\ref{ax-lkr11} could be repeated almost word by word. But
the 7-th Step of the proof of Lemma~\ref{ax-lkr11} needs some
modifications. Recall that the main idea there was to use the
identity
\begin{equation}\label{i1_u0}
\begin{array}{lcr}
\int\limits_{\partial\Omega_k(t)}\nabla\Phi_k\cdot\n\,dS=\int\limits_{S_k(t_1,t_2;t)}\nabla\Phi_k\cdot\n\,dS+
\int\limits_{{ S_{R_k }\cup S_k(\e)}}\nabla\Phi_k\cdot\n\,dS\\
\\
= \int\limits_{\Omega_k(t)}\Delta\Phi_k\,dx.
\end{array}
\end{equation}
Since now
\begin{equation}
\label{i2_u0}
\nabla\Phi_k=-\nu_k\curl\,\ov_k+\ov_k\times\ue_k+\fe_k+\tilde\fe_k,
\end{equation}
\begin{equation}
\label{cle_laps_u0}
\Delta\Phi_k=\omega_k^2+\frac1{\nu_k}\div(\Phi_k\ue_k)-\frac{1}{\nu_k}{\bf
f}_k\cdot\ue_k- \frac{1}{\nu_k}\tilde\fe_k\cdot\ue_k,
\end{equation}
we need  to prove the smallness of the following integrals
generated by the additional term~$\tilde\fe_k$:
\begin{eqnarray}
\label{int1_u0} \int_{S_{R_k }\cup S_k(\e)}\tilde\fe_k\cdot{\bf
n}\,dS,\\
\label{int2_u0}\frac1{\nu_k}\int_{\Omega_{k}(t)}\tilde\fe_k\cdot\ue_k\,dx=
\alpha\int_{\Omega_{k}(t)}\ue_k\cdot\partial_z\ue_k\,dx.
\end{eqnarray}
The difficulty is, that the term
$\tilde\fe_k=\nu_k\alpha\partial_z\ue_k$ in the right-hand side
of~(\ref{NSk-u0}) is not "small enough" (it is of order~$O(\nu_k)$
only, not of order $O(\nu_k^2)$ as~$\fe_k$). However, this term
has very good symmetry properties, and by Lemma~\ref{lemu-Eu2-01}
we immediately obtain that the  integral \eqref{int1_u0} is
negligible.

\begin{lemA}
\label{lemu01} {\sl The identity
$$\int_{S_{R_k }\cup S_k(\e)}\tilde\fe_k\cdot{\bf
n}\,dS=0$$ holds.}
\end{lemA}

Let us estimate the integral which is in formula~(\ref{int2_u0}).
We need to use the limit solution of the Euler equations.

For $t\in\Ti$ denote $S(t)=\bigcup_{j=M+1}^N A_j(t)$ (recall, that
the set $\Ti$ and the~$C^1$-cycles $A_j(t)$ were defined in the
beginning of Section~5). Denote by $\Omega_{S(t)}$ the bounded
open set in $\R^3$ such that $\partial\Omega_{S(t)}=\widetilde
S(t)$. Further, put $\Omega'_{S(t)}=\Omega\cap\Omega_{S(t)}$.
Convergence~(\ref{E-NS-u0}) implies, in particular, that
\begin{equation}
\label{sc1} \ue_k\to \ve\mbox{ \ in \
}L^q_{\loc}(\overline\Omega)\quad\mbox{ for any }q\in[1,6).
\end{equation}
Therefore, by Lemma~\ref{lemu-Eu-01} we obtain immediately
\begin{lem}
\label{lemu-2i-01} {\sl For any $t\in\Ti$ the convergence
$$\int_{\Omega'_{S(t)}}\ue_k\cdot\partial_z\ue_k\,dx\to0$$
holds.}
\end{lem}
For $t\in\Ti$ denote $\mu(t)=\meas\Omega'_{S(t)}$. By
construction, the function $\mu(t)$ is strictly decreasing,
therefore it is continuous on~$\Ti$ except for at most countable
set. Removing the discontinuity points from~$\Ti$, we could assume
without loss of generality that the
function~$\mu:\Ti\to(0,+\infty)$ is continuous. By constructions
of Section~5, it is easy to see, that if $t_1,
\tau_1,t,\tau_2,t_2\in\Ti$ and $t_1<\tau_1<t<\tau_2<t_2$, then
\begin{equation}
\label{sc2} S_k(t_1,t_2;t)=S_k(\tau_1,\tau_2;t)\qquad \mbox{ and
}\qquad
\Omega'_{S(\tau_2)}\subset\Omega'_{S_k(t_1,t_2;t)}\subset\Omega'_{S(\tau_1)}
\end{equation}
for sufficiently large~$k$. Using these facts, the continuity
of~$\mu(t)$, the convergence~(\ref{sc1}) and
Lemma~\ref{lemu-2i-01}, we obtain
\begin{lem}
\label{lemu-2i1-01} {\sl For any $t_1,t,t_2\in\Ti$ with
$t_1<t<t_2$ the convergence
$$\int_{\Omega'_{S_k(t_1,t_2,t)}}\ue_k\cdot\partial_z\ue_k\,dx\to0$$ holds as $k\to\infty$.}
\end{lem}
Moreover, since the function $\mu(t)$ is {\it uniformly}
continuous on each compact subset of~$\Ti$, we have by the same
reasons that
\begin{lem}
\label{lemu-2i2-01} {\sl For any $\epsilon>0$ and $t_1,t_2\in\Ti$
with $t_1<t_2$ the convergence
$$\meas\bigl\{t\in(t_1,t_2):\biggl|\int_{\Omega'_{S_k(t_1,t_2,t)}}\ue_k\cdot\partial_z\ue_k\,dx
\biggr|>\epsilon\bigl\}\to0$$ holds as $k\to\infty$.}
\end{lem}
On the other hand, by constructions of Section~5
$$\partial\bigl(\Omega_k(t)\cup
\Omega'_{S_k(t_1,t_2,t)}\bigr)\subset S_k(\e)\cup
\partial\Omega_k$$ and
$$\|\ue_k\|_{L^2(S_k(\e)\cup \partial\Omega_k)}\to0$$
as $k\to\infty$ (see~(\ref{rem0}), (\ref{NSk-u0}${}_3$)\,).
Therefore, for any $t_1,t,t_2\in\Ti$ with $t_1<t<t_2$  we have
$$\int_{\Omega_k(t)\cup \Omega'_{S_k(t_1,t_2,t)}}\ue_k\cdot\partial_z\ue_k\,dx\to0$$
as $k\to\infty$. This fact together with Lemma~\ref{lemu-2i2-01}
implies the required assertion (the~smallness of the
integral~(\ref{int2_u0})\,):
\begin{lem}
\label{lemu-2i3-01} {\sl For any $\epsilon>0$ and $t_1,t_2\in\Ti$
with $t_1<t_2$ the convergence
$$\meas\bigl\{t\in(t_1,t_2):\biggl|\int_{\Omega_k(t)}\ue_k\cdot\partial_z\ue_k\,dx
\biggr|>\epsilon\bigl\}\to0$$ holds as $k\to\infty$.}
\end{lem}
Now repeating the arguments of Step~7 of the proof of
Lemma~\ref{ax-lkr11}, we obtained its assertion with the following
modification
\begin{lem}
\label{ax-lkr11-u0}{\sl There exists a constant~$\F>0$ such that
for  any $t_1,t_2\in\Ti$ with $t_1<t_2$ the convergence
$$\meas\bigl\{t\in(t_1,t_2):\int\limits_{\widetilde S_k(t_1,t_2;t)}|\nabla\Phi_k|\,dS\ge\F t
\bigl\}\to0$$ holds as $k\to\infty$. }
\end{lem}

It is easy to see that Lemma~\ref{ax-lkr11-u0} allow us to obtain
the assertion of Lemma~\ref{ax-lkr12} (without any changes) with
almost the same proof. We need only to define the corresponding
set $E_k$ as
$$E_k=\bigcup\limits_{t\in\Ti_k}\widetilde S_k(t_1,t_2;t),$$
where $$\Ti_k=\bigl\{t\in[t',t'']:\int\limits_{\widetilde
S_k(t_1,t_2;t)}|\nabla\Phi_k|\,dS\le\F t\bigr\},$$ and to use the
fact that, by Lemma~\ref{ax-lkr11-u0},
$\meas\bigl([t',t'']\setminus \Ti_k\bigr)\to0$ as $k\to\infty$.
Furthermore, Lemma~\ref{ax-lkr11-u0} together with
Lemma~\ref{real-an2} give us the required contradiction. This
contradiction finishes the proof of the Existence Theorem for the
case~$\ue_0\ne{\bf0}$.

\section{Appendix: Lemma from real analysis}
\setcounter{theo}{0} \setcounter{lem}{0}
\setcounter{lemr}{0}\setcounter{equation}{0}
\setcounter{lemA}{0}

We need the following elementary fact.

\begin{lem}
\label{real_an} {\sl Let $f:(0,\delta]\to\R$ be a positive
decreasing  function. Then
\begin{equation}\label{real1}
\esssup_{t\in(0,\delta]}\frac{[f(t)]^{\frac43}}{t|f'(t)|}=\infty.
\end{equation}
}
\end{lem}
\pr  Recall that by
the~Lebesgue theorem,  the derivative $f'(t)$ exists almost everywhere. Suppose
that the assertion~(\ref{real1}) fails. Then, taking into account
that~$f'(t)\le0$, we have
\begin{equation}\label{real2}
-f'(t)[f(t)]^{-\frac43}\ge \frac{C}{t}\quad\mbox{ for almost all
}t\in(0,\delta],
\end{equation}
with some positive constant $C$ independent of~$t$. Put
$g(t)=[f(t)]^{-\frac13}$. Then $g(t)$ is positive increasing
function on~$(0,\delta]$. By the Lebesgue theorem,
\begin{equation}\label{real-1}
g(t_2)-g(t_1)\ge\int_{t_1}^{t_2}g'(t)\,dt
\end{equation}
for any pair $t_1,t_2\in(0,\delta]$ with $t_1<t_2$. On the other
hand, (\ref{real2}) implies
\begin{equation}\label{real3}
g'(t)\ge \frac{C}{3t}\quad\mbox{ for almost all }t\in(0,\delta].
\end{equation}
The estimates~(\ref{real-1})--(\ref{real3}) contradict the
boundeness of $g$: $0<g(t)\le g(1)$ for all $t\in(0,\delta]$.
This contradiction finishes the proof. $\qed$

\medskip

{\bf Proof of Lemma~\ref{real-an2}}. We can extend the function~$f$ by one-sided continuity rule to
the whole interval $(0,\delta]$ and to apply
Lemma~\ref{real_an}.\qed

\section*{Acknowledgements}
$\quad\,\,$The authors are deeply indebted to V.V.~Pukhnachev for
valuable discussions.

The research of M. Korobkov was partially supported by the Russian
Foundation for Basic Research (project No.~12-01-00390-a).
M. Korobkov thanks also the Gruppo Nazionale per la Fisica Matematica of the Istituto Nazionale di Alta Matematica  for the financial support during his stays
in the Department of Mathematics and Physics of the Second University of Naples (Italy).

The research of K. Pileckas was funded by the Lithuanian-Swiss cooperation programme to
reduce economic and social disparities within the enlarged European Union under the project agreement No. CH-3-\v{S}MM-01/01.

\end{document}